%&bigtex
\input amstex
\input xy
\xyoption{all}
\SelectTips {cm}{}

\documentstyle{amsppt}
\pagewidth{35pc}

\topmatter
\title
Strict $\infty $-categories. Concrete Duality
\endtitle
\author
G.V. Kondratiev
\endauthor
\address
Department of Mathematics, Ottawa University
\endaddress
\email
gennadii\@\, hotmail.com
\endemail
\subjclass
CT Category Theory
\endsubjclass
\abstract
An elementary theory of strict $\infty $-categories with application to concrete duality is given. New examples of first and second order
concrete duality are presented.
\endabstract

\endtopmatter

\document

\head {} {\bf 1. Categories, functors, natural transformations, modifications}\endhead

\vskip 0.2cm
There are two kinds of weakness happenning to $\infty $-categories. One is changing all occurences of equality $=$ with a weaker equivalence
realtion $\sim $\, . The other one is a weak naturality condition. The first one is not proper and implies strict category theory. The second 
one is proper and gives a weak category theory. Below we use $\sim $ instead of $=$\, . It is not necessary but has an advantage to treat 
directly the classification problem (up to $\sim $).

\vskip 0.2cm
{\bf Definition 1.1.} 
\item{$\bullet $} {\bf $\infty$-precategory} is a (big) set $L$ endowed with
\roster
\item grading $L=\coprod\limits _{n\ge 0} L^n$
\item unary operations $d,c:\coprod\limits _{n\ge 1} L^n\to \coprod\limits _{n\ge 0} L^n$, $deg(d)=deg(c)=-1$, $dc=d^2$, $cd=c^2$
\item unary operation $e:\coprod\limits _{n\ge 0} L^n\to \coprod\limits _{n\ge 0} L^n$, $deg(e)=1$, $de=1$, $ce=1$
\item partial binary operations $\circ _k$, $k=1,2,...$, of degree $0$. $f\circ _k g$ is determined iff $d^kf=c^kg$
\endroster
\item{}such that each {\bf hom-set} $L(a,a'):=\{f\in L\ |\ \exists \, k\in {\Bbb N} \ d^kf=a, \ c^kf=a'\}$, $deg(a)=deg(a')$, inherits all properties (1)-(4).
\item{$\bullet $} $\forall \, a,a',a'' \in L^m$ there are maps $\mu _{a,a',a''} :\coprod\limits _{n\ge 0}L^n(a',a'')\times L^n(a,a')\to L(a,a'')$ 
such that if the bottom composite is determined then
$$\xymatrix{
\coprod\limits _{n\ge 0}L^n(a',a'')\times L^n(a,a') \ar[rr]^-{\mu _{a,a',a''}} && L(a,a'')\\
L^n(a',a'')\times L^n(a,a') \ar[u]^{i\times i} \ar[rr]_-{\circ _{n+1}} && L^n(a,a'') \ar[u]_{i}\\
}$$
$\mu _{a,a',a''}$ are called {\bf horizontal composites} on level $deg(a)$, all composites inside of $L(a,a')$ are {\bf vertical}. \hfill $\square $
\vskip 0.2cm
{\bf Remarks.} \item{$\bullet $} If $\alpha ^n, \beta ^n \in L^n, n>0$, such that $d\alpha ^n\ne d\beta ^n$ or $c\alpha ^n\ne c\beta ^n$, then 
$L(\alpha ^n,\beta ^n)=\emptyset$ (because of $d^2=dc, \, c^2=cd$). So, $\mu _{a,a',a''}$ can be empty map $\emptyset :\emptyset \to \emptyset$. 
\item{$\bullet $} It is convenient to use a letter with appropriate superscript, like $x^m, \alpha ^k$, etc., as an element (or sometimes as a variable) 
with domain $L^m, L^k$, etc. respectively (or with domain $L^m(a,b), L^k(x,y)$, etc.) \, Also, grading can be taken over all $\Bbb Z$ 
under assumption $L^{-m}:=\emptyset , \ m>0$. 
\item{$\bullet $} Call elements $a\in L^0$ of degree $0$ by {\bf objects} of $L$, elements $f^n\in L^n(a,a'), \, a,a'\in L^0$, by {\bf arrows of degree $n+1$ from $a$ to $a'$}. 
\item{$\bullet $} Denote {\bf horizontal composites} by $*$\, , and extend it over arrows {\bf of different degrees} 
by the rule $*:L(b,c)\times L(a,b)\to L(a,c):(g^n,f^m)\mapsto \mu _{a,b,c}(e^{max(m,n)-n}g^n,e^{max(m,n)-m}f^m)=:g^n*f^m$
($f^m\in L^m(a,b), g^n\in L^n(b,c)$).  \hfill $\square $

{\bf Definition 1.2.} For $a,b\in L^n$ $a\sim b$ iff \ $\exists \
\xymatrix{a \ar@/^/[r]^f & b \ar@/^/[l]^g}$ such that $e(a)\sim g\circ _1f$ \ $f\circ _1g \sim e(b)$
(it means that $\exists $ $f\in L^0(a,b)$, $g\in L^0(b,a)$ and two certain infinite sequences of arrows of higher order,
one in $L(a,a)$ and the other in $L(b,b)$).        \hfill $\square $

\vskip 0.1cm
$\sim $ \, is reflexive and symmetric, but may be not transitive.

\proclaim{\bf {Lemma 1.1}} If $L$ is an $\infty $-precategory such that 
\item{}$\circ _1$ is weakly associative $f\circ _1(g\circ _1h)\sim (f\circ _1g)\circ _1h$ (for composable arrows),
\item{}$\circ _1$ satisfies weak unit law  
\ $\forall f\in \coprod\limits _{n\ge 1}L^n$ $\cases f\circ _1 edf\sim f & \\
ecf\circ _1 f\sim f & \endcases $\hskip -0.3cm,
\item{}$\sim $ \, is compatible with $\circ _1$ \, $(f\sim g)\, \& \, (h\sim k)\Rightarrow (f\circ _1h)\sim (g\circ _1k)$ (for composable arrows), 
\item{}$\sim $ \, is transitive in higher orders, i.e. $\exists \, m>0$ such that $\sim $ is transitive for $\coprod\limits _{n\ge m}L^n$, 
\item{}then $\sim $ \, is transitive in all orders.
\endproclaim
\demo{Proof} Let $\vcenter{\xymatrix{a \ar@/^1ex/[r]^-{f} & b \ar@/^1ex/[l]^-{g} \ar@/^1ex/[r]^-{f'} & c \ar@/^1ex/[l]^-{g'}}}$ be given 
equivalences, i.e. $ea\sim g\circ _1f$, \, $eb\sim f\circ _1g$, \, $eb\sim g'\circ _1f'$, \, $ec\sim f'\circ _1g'$. Then 
$\vcenter{\xymatrix{a \ar@/^1ex/[r]^-{f'\circ _1f} & c \ar@/^1ex/[l]^-{g\circ _1g'}}}$ is the required equivalence since 
$ea\sim g\circ _1f\sim g\circ _1(eb\circ _1f)\sim g\circ _1((g'\circ _1f')\circ _1 f)\sim (g\circ _1g')\circ _1(f'\circ _1f)$ and similarly
$ec\sim (f'\circ _1f)\circ _1(g\circ _1g')$. \hfill $\square $
\enddemo
{\bf Remarks.} 
\item{$\bullet $} Transitivity in higher orders trivially holds for $n$-categories (starting from level $n$). For proper 
$\infty $-categories it is better to make assumption '\, $\sim $ \, is transitive in all orders' from the beginning.
\item{$\bullet $} This lemma shows that although transitivity of \, $\sim $ \, is not automatic for $\infty $-precategories, it is very 
consistent with (weak) associativity, unit law, and compatibility of \, $\sim $ \, with composites.   \vskip 0.0cm  \hfill    $\square $
\vskip 0.1cm
{\bf Definition 1.3.} $\infty $-precategory $L$ with relation \, $\sim $ \, as above is called an {\bf $\infty $-category} iff
\item{$\bullet $} \, $\sim $ \, is transitive \, $\alpha \sim \beta \sim \gamma \Rightarrow \alpha \sim \gamma $,
\item{$\bullet $} \, $\sim $ \, is compatible with all composites \, $(f\sim g)\, \& \, (h\sim k)\Rightarrow (f\circ _nh)\sim (g\circ _nk)$ 
(when they are defined),
\item{$\bullet $} horizontal composites preserve properties (1)-(2) and weakly preserve properties (3)-(4) of $\infty $-precategories:
\roster
\item grading $deg_{L(a,a'')}(\mu _{a,a',a''}(f,g))=deg_{L(a',a'')}(f)=deg_{L(a,a')}(g)$
\item $\mu _{a,a',a''}(df,dg)=d\mu _{a,a',a''}(f,g)$, \ $\mu _{a,a',a''}(cf,cg)=c\mu _{a,a',a''}(f,g)$
\item $\mu _{a,a',a''}(ef,eg)\sim \, e\mu _{a,a',a''}(f,g)$
\item $\mu _{a,a',a''}(f\circ _k f',g\circ _k g')\sim \mu _{a,a',a''}(f,g)\circ _k \mu _{a,a',a''}(f',g')$ \ ("interchange law")
\endroster
\item{$\bullet $} each $\circ _k, k\in \Bbb N$, is {\bf associative} $(f\circ _k g)\circ _k h\sim f\circ _k(g\circ _k h)$ 
(for composable elements),
\item{$\bullet $} {\bf unit law} holds $e^{k}c^kf\circ _k f\sim f$, \ $f\circ _k e^kd^kf\sim f$ (when all operations are defined).   \hfill $\square $

\vskip 0.2cm
{\bf Remarks.} 
\item{$\bullet $} By lemma 1.1, \, for $n$-categories transitivity condition on \, $\sim $ \, follows from the others.
\item{$\bullet $} Hom-sets in $\infty $-category $L$ are $\infty $-categories themselves, and horizontal composites \,  
$*:L(b,c)\times L(a,b)\to L(a,c)$, are $\infty $-functors.    
\item{$\bullet $} Since strict functors preserve equivalences $\sim $ for categories in which horizontal composites preserve identity 
and composites strictly the compatibility condition on $\sim $ with composites holds automatically.   \hfill   $\square $

\vskip 0.2cm
A category is called {\bf strict} if associativity and unit laws hold strictly up to \, $=$\, , and 
horizontal composites preserve identities and composites strictly. $\sim $ \, still makes sence for strict categories.

\vskip 0.2cm
\proclaim{\bf {Proposition 1.1}} In a strict $\infty $-category $L$ arrows of degree $n$ (i.e., $L^n$) form $1$-category with objects 
$L^0$, arrows $L^n$, domain function $d^n$, codomain function $c^n$. $d,c:L^n\to L^{n-1}$ are $1$-functors. \hfill $\square $
\endproclaim

\proclaim{\bf {Lemma 1.2}} 
\item{$\bullet $} In $\infty $-category $L$ \ $e^k(f\circ _ng)\sim e^kf\circ _{n+k}e^kg$ (when either side is defined).
\item{$\bullet $} $\sim $ is preserved under $\sim $\, , i.e., if $\xymatrix{a \ar[r]^-{\sim }_-{f} & a'}$ is an equivalence with
$\xymatrix{a' \ar[r]^-{\sim }_-{g} & a}$, its quasiinverse ($ea\sim g\circ f$, $ea'\sim f\circ g$), and $f'\sim f$ then $g$ is 
quasiinverse of $f'$ as well.
\item{$\bullet $} A quasiinverse is determined up to $\sim $\, , i.e. if 
$\xymatrix{a \ar@/^0.6pc/[r]^-{f}  &   b   \ar@/^0.5pc/[l]_-{g} \ar@/^1pc/[l]^-{g'} }$ and $g'\circ _1f\sim ea\sim g\circ _1f$ and
$f\circ _1g'\sim eb\sim f\circ _1 g$ then $g'\sim g$.
\vskip 0.1cm\item{$\bullet $} All $n+1$ composites in $\bold{End}(e^na):=L^{0}(e^na,e^na), \, n\ge 0$ coincide up to equivalence $\sim $\, .
\endproclaim
\demo{Proof}
\item{$\bullet $} Assume $f,g\in L^m, \ m\ge n$, then $f\circ _ng=\mu _{d^ng,c^ng,c^nf}(f,g)$ which weakly preserves $e$.
\item{$\bullet $} $ea\sim g\circ f\sim g\circ f'$, \ $ea'\sim f\circ g\sim f'\circ g$.
\item{$\bullet $} $g'=g'\circ _1eb\sim g'\circ _1f\circ _1g\sim g\circ _1f\circ _1g\sim g\circ _1eb=b$.
\item{$\bullet $} $f\circ _{n+1}g=\mu _{a,a,a}(f,g)\sim \mu _{a,a,a}(f\circ _k e^{n+1}a,e^{n+1}a\circ _k g)\sim \mu _{a,a,a}(f,e^{n+1}a)\circ _k\mu _{a,a,a}(e^{n+1}a,g)\sim f\circ _kg$, \, $1\le k\le n+1$.
\hfill $\square $
\enddemo

\vskip 0.1cm
{\bf Definition 1.4.} An arrow $(f:a\to a')\in L^0(a,a')$, \, $deg(a)=deg(a')=m\ge 0$, \, is called
\item{$\bullet $} {\bf monic} if $\forall g,h:z\to a$ if $f\circ _1g\sim f\circ _1h$ then $g\sim h$
\item{$\bullet $} {\bf epic} if $\forall g',h':a'\to w$ if $g'\circ _1f\sim h'\circ _1f$ then $g'\sim h'$
\item{$\bullet $} {\bf equivalence} if $\exists f':a'\to a$ such that $edf\sim f'\circ _1f$ and $edf'\sim f\circ _1f'$ \hfill $\square $

\proclaim{\bf Proposition 1.2} For composable arrows 
\item{$\bullet $} If $f,g$ are monics then $f\circ _1g$ is monic. If $f\circ _1g$ is monic then $g$ is monic
\item{$\bullet $} If $f,g$ are epics then $f\circ _1g$ is epic. If $f\circ _1g$ is epic then $f$ is epic
\item{$\bullet $} If $f,g$ are equivalences then $f\circ _1g$ is equivalence \hfill $\square $
\endproclaim

\vskip 0.1cm
\proclaim{\bf Proposition 1.3} All arrows representing equivalence $a\sim b$ are equivalences.     \hfill    $\square $
\endproclaim 

\vskip 0.1cm
{\bf Definition 1.5.} {\bf $\infty $-functor} $F:L\to L'$ is a function which 
strictly preserves properties of precategories (1)-(2)
\roster
\item if $a\in L^n$ then $F(a)\in L^{'n}$
\item $F(da)=dF(a), F(ca)=cF(a)$
\endroster
and weakly preserves properties (3)-(4)
\roster
\item[3] $F(ea)\sim eF(a)$
\item $F(a\circ _k b)\sim F(a)\circ _k F(b)$ \hfill $\square $
\endroster

\vskip 0.1cm
{\bf Remarks.} 
\item{$\bullet $} We do not require functor $F$ to preserve equivalences $\sim $ because it is not automatic and can be too 
restrictive. However, namely the functors preserving $\sim $ are most important (e.g., see point 1.2).  
\item{$\bullet $} Inverse map $F'$ for a bijective weak functor $F$ is not a functor, in general. If $F$ preserves $\sim $ then for 
the inverse map $F'$ to be a (weak) functor is equivalent to preserving $\sim $ by $F'$. Inverse for a strict functor is 
always a strict functor.    \hfill   $\square $  

\vskip 0.2cm
\proclaim{\bf Lemma 1.3}
\item{$\bullet $} Strict functors preserve equivalences $\sim $\, .
\item{$\bullet $} If functor $F\hskip -0.035cm:\hskip -0.035cmL\hskip -0.035cm\to \hskip -0.035cmL'$ is such that each its 
restriction on hom-sets $F_{a,b}\hskip -0.035cm:\hskip -0.035cmL(a,b)\hskip -0.035cm\to \hskip -0.035cmL'(F(a),F(b))$, 
$a,b\in L^0$, preserves equivalences $\sim $\, , then $F$ preserves equivalences $\sim $\, . 
\item{$\bullet $} If $F\hskip -0.035cm:\hskip -0.035cmL\hskip -0.035cm\to \hskip -0.035cmL'$ is an embedding (injective map) 
such that $\forall a,b\hskip -0.0125cm\in \hskip -0.025cmL^0$ $F_{a,b}\hskip -0.035cm:\hskip -0.035cmL(a,b)\hskip -0.035cm\to \hskip -0.035cmL'(F(a),F(b))$ 
is a strict isomorphism and inverse $F'$ to codomain restriction of 
$\xymatrix{F:L \ar[r]_-{F\hskip -0.3cm\underset {Im(F)}\to |} & Im(F) \ar@{-->}@/_/[l]_-{F'} \ar@{^{(}->}[r] & L'}$
is a functor then $F$ reflects $\sim $\, .
\endproclaim 
\demo{Proof}
\item{$\bullet $} Each arrow presenting a given equivalence $x\sim y$ is between a domain and a codomain which are constructed in a certain way
only by composites and identity operations from arrows of smaller degree presenting the given equivalence and from elements $x$ and $y$. 
A strict functor keeps the structure of the domains and codomains of arrows presenting equivalence $x\sim y$. So, the image of 
arrows presenting equivalence $x\sim y$ will be a family of arrows presenting equivalence $F(x)\sim F(y)$.
\item{$\bullet $} For arrows of degree $>0$ equivalences are preserved by assumption.
Let $\xymatrix{a \ar@/^0.58pc/[r]^-{f}_-{\sim } & b \ar@/^0.58pc/[l]^-{g}_-{\sim } }$, $a,b\in L^0$, be an equivalence for objects in $L$, i.e. 
$ea\sim g\circ _1f$, $eb\sim f\circ _1g$. Then there are two opposite arrows $\xymatrix{F(a) \ar@/^0.58pc/[r]^-{F(f)} & F(b) \ar@/^0.58pc/[l]^-{F(g)} }$.
By assumption, $F(ea)\sim F(g\circ _1f)$, $F(eb)\sim F(f\circ _1g)$. So, $eF(a)\sim F(ea)\sim F(g\circ _1f)\sim F(g)\circ _1F(f)$ 
and $eF(b)\sim F(eb)\sim F(f\circ _1g)\sim F(f)\circ _1F(g)$. Therefore, $\xymatrix{F(a) \ar@/^0.58pc/[r]^-{F(f)}_-{\sim } & F(b) \ar@/^0.58pc/[l]^-{F(g)}_-{\sim } }$ 
is an equivalence.
\item{$\bullet $} Inverse to a strict isomorphism is a strict isomorphism, i.e. preserves equivalences. So, $F'$ is a functor 
which preserves equivalences in all hom-sets and, consequently, preserves all equivalences. Preservation of equivalences for $F'$ 
is exactly reflection of equivalences for $F$.         \hfill    $\square $
\enddemo 

\proclaim{\bf Lemma 1.4}
\item{$\bullet $} $x=y$ \, iff \, $ex\sim ey$ \, [in particular, \, $=$ \, is definable via \, $\sim $\, ].
\item{$\bullet $} Functors, preserving \, $\sim $\, , strictly preserve all composites \, $\circ _k$, $k\ge 1$.
\item{$\bullet $} Functors, weakly preserving $e^2$, strictly preserve $e$, i.e. $e^2F(a)\sim F(e^2a)$ $\Rightarrow $ $eF(a)=F(ea)$.
\item{$\bullet $} {\bf Quasiiequal} functors (i.e. $F(f^n)\sim G(f^n)$ for all $f^n\in L^n$, $n\ge 0$) are equal.
\endproclaim 
\demo{Proof}
\item{$\bullet $} $x=y$ $\Rightarrow $ $ex=ey$ $\Rightarrow $ $ex\sim ey$. Conversely, $ex\sim ey$ $\Rightarrow $ $dex=dey$ 
$\Rightarrow $ $x=y$. 
\item{$\bullet $} Sufficient to prove $eF(f\circ _kg)\sim e(F(f)\circ _kF(g))$, but it holds 
$eF(f\circ _kg)\sim F(e(f\circ _kg))\sim \text{($F$ preserves $\sim $) }F((ef)\circ _{k+1}(eg))\sim F(ef)\circ _{k+1}F(eg)\sim eF(f)\circ _{k+1}eF(g)\sim e(F(f)\circ _kF(g))$.
\item{$\bullet $} $e^2F(a)\sim F(e^2a)$ $\Rightarrow $ $de^2F(a)=dF(e^2a)$ $\Rightarrow $ $eF(a)=F(ea)$.
\item{$\bullet $} Again, it is sufficient to prove $eF(f^n)\sim eG(f^n)$. 
\item{}$eF(f^n)\sim F(ef^n)\sim \text{(by assumption) }G(ef^n)\sim eG(f^n)$.   \hfill   $\square $
\enddemo 

\vskip 0.0cm
{\bf Corollary.} {\bf $\infty $-categories in the sense of definition 1.3 are {\rm almost} strict}, namely, with strict associativity, 
identity, and interchange laws.
\demo{Proof} Strict associativity and strict identity laws hold because by the axioms functors $L(x,y)\times L(y,z)\times L(z,t)\to L(x,t):(f^n,g^n,h^n)\mapsto (h^n*g^n)*f^n$ and
$L(x,y)\times L(y,z)\times L(z,t)\to L(x,t):(f^n,g^n,h^n)\mapsto h^n*(g^n*f^n)$, $deg(x)=deg(y)=deg(z)=deg(t)$, are quasiequal, and, respectively, functors 
$L(x,y)\to L(x,y):f\mapsto f$ and $L(x,y)\to L(x,y):f\mapsto ey*f$, $deg(x)=deg(y)$ (similar for the right identity), are quasiequal. Strict interchange law
is because functor $L(x,y)\times L(y,z):(f,g)\mapsto g*f$ preserves $\sim $.     \hfill  $\square $
\enddemo 

\vskip 0.1cm
{\bf Definition 1.6.} For given two functors $F$, $G$ \ {\bf $\infty $-natural transformation} $\alpha :F\to G$ \ is a function
$\alpha :L^0\to L^{'1}:a\mapsto (\xymatrix{F(a)\ar[r]^-{\alpha (a)} & G(a)})$ such that
$$\mu _{F(a),F(b),G(b)}(e^k\alpha (b), F(f))\sim \mu _{F(a),G(a),G(b)}(G(f),e^k\alpha (a))$$ for all $f\in L^{k}(a,b)$,
$k=0,1,...$ \hfill $\square $

{\bf Definition 1.7.} For given two functors $F$, $G$ and two natural transformations
$\xymatrix{F \ar@<1ex>[r]^-{\alpha }
\ar@<-1ex>[r]_-{\beta } & G}$
{\bf $\infty $-modification} $\lambda :\alpha \to \beta $ is a function $\lambda :L^0 \to L^{'2}:a\mapsto (\xymatrix{\alpha (a) \ar[r]^-{\lambda (a)} & \beta (a)})$ such that
$$\mu _{F(a),F(b),G(b)}(e^k\lambda (b), F(f))\sim \mu _{F(a),G(a),G(b)}(G(f),e^k\lambda (a))$$
for all $f\in L^{k+1}(a,b)$,
$k=0,1,...$ \hfill $\square $
\vskip 0.3cm
\hskip -0.39cmAnalogously, modifications of higher order are introduced. Call modifications by $1$-modifications, natural transformations 
by $0$-modifications.

{\bf Definition 1.8.} Given two functors $F$, $G$, two $0$-modifications $\xymatrix{F \ar@<1ex>[r]^-{\alpha ^0_1} \ar@<-1ex>[r]_-{\alpha ^0_2} & G}$, \newline
two $1$-modifications
$\xymatrix{\alpha ^0_1 \ar@<1ex>[r]^-{\alpha ^1_1} \ar@<-1ex>[r]_-{\alpha ^1_2} & \alpha ^0_2}$,..., two $n-1$-modifications
$\xymatrix{\alpha ^{n-2}_1 \ar@<1ex>[r]^-{\alpha ^{n-1}_1} \ar@<-1ex>[r]_-{\alpha ^{n-1}_2} & \alpha ^{n-2}_2}$ \newline
{\bf $\infty $-$n$-modification} $\alpha ^n:\alpha ^{n-1}_1\to \alpha ^{n-1}_2$ is a function
$\alpha ^n:L^0\to L^{'n+1}:$\newline
$a\mapsto (\xymatrix{\alpha ^{n-1}_1(a) \ar[r]^-{\alpha ^n(a)} & \alpha ^{n-1}_2(a)})$ such that
$$\mu _{F(a),F(b),G(b)}(e^k\alpha ^n (b), F(f))\sim \mu _{F(a),G(a),G(b)}(G(f),e^k\alpha ^n (a))$$
for all $f\in L^{k+n}(a,b)$,
$k=0,1,...$ \hfill $\square $

\vskip 0.2cm
{\bf Corollary.} {\bf All $n$-modification in the sense of definition 1.8 are strict.}
\demo{Proof} By condition two functors $\alpha ^n(b)*F(-):L^{\ge n}(a,b)\to L^{'\ge n}(F(a),G(b))$ and 
$G(-)*\alpha ^n(a):L^{\ge n}(a,b)\to L^{'\ge n}(F(a),G(b))$ are quasiequal and, so, equal.   \hfill   $\square $
\enddemo

\vskip 0.0cm
{\bf Definition 1.9.} {\bf $\infty $-CAT} is an $\infty $-category consisting of
\item{$\bullet $} graded set $C=\coprod\limits _{n\ge 0}C^n$, where $C^0$ are categories, $C^1$ functors, $C^n$ $(n-2)$-modifications
\item{$\bullet $} if $\alpha ^n:\alpha ^{n-1}_1\to \alpha ^{n-1}_2\in C^n$ then $d\,\alpha ^n=\alpha ^{n-1}_1$, $c\,\alpha ^n=\alpha ^{n-1}_2$
\item{$\bullet $} $e\,\alpha ^n\in C^{n+1}$ is the map $L^0\to L^{'(n+1)}:a\mapsto e(\alpha ^n(a))$
\item{$\bullet $} for given two n-modifications $\alpha ^n_1$, $\alpha ^n_2$ such that $d^k\alpha ^n_1=c^k\alpha ^n_2$
$$\hskip -0.05cm\alpha ^n_1\circ _k\alpha ^n_2:=\cases
a\mapsto (\alpha ^n_1(a)\circ _k \alpha ^n_2(a))  & \text{if} \ \ k<n+2 \\
a\mapsto (\alpha ^n_1(F'(a))\circ _{(n+1)}G(\alpha ^n_2(a)))  & \text{if}
\ \ k=n+2, \ F'=c^{(n+1)}\alpha ^n_2, \ G=d^{(n+1)}\alpha ^n_1
\endcases $$
First composite works when $\alpha ^n_1, \alpha ^n_2\in \infty \text{\bf -CAT}(L,L')$, second when $\alpha ^n_1\in \infty \text{\bf -CAT}(L',L'')$
$\alpha ^n_2\in \infty \text{\bf -CAT}(L,L')$, where $L,L',L''$ are categories. \hfill $\square $

\vskip 0.35cm
\proclaim{\bf {Proposition 1.4}} Categories, functors, natural transformations, modifications, etc. constitute $\infty $-category {\bf $\infty $-CAT} of $\infty $-categories. \hfill $\square $
\endproclaim

\vskip 0.0cm
{\bf Definition 1.10.} Category $L$ is called {\bf $\infty $-$n$-category} if $L^{j+1}=e(L^j)$ for $j\ge n$.
\hfill $\square $

$L/\hskip -0.1cm\sim $ is not a category in general since $\sim $ is not compatible with $e$. However, if we take quotient only on a fixed 
level $n$ and make all higher arrows identities we get $\infty $-$n$-category $L^{(n)}$, $n$-th approximation of $L$. 
Generally there are no functors $\xymatrix{L^{(n)} \ar@{^{(}->}[r] & L}$, $\xymatrix{L \ar@{->>}[r] & L^{(n)}}$ (except for the last 
surjection if $L$ is a weak $\infty $-$(n+1)$-category and all $(n+1)$-arrows are iso's).

\subhead {} 1.a. Weak categories, functors, natural transformations, modifications \endsubhead 

As we saw above, using a weak language (substitution $\sim $ instead of $=$) does not give a weak category theory. The only 
advantage was that we could deal with $\sim $ instead of $=$ (which is important for the classification problem that still makes sense 
for strict $\infty $-categories). All known definitions of weak categories \cite{C-L, Lei, Koc, etc.} are nonelementary (at least, they use functors, natural transformations, operads, monads just 
for the very definition). Probably, this is a fundamental feature of weak categories. To introduce them we also need the whole universe 
$\infty \text{-}\bold{PreCat}$ of $\infty $-precategories.

\vskip 0.2cm
{\bf Definition 1.a.1.} $\infty \text{-}\bold{PreCat}$ consists of 
\item{$\bullet $} {\bf $\infty $-precategories} (definition 1.1) together with $\sim $-relation in each [$\sim $ may be not transitive],
\item{$\bullet $} {\bf $\infty $-functors} (definition is like 1.5 for $\infty $-categories), i.e. functions $F:L\to L'$ of degree $0$
preserving $d$ and $c$ strictly, and $e$ and $\circ _k$, $k\ge 1$, weakly, 
\item{$\bullet $} {\bf lax $\infty \text{-}n$-modifications}, $n\ge 0$, i.e. {\bf total} maps $\alpha ^n:L\to L'$ (with variable degree on 
different elements, but $\le n+1$, more precisely, the induced map $\Bbb N\to \Bbb N:deg(x)\mapsto (deg(\alpha ^n(x))-deg(x))$ 
is an antimonotone map, decreasing by 1 at each step from $n+1$ at $deg(x)=0$ to 1 at $deg(x)=n$ and remaining constant 1 after) 
being defined for a given sequence of two functors $F,G:L\to L'$, two $0$-modifications (natural transformations) $\alpha ^0_1,\alpha ^0_2:F\to G$, ...,
two $(n-1)$-modifications $\alpha ^{n-1}_1,\alpha ^{n-1}_2:\alpha ^{n-2}_1\to \alpha ^{n-2}_2$ \ \ \ as \hskip 0.5cm $\alpha ^n:=$
\item{}\hskip -0.25cm$\cases  \hskip -0.05cm(\alpha ^n(x):\alpha ^{n-1}_1(x)\to \alpha ^{n-1}_2(x))\in L^{'n}(F(x),G(x))  &  \hskip 0.25cmx\in L^0 \hskip -0.3cm \\
\hskip -0.05cm\alpha ^n(x):=\alpha ^n(e^{n+1-k}x)\in L^{'n+1}(F(d^kx),G(c^kx))  &   \hskip 0.15cm\underset {0<k<n+1}\to {x\in L^k}     \hskip -0.3cm  \\
\hskip -0.05cm(\alpha ^n(x):\alpha ^n(c^{n+1}x)\circ _{n+1}F(x)\to G(x)\circ _{n+1}\alpha ^n(d^{n+1}x))\in & \hskip 0.0cmx\in L^{n+1} \hskip -0.3cm \\        
\hskip 4.4cm\in L^{'n+1}(F(d^{n+1}x),G(c^{n+1}x))	&     \hskip -0.3cm  \\
\hskip -0.05cm(\alpha ^n(x):\alpha ^n(cx)\circ _1(e\alpha ^n(c^{n+2}x)\circ _{n+2}F(x))\to     & \hskip 0.0cmx\in L^{n+2} \hskip -0.3cm \\
\hskip -0.05cm(G(x)\circ _{n+2}e\alpha ^n(d^{n+2}x))\circ _1\alpha ^n(dx))\in L^{'n+2}(F(d^{n+2}x),G(c^{n+2}x))   &          \hskip -0.3cm  \\
\hskip -0.05cm\alpha ^n(x):\alpha ^n(cx)\circ _1(e\alpha ^n(c^2x)\circ _2(e^2\alpha ^n(c^{n+3}x)\circ _{n+3}F(x)))\to     &  \hskip 0.0cmx\in L^{n+3}   \hskip -0.3cm  \\
\hskip -0.05cm((G(x)\circ _{n+3}e^2\alpha ^n(d^{n+3}x))\circ _2e\alpha ^n(d^2x))\circ _1\alpha ^n(dx)\in L^{'n+3}(F(d^{n+3}x),G(c^{n+3}x)) &             \hskip -0.3cm  \\
  \dots                      &                            \\
\hskip -0.05cm\alpha ^n(x):    &   \hskip 0.0cmx\in L^{n+m}   \hskip -0.3cm  \\
\hskip -0.05cm\alpha ^n(cx)\circ _1\dots \circ _{m-2}(e^{m-2}\alpha ^n(c^{m-1}x)\circ _{m-1}(e^{m-1}\alpha ^n(c^{n+m}x)\circ _{n+m}F(x)\underbrace{)\dots )}_{m-1}\to &    \hskip -0.3cm  \\
\hskip -0.05cm\underbrace{(\dots (}_{m-1}G(x)\circ _{n+m}e^{m-1}\alpha ^n(d^{n+m}x))\circ _{m-1}e^{m-2}\alpha ^n(d^{m-1}x))\circ _{m-2}\dots \circ _1\alpha ^n(dx)\in  &       \hskip -0.3cm  \\
\hskip 4.4cm\in L^{'n+m}(F(d^{n+m}x),G(c^{n+m}x))  &       \hskip -0.3cm \\
\dots               &                \hskip -0.3cm \endcases $   
\item{}$d\alpha ^n:=\alpha ^{n-1}_1, \ c\alpha ^n:=\alpha ^{n-1}_2$ \ [$(d\alpha ^n)(x)\ne d(\alpha ^n(x)), \ (c\alpha ^n)(x)\ne c(\alpha ^n(x))$ if $deg(x)>0$].   \hfill    $\square $

\vskip0.2cm
{\bf Remarks.} 
\item{$\bullet $} $\infty \text{-}n$-modifications look terrible but it is the weakest form of naturality (infinite sequence 
of naturality squares for naturality squares). To deal with such entities a kind of operads is needed.
\item{$\bullet $} To give $n$-modification $\alpha ^n$ is the same as to give a map $\alpha ^n\underset {L^0}\to |:L^0\to L'$ of degree
$n+1$ and $\forall a,b\in L^0$ a natural transformation $\nu ^{\alpha ^n}_{a,b}:\alpha ^n(b)*F(-)\to G(-)*\alpha ^n(a):L^{\ge n}(a,b)\to L^{'\ge n}(F(a),G(b))$, 
where $F=d^{n+1}\alpha ^n$, $G=c^{n+1}\alpha ^n$.
\item{$\bullet $} When $\alpha ^n(x)$, $deg(x)>0$, are all identities (of the required types) $\infty \text{-}n$-modifications are 
called {\bf strict}. They are usual modifications and composable like in definition 1.9 when universe $\infty \text{-}\bold{CAT}$ 
is strict (in that case strict modifications are weak as well). In a weak universe $\infty \text{-}\bold{CAT}$ strict modifications
need not to be weak (i.e. to be modifications at all).
\item{$\bullet $} $\infty \text{-}\bold{PreCAT}$ is not an $\infty $-precategory itself because there are no identities and composites 
for weak $n$-modifications. 
The problem with identities and composites is not clear, if they exist at all without making either naturality condition or 
$\infty $-categories stricter. 
\item{$\bullet $} In general, these two sides 'categories and functors' and '$n$-modifications' form a strange pair. If we weaken
one of these sides the other one becomes stricter (under condition that $\infty \text{-}\bold{CAT}$ is a (let it be very weak) {\bf category}). 
So, the following {\bf hypothesis} holds:\vskip 0.1cm
\item{}{\it There is no $\infty \text{-}\bold{CAT}$ with {\bf simultaneously weak} categories, functors, and $n$-modification}. 
\vskip 0.1cm\item{} For example, if we want weak modifications and want them to be composable we need to introduce several axioms on categories, 
one of which is like '\, {\it $\forall a,b\in L^0$ and  $\forall $ functors $F,G:L\to L'$ if $\exists $ natural transformations 
$\alpha :f_1*F(-)\to G(-)*g_1:L^{\ge n}(a,b)\to L^{'\ge n}(F(a),G(b))$ and 
$\beta :f_2*F(-)\to G(-)*g_2:L^{\ge n}(a,b)\to L^{'\ge n}(F(a),G(b))$ and $n+1$-cells 
$f_1,f_2$ and $g_1,g_2$ are $\circ _k$-composable then $\exists $ a natural transformation ($k$-composite)
$\gamma :(f_1\circ _kf_2)*F(-)\to G(-)*(g_1\circ _kg_2):L^{\ge n}(a,b)\to L^{'\ge n}(F(a),G(b))$}\, '. But such axioms make 
categories very special. From the other side if we want categories to be weak we need to make stricter (maybe, strict) $n$-modifications
in order they would be composable. The problem is in existence of composites (and units) for weak $n$-modifications. 
\item{$\bullet $} Instead of lax $n$-modifications we could use modifications with $\alpha ^n(x)$ being $\sim $ for $deg(x)>0$ in $L'$. In both cases 
in order to make horizontal composites (at least, $F*\alpha ^n:=F\circ _{\bold{SET}}\alpha ^n$) we need functors preserving composites
(or composites and $\sim $), i.e. 'weak modifications' $\Rightarrow $ 'strict functors'.
\item{$\bullet $} If the above hypothesis was true it could be nice, e.g. a universe where $\infty \text{-}\bold{Top}$ 
lives would contain only strict $n$-modifications.  \hfill   $\square $

\vskip 0.2cm
{\bf Definition 1.a.2.} {\bf Weak $\infty $-category $L$} is an $\infty $-precategory (see definition 1.1)
such that \, [all elements below are supposed to be composable when we write composite for them]
\item{$\bullet $} $\sim $ is transitive \, $x\sim y\sim z$ $\Rightarrow $ $x\sim z$,
\item{$\bullet $} horizontal composites \, $*$ \, strictly preserve properties (1)-(2) of precategories
\itemitem{(1)} $deg(x*y)=deg(x)=deg(y)$ \, if \, $deg(x)=deg(y)$ \, (interchange law for degree)
\itemitem{(2)} $d(x*y)=(dx)*(dy)$, $c(x*y)=(cx)*(cy)$ \, if \, $deg(x)=deg(y)$ \, (interchange law for domain and codomain)
\item{}and weakly preserve properties (3)-(4) of precategories
\itemitem{(3)} $e(x*y)\sim (ex)*(ey)$ \, if \, $deg(x)=deg(y)$ \, (interchange law for identity)
\itemitem{(4)} $(x\circ _ky)*(z\circ _kt)\sim (x*z)\circ _k(y*t)$ \, if \, $deg(x)=deg(y)=deg(z)=deg(t)$ \, 
(interchange law for composites) \, [$\circ _k$ has smaller 'deepness' $k$ than the given $*=\circ _n$, $n>k$],
\item{$\bullet $} ({\bf weak associativity}) 
\item{}$\forall x,y,z,t\in L^n$ for two functors 
$l_{x,y,z,t}:L(x,y)\times L(y,z)\times L(z,t)\to L(x,t):(f,g,h)\mapsto (h*g)*f$ and 
$r_{x,y,z,t}:L(x,y)\times L(y,z)\times L(z,t)\to L(x,t):(f,g,h)\mapsto h*(g*f)$ 
$\exists $ natural transformation $\alpha _{x,y,z,t}:l_{x,y,z,t}\to r_{x,y,z,t}$, 
\item{$\bullet $} ({\bf weak unit}) 
\item{}$\forall x,y\in L^n$ and functors $u^{l}_{x,y}:L(x,y)\to L(x,y):f\mapsto ey*f$ and 
$u^{r}_{x,y}:L(x,y)\to L(x,y):f\mapsto f*ex$ $\exists $ natural transformations
$\epsilon ^l_{x,y}:u^l_{x,y}\to Id$ and $\epsilon ^r_{x,y}:Id\to u^r_{x,y}$.      \hfill   $\square $

\vskip 0.2cm
{\bf Remarks.} 
\item{$\bullet $} We do not introduce a universe $\infty \text{-}\bold{CAT}$ with weak categories, functors and $n$-modifications because
there is no (at least, obvious) units and composites for $n$-modifications (however, identity natural transformations exist if only 
vertical composites of natural transformations are defined, for if $F:L\to L'$ is a functor take $(eF)(a):=e(F(a)), \, a\in L^0$ and by the weak unit law $\forall a,b\in L^0$ 
$\exists $ a natural transformation $\nu _{a,b}:e(F(b))*F(-)\to F(-)*e(F(a)):L^{\ge 0}(a,b)\to L^{'\ge 0}(F(a),F(b))$, take 
$\nu _{a,b}:=(\epsilon ^{u^r}_{F(a),F(b)}\circ _1\epsilon ^{u^l}_{F(a),F(b)})*F:=(\epsilon ^{u^r}_{F(a),F(b)}\circ _1\epsilon ^{u^l}_{F(a),F(b)})\circ _{\bold{SET}}F$). 
The problem is what are the weakest conditions on categoies,
functors and $n$-modifications in order they form a category. Maybe, there are several independent such conditions and, so, several
categories living in $\infty \text{-}\bold{CAT}$ with weakest entities.
\item{$\bullet $} To keep a usual form of (weak) associativity and (weak) unit we could introduce relations $\underset k\to {\sim }$
for elements of images of two functors $F,G:L\to L'$ connected by a natural transformation $\alpha :F\to G$, namely, 
$x\underset k\to {\sim }y$ if $\exists $ $z\in L^k$ such that $x=F(z), \, y=G(z)$. These relations are not refrexive, symmetric and 
transitive. Then we could write associativity and unit laws as $(x\circ _ky)\circ _kz\underset {k-1}\to {\sim }x\circ _k(y\circ _kz)$
and $e^kc^kx\circ _kx\underset {k-1}\to {\sim }x$, \, $x\underset {k-1}\to {\sim }x\circ _ke^kd^kx$. Under assumption that composites and 
units exist in an $\infty \text{-}\bold{CAT}$ we could choose more sensible piece of $\infty \text{-}\bold{CAT}$ with categories in 
which $\underset 0\to {\sim }\, \equiv \, \sim $ and all $\underset k\to {\sim }$ are symmetric and transitive by the requirement that
$\alpha _{x,y,z,t}, \, \epsilon ^l_{x,y}, \, \epsilon ^r_{x,y}$ are equivalences.          \hfill    $\square $

\vskip 0.35cm
\centerline{\bf Examples}
\vskip 0.2cm
\item{1.} {\bf $\infty $-Top} is an $\infty $-category with homotopies for homotopies as higher order cells.
\item{2.} {\bf $\infty $-Diff} is an $\infty $-category of differentiable manifolds in the same way as {\bf $\infty $-Top}.
\item{3.} {\bf $\infty $-TopALg} is an $\infty $-category of topological algebras in the same way as {\bf $\infty $-Top} where each instance of
homotopy is a homomorphism of topological algebras.
\item{4.} {\bf $2$-Top} is a strict $\infty $-$2$-category with 2-cells, homotopy classes of homotopies, and just identities in 
higher order (\, $\sim $ \, on the level of objects means homotopy equivalence of spaces, on the level of $1$-arrows homotopness of maps, and on the
level $\ge 2$ coincidence). $2\text{-}\bold{Cat}$ is similar.
\item{5.} {\bf $\infty $-Compl} is an $\infty $-category of (co)chain complexes with (algebraic) homotopies for homotopies as higher order cells.
\item{6.} For $1$-category $A$, \ $A_{equiv}$ is a strict $\infty $-$2$-category such that $A^0_{equiv}=A^0$, $A^1_{equiv}=\Bigl\{f\in A \ \Bigl|\Bigr.\vcenter{\xymatrix{{\bullet } \ar[r]^-{\exists H}_-{\sim } & {\bullet }\\
{\bullet } \ar[u]^{f} \ar[r]_-{\forall h}^-{\sim } & {\bullet } \ar[u]_{f}}}
\Bigr\}$, $A^2_{equiv}=\Bigl\{\text{iso's} \ \Bigl|\Bigr. \ \forall f,g\in A^1_{equiv} \ \exists! \ \xymatrix{f \ar[r]^-{\gamma }_-{\sim } & g} \ \text{iff} \ \vcenter{\xymatrix{{\bullet } \ar[r]^-{\exists H}_-{\sim } & {\bullet }\\
{\bullet } \ar[u]^{f} \ar[r]_-{\forall h}^-{\sim } & {\bullet } \ar[u]_{g}}} \Bigr\}$. \newline
$A_{equiv}$ contains all equivariant maps $f:X\to Y$ with respect to a group homomorphism $\rho :\bold{Aut}(X)\to \bold{Aut}(Y)$.
\item{7.} (weak) {\bf covariant} {\bf $\infty $-Hom}-functor $L(a,-):L\to \infty \text{\bf -CAT}:$
\vskip -0.3cm
$$\cases
b\mapsto L(a,b) & \ b\in L^0 \\
(f:b\to b')\mapsto (L(a,f):g\mapsto \mu (e^kf,g))& \ f\in L^0(b,b'), \ g\in L^{k}(a,b) \\
(\alpha :f\to f') \mapsto (L(a,\alpha ):x\mapsto \mu (\alpha ,ex ) & \alpha \in L^1(b,b'), \ x\in L^0(a,b) \\
(\delta :\alpha \to \alpha ') \mapsto (L(a,\delta ):x\mapsto \mu (\delta ,e^2x)) & \delta \in L^2(b,b'), \ x\in L^0(a,b) \\
\dots & \\
(\alpha ^n:\alpha ^{(n-1)}_1\to \alpha ^{(n-1)}_2) \mapsto (L(a,\alpha ^n):x\mapsto \mu (\alpha ^n,e^nx)) & \alpha ^n \in L^n(b,b'), \ x\in L^0(a,b)\\
\dots &
\endcases$$
\vskip -0.05cm
\item{8.} {\bf opposite category} $L^{op}$ is an $\infty $-category such that
\itemitem{$\bullet $} $(L^{op})^n=L^n$, $n\ge 0$
\itemitem{$\bullet $} \hbox{$d^{op}(\alpha ^n)=\cases d(\alpha ^n) & \text{if} \ n\ge 2\\
c(\alpha ^n) & \text{if} \ n=1 \endcases $} \hskip 1cm
\hbox{$c^{op}(\alpha ^n)=\cases c(\alpha ^n) & \text{if} \ n\ge 2\\
d(\alpha ^n) & \text{if} \ n=1 \endcases $}
\itemitem{$\bullet $} $e^{op}=e$
\itemitem{$\bullet $} \hbox{$\beta ^n\circ ^{op}_k\alpha ^n=\cases
\beta ^n\circ _k\alpha ^n & \text{if } \ \alpha ^n,\beta ^n\in L^n, \ k<n \\
\alpha ^n\circ _k\beta ^n & \text{if } \ \alpha ^n,\beta ^n\in L^n, \ k=n
\endcases $}
\hskip 0.3cm
\hbox{(for composable elements)}
\vskip 0.2cm
\item{9.} (weak) {\bf contravariant} {\bf $\infty $-Hom}-functor $L(-,b):L^{op}\to \infty \text{\bf -CAT}:$
$$\cases
a\mapsto L(a,b) & \ a\in L^0 \\
(f:a\to a')\mapsto (L(f,b):g\mapsto \mu (g,e^kf))& \ f\in L^0(a,a'), \ g\in L^{k}(a',b) \\
(\alpha :f\to f') \mapsto (L(\alpha ,b):x\mapsto \mu (ex,\alpha ) & \alpha \in L^1(a,a'), \ x\in L^0(a',b) \\
(\delta :\alpha \to \alpha ') \mapsto (L(\delta ,b):x\mapsto \mu (e^2x,\delta )) & \delta \in L^2(a,a'), \ x\in L^0(a',b) \\
\dots & \\
(\alpha ^n:\alpha ^{(n-1)}_1\to \alpha ^{(n-1)}_2) \mapsto (L(\alpha ^n ,b):x\mapsto \mu (e^nx,\alpha ^n )) & \alpha ^n \in L^n(a,a'), \ x\in L^0(a',b)\\
\dots &
\endcases$$
\item{10.} {\bf Yoneda embedding} \ $\bold Y:L\to \infty\text{\bf -CAT}^{L^{op}}:\alpha \mapsto L(-,\alpha )$, \ $\alpha \in L$.
\item{11.} $\bold{Set}$ is simultaneously an object and a full subcategory of $\infty \text{-}\bold{CAT}$.
\item{12.} A (big) set $L_{\sim }:=\coprod\limits _{n\ge 0}L^n_{\sim }$, where $L^n_{\sim }$ are defined recursively as 
$L^0_{\sim }:=L^0$ and $L^n_{\sim }$ are all equivalences from $L^n$ with domain and codomain in $L^{n-1}_{\sim }$, is a subcategory
of $L$. Similarly, $L_{k\sim }:=\coprod\limits _{n\ge 0}L^n_{k\sim }, \, k\ge 0$, where 
$L^n_{k\sim }:=\cases \hskip -0.05cmL^n & n\le k \hskip -0.05cm\\
\hskip -0.05cm\text{equivalences from }L^n\text{ with dom and codom in }L^{n-1}_{k\sim } & n>k \hskip -0.05cm\endcases $, is a subcategory of $L$. From this point
$L_{\sim }=L_{0\sim }$. Such categories are most important for the classification problem (up to $\sim $).
Sometimes, 'invariants' can be constructed only for $L_{\sim }$ (see point 2.1).
\item{13.} {\bf Higher order concepts} can {\bf simplify} proof of first order facts. E.g., each strict $2$-functor 
$\Phi :2\text{-}\bold{CAT}\to 2\text{-}\bold{CAT}$, where $2\text{-}\bold{CAT}$ is the usual strict category of categories, functors, and natural transformations, preserves adjunction $\biggl ($indeed, triangle identities $\cases \varepsilon G\circ G\eta =1_G   & \hskip -0.35cm\\
F\varepsilon \circ \eta F=1_F  & \hskip -0.35cm\endcases $ are respected by $\Phi $ $\cases \Phi (\varepsilon)\Phi (G)\circ \Phi (G)\Phi (\eta )=1_{\Phi (G)}   & \hskip -0.35cm\\
\Phi (F)\Phi (\varepsilon )\circ \Phi (\eta )\Phi (F)=1_{\Phi (F)}  & \hskip -0.35cm\endcases $ $\biggr )$. It gives short proofs of the following results.
\item{a)} {\it Right adjoints preserve limits (left adjoints preserve colimits).}
\item{}{\smc Proof}. \hskip 1cm$\vcenter
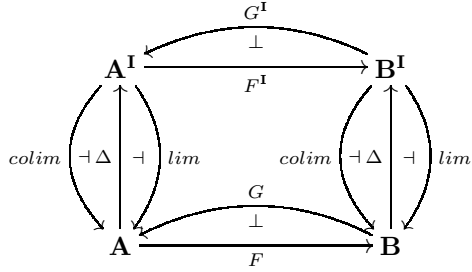
{\xymatrix{\bold{A}^{\bold{I}} \ar[rrr]_-{F^{\bold{I}}} \ar@/^1.25pc/[dd]^-{lim} \ar@/_1.6pc/[dd]_-{colim}^-{\dashv } & & & \bold{B}^{\bold{I}} \ar@/^1.25pc/[dd]^-{lim} \ar@/_1.6pc/[dd]_-{colim}^-{\dashv } \ar@/_1.25pc/[lll]^-{\perp }_-{G^{\bold{I}}} \\
    & & & \\
\bold{A} \ar[uu]^-{\Delta }_-{\, \dashv } \ar[rrr]_-{F} & & & \bold{B} \ar[uu]^-{\Delta }_-{\, \dashv } \ar@/_1.25pc/[lll]_-{G}^-{\perp } }}$
\item{}where \, $(-)^{\bold{I}}\equiv 2\text{-}\bold{CAT}(\bold{I},-):2\text{-}\bold{CAT}\to 2\text{-}\bold{CAT}$ is a hom-$2$-functor.
\item{}Now, $G^{\bold{I}}\circ \Delta =\Delta \circ G$ (obvious). Taking right adjoints of both sides completes the proof 
$lim\circ F^{\bold{I}}\simeq F\circ lim$ (for colimits the same argument works $F^{\bold{I}}\circ \Delta =\Delta \circ F\, \Rightarrow \, colim \circ G^{\bold{I}}\simeq G\circ colim$).   \hfill $\square $
\item{b)} {\it Each $1\text{-}\bold{Cat}$-valued presheaf admits a sheafification ($1\text{-}\bold{Cat}$ is a category of small categories and functors between them)}.
\item{}{\smc Proof}. $1\text{-}\bold{Cat}$-valued presheaf on $\bold{C}$ is the same as an internal category object in $\bold{Set}^{\bold{C}^{op}}$.
There is an adjoint situation $\xymatrix{\bold{Sh}(\bold{C}) \ar@{^{(}->}[r] & \bold{Set}^{\bold{C}^{op}} \ar@/_1.1pc/[l]^-{\perp } }$
in $\bold{LEX}$, where $\bold{LEX}\hookrightarrow 2\text{-}\bold{CAT}$ is a $2$-category of finitely complete categories, functors 
preserving finite limits, and (arbitrary) natural transformations. There is a $2$-functor $\bold{CAT}(-):\bold{LEX}\to 2\text{-}\bold{CAT}$
assigning to each category in $\bold{LEX}$ the category of its internal category objects and to each functor and natural transformation 
the induced ones. Then $\exists $ an adjunction $\xymatrix{\bold{CAT}(\bold{Sh}(\bold{C})) \ar[r] & \bold{CAT}(\bold{Set}^{\bold{C}^{op}}) \ar@/_1.1pc/[l]^-{\perp } }$
which means that each $1\text{-}\bold{Cat}$-valued presheaf can be sheafified by the top curved arrow.    \hfill  $\square $

\subhead {} 1.1. Fractal organization of the new universum\endsubhead

{\bf Fractal Principle.} Object $A$ with properties $\{P_i\}_I$ has fractal structure if there are subobjects $\{A_j\}_J$ which relate
to each other in a certain way (express it by additional property $P=$'to have $|J|$ subobjects which relate in the certain way') and
each $A_j$ inherits all properties $\{P_i\}_I\&P$.

\hfill $\square $

It can be useful to see that in spite of complicated structure each $\infty $-category and, moreover, {\bf $\infty $-CAT} has a regular structure
which is repeated for certain arbitrary small pieces. Such pieces are, of course, hom-sets $L(a,b)$ which inherit all properties (1)-(4),
associativity and identity laws, and each piece of which still has the same structure. In particular, $L(a,b)(c,d)=L(c,d)$. $\infty $-functor restricted to such a piece is again $\infty $-functor.
Moreover, each $\infty $-category can be regarded as a hom-set of a little bit bigger category if we formally attach two distinct elements 
$\alpha , \beta \in L^{-1}$ with their identities of higher order $e^n(\alpha ),\ e^n(\beta ), \ n\ge 1$ 
(such that $d(L^0)=\alpha ,\ c(L^0)=\beta $ and composites with these identities of other elements hold strictly).
Other natural pieces of $L$ which inherit all properties and are $\infty $-categories are $L^{\ge n}$, $L^{\ge n}(a,b)$ (elements of degree not lower than $n$).

\subhead {} 1.2. Notes on Coherence Principle\endsubhead

This principle is an axiom to deal with equivalence relation $\sim $. It is not logically necessary for higher order category theory itself. There
can be categories in which it does not hold.

{\bf Coherence Principle.} For a given set of cells $\{a_i\}_I$ and a given set of base equivalences
$\{t_j(\{a_i\}_I)\sim s_j(\{a_i\}_I)\}_J$ for any two constructions $F_1(\{a_i\}_I)$ and $F_2(\{a_i\}_I)$ and any two derived
equivalences $\varepsilon ^0_i:F_1(\{a_i\}_I)\sim F_2(\{a_i\}_I)$, $i=1,2$ there are derived equivalences $\varepsilon ^1_{m}:
\varepsilon ^0_1\sim \varepsilon ^0_2$, $m\in M^1$, such that for any two of them $\varepsilon ^1_{m_1}$, $\varepsilon ^1_{m_2}$ there are derived
equivalences $\varepsilon ^2_{m}:\varepsilon ^1_{m_1}\sim \varepsilon ^1_{m_2}$, $m\in M^2$ again such that for any pair of them $\varepsilon ^2_{m_1}$, $\varepsilon ^2_{m_2}$
there are derived equivalences of higher order, etc. \hfill $\square $

Here constructions mean application of composites, functors, natural transformations,.. to $\{a_i\}_I$. Derived equivalences mean
eqivalences obtained from base ones by virtue of categorical axioms.

\head {\bf 2. $(m,n)$-invariants} \endhead 

\vskip 0.2cm
{\bf Definition 2.1.} 
\item{$\bullet $} {\bf Equivalence} \, $x^k\sim y^k$, \, $x^k,y^k\in L^k$, \, $k\ge 0$, is called {\bf of degree $l$}, \, $deg(\sim ):=l$, $l\ge 0$, if all arrows 
representing it starting from order $k+l+1$ and higher are identities and for $l>0$ there is at least one nonidentity arrow on level 
$k+l$. If there is no such \, $l\in \Bbb N$, \, $deg(\sim ):=\infty $. Denote $\sim $ of degree $l$ by $\sim _l$.
\item{$\bullet $} {\bf Pair of equivalent elements} \, $x^k\sim y^k$, \, $k\ge 0$, is called {\bf of degree $l$}, \, $deg(x^k\sim y^k):=l$, $l\ge 0$,  
if the lowest degree of equivalences existing between $x^k$ and $y^k$ is \, $l$.
\item{$\bullet $} {\bf $\infty $-category} $L$ is called {\bf of degree l}, \, $deg(L)=l$, \, $l\ge 0$, if for any pair of 
equivalent objects $a\sim a'$, $a,a'\in L^0$, there exists an equivalence $a\sim _ka'$ of degree $k\le l$ and there exists at least
one pair of equivalent objects from $L$ of degree $l$.
\item{$\bullet $} {\bf Functor} $F:L\to L'$ is called {\bf $(m,n)$-invariant} if $F$ preserves equivalences $\sim $\, ,
$m=deg(L)$, \ $0\le n\le deg(L')$ and $F$ maps every
pair of equivalent objects of degree $\le m$ to a pair of equivalent objects of degree $\le n$, i.e. $deg(a\sim a')\le m\Rightarrow deg(F(a)\sim F(a'))\le n$,
and boundary $n$ is actually achieved on a pair of equivalent objects of $L$.       \hfill        $\square $

\vskip 0.2cm
{\bf Remarks.} 
\item{$\bullet $} $(m,n)$-invariants are important for the classification problem (up to $\sim $). If $n<m$ an $(m,n)$-invariant decreases
complexity of the equivalence relation, i.e. partially resolves it.   
\item{$\bullet $} There can be trivial invariants which do not distinguish anything and do not carry any information 
such as constant functors $c:L\to L'$ (although they are $(deg(L),0)$-invariants).  \vskip -0.35cm   \hfill   $\square $

\vskip 0.3cm
\centerline{\bf Examples}
\vskip 0.25cm
\item{1.} $deg(ea)=0$; \, $deg(f:a@>\sim >iso>a')=1$; \, $deg(\bold{Set})=1$; \, $deg(\infty \text{-}\bold{Top})=2$; \, $deg(\infty \text{-}\bold{CAT})=\infty \, (?)$. 
\vskip 0.1cm
\item{2.} Homology and cohomology functors $H_*,H^*:\infty \text{-}\bold{Top}\to \bold{Ab}$ (trivially extended over higher order cells) are $(2,1)$-invariants.
\vskip 0.1cm
\item{3.} $\pi ^I_n/\hskip -0.1cm\sim \hskip 0.05cm:L^*_{1\sim }\to \bold{Grp}$ is an $(\infty ,1)$-invariant (see proposition 2.1.2).
\vskip 0.1cm
\item{4.} Let $X$ be a smooth manifold with Lie group action $\rho :G\times X\to X$, \, $L$ be a category with $L^0$, the set of 
submanifolds of $X$, $L^1(a,b):=\{(a,g,b)\in L^0\times G\times L^0\, |\, \rho (g,a)=b\, \}$, \, $L^n:=eL^{n-1}$ for  $n\ge 2$, \, 
$L'$ be a category with $L^{'0}:=C^{\infty }(X,\Bbb R)$ (smooth functions), 
$L^{'1}(f,h):=\{(f,g,h)\in L^{'0}\times G\times L^{'0}\, |\, f\circ \rho (g^{-1},-)=h\, \}$, \, $L^{'n}:=eL^{'n-1}$ for  $n\ge 2$.
If $F:L\to L'$ is a construction (functor) assigning invariant functions to objects from $L$ then $F$ is an $(1,0)$-invariant.  
\item{5.} Each equivalence $L@>\sim >>L'$ is \, $(deg(L),deg(L'))$-invariant with \, $deg(L)=deg(L')$. 

\vskip 0.3cm
\subhead 2.1. Homotopy groups associated to $\infty $-categories \endsubhead

\vskip 0.2cm
Let $L$ be an $\infty $-category in which $*$ strictly preserves $e$ and $\sim $ (i.e. $*$ is a strict functor).
Denote by $eqL:=\{f\in L\ |\ \exists \, g. \ edf\sim g\circ _1f, edg\sim f\circ _1g\}$ subset of eqivalences of $\infty $-category $L$. 
It can be not a category (because it is not closed under $d,c$, in general).
\vskip 0.2cm
{\bf Definition 2.1.1.} Assume, $L(I,-):L\to \infty \text{\bf -CAT}$, \, $x\in L^0(I,a)$. Then $\pi ^I_{n}(a,x):=$\newline
$$\cases (L^0(I,a),x) & \ \ \text{if} \ $n=0$ \\
\bold{Aut}_{L(I,a)}(e^{n-1}x):=eqL(I,a)(e^{n-1}x,e^{n-1}x)\cap (L(I,a))^0(e^{n-1}x,e^{n-1}x)= & \\
=eqL(e^{n-1}x,e^{n-1}x)\cap L^{n+1} & \ \ \text{if} \ n>0
\endcases$$
are (weak) {\bf homotopy groups} of object $a$ at point $x$ with representing object $I\in L^0$.  \hfill $\square $
\vskip 0.2cm
$\pi ^I_{0}(a,x)$ or $\pi ^I_{0}(a,x)/\sim $ are just pointed sets, \ $\pi ^I_{n}(a,x)/\sim , n>0$ are strict groups.
\vskip 0.2cm
{\bf Remark.} If $L=\infty \text{\bf -Top}$, \ $I=\bold 1$ then the above homotopy groups are usual ones.  \hfill  $\square $
\vskip 0.2cm
{\bf Definition 2.1.2.} For a map \ $f:a\to b$ \ such that $f\circ x=y$, \ $x\in L^0(I,a)$, \ $y\in L^0(I,b)$ the {\bf induced map} \ $f_*\equiv \pi ^I_n(f):\pi ^I_{n}(a,x)\to \pi ^I_{n}(b,y)$ \ is determined by restriction of functor $L(I,f):$\newline
$$\cases L^0(I,a)\to L^0(I,b):x'\mapsto f\circ _1x' & \ \ \text{if} \ n=0 \\
\bold{Aut}_{L(I,a)}(e^{n-1}x)\to \bold{Aut}_{L(I,b)}(e^{n-1}y):g\mapsto \mu _{I,a,b}(e^nf,g) & \ \ \text{if} \ n>0
\endcases $$
\vskip -0.6cm
\hfill $\square $

\vskip 0.2cm
{\bf Remark.} To be correctly defined induced maps $\pi ^I_n(f)$ for $n>1$ need commutativity of $*$ with $e$. First two 'groups'
$\pi ^I_0(a,x), \pi ^I_1(a,x)$ always make sense and depend functorially on objects.    \hfill  $\square $

\vskip 0.2cm
\proclaim{\bf {Proposition 2.1.1 (homotopy invariance of homotopy groups)}} If $x:I\to a$, $f\sim f'\in L^0(a,b)$ such that
$f\circ _1x\sim f'\circ _1x$ is trivial equivalence (all arrows for $\sim $ are identities) then 
$\pi ^I_n(f)/\sim \, =\pi ^I_n(f')/\sim \, :\pi^I_n(a,x)/\sim \, \to \pi ^I_n(b,f\circ x)/\sim \, $.
\endproclaim
\demo\nofrills{Proof\ \ } is immediate. \hfill $\square $
\enddemo

\proclaim{\bf Proposition 2.1.2} $\pi ^I_n/\hskip -0.1cm\sim \hskip 0.05cm:L^*_{1\sim }\to \bold{Grp}$ is an $(\infty ,1)$-invariant, 
where $L^*_{1\sim }:=\coprod\limits _{n\ge 0}L^{*n}_{1\sim }$, \, $L^{*n}_{1\sim }:=\cases  L^{*n}\text{ (pointed objects and maps) } & n=0,1  \\
\text{equivalences from }L^{n} \text{ with dom and codom in }L^{*(n-1)}_{1\sim }    & n>1     \endcases $. 
\endproclaim 
\demo{Proof} Partial functor $\pi ^I_{n}/\hskip -0.1cm\sim \hskip 0.05cm:L^{*0}\coprod L^{*1}\to \bold{Grp}$ is trivially extendable starting from 
equivalences on level 2 (because of proposition 2.1.1).    \hfill   $\square $
\enddemo 

\vskip 0.1cm
\centerline{{\bf Example} (Fundamental Group)}
\vskip 0.2cm
Let $2\text{-}\bold{Top}$ be usual $\bold{Top}$ with homotopy classes of homotopies as $2$-cells. 
Define {\bf fundamental groupoid} $2$-functor as representable 
\, $\Pi (-):=Hom_{2\text{-}\bold{Top}}(1,-):2\text{-}\bold{Top}\to 2\text{-}\bold{Cat}:$\vskip 0.2cm \hskip -0.39cm
$\cases  X\to \Pi (X) & Ob\, (\Pi (X)) \text{ are its points, \, } Ar\, (\Pi (X)) \text{ are homotopy classes of pathes} \\
(X@>f>>Y)\mapsto \Pi (f) & \text{transformation of fundamental groupoids, \ } \Pi (f):\cases x\mapsto f(x) & \\ 
[\gamma ]\mapsto [f\circ \gamma ] & \endcases \\
(f@>[H]>>f')\mapsto \Pi ([H]) &  \text{nat. trans. } \Pi ([H])=\underset {{x\in X}}\to {\{[H]*i_x\}}:Hom_{2\text{-}\bold{Top}}(1,f)@>\sim>>Hom_{2\text{-}\bold{Top}}(1,f')  \hskip -0.2cm  \endcases $ 
\vskip 0.0cm
(where \, $\underset {{x\in X}}\to {\{[H]*i_x\}}=\underset {x\in X}\to {\{[H(x,-)]\}}$ \, are homotopy classes of pathes between $f(x)$ and $f'(x)$ natural in $x\in X$).
\vskip 0.15cm\hskip -0.39cm 
$\pi _1(X,x_0):=\bold{Aut}_{\Pi (X)}(x_0)\hookrightarrow \Pi (X)$ is {\bf fundamental group} of space $X$ at point $x_0\in X$, 
$\pi _1((X,x_0)@>f>>(Y,y_0)):=\bold{Aut}_{\Pi (X)}(x_0)@>\Pi (f)>>\bold{Aut}_{\Pi (Y)}(y_0)$.

\proclaim{\bf Proposition 2.1.3} \vskip 0.05cm
\item{$\bullet $} If \, $[H]:f@>\sim >>f':X\to Y$ \, is a $2$-cell in \, $2\text{-}\bold{Top}$ \, then \, $\pi _1(f')([\gamma ])=[H(x_0,-)]\circ \pi _1(f)([\gamma ])\circ [H(x_0,-)]^{-1}$, \, $[\gamma ]\in \pi _1(X,x_0)$.
\item{$\bullet $} In the case \, $[H]:f@>\sim >>f':(X,x_0)\to (Y,y_0)$ \, is a pointed $2$-cell \, (\, $[H(x_0,-)]=1_{f(x_0)}:f(x_0)\to f(x_0)=f'(x_0)$) \, then \, $\pi _1(f)=\pi _1(f')$. 
\endproclaim
\demo\nofrills{Proof \ \ } follows from the naturality square\hskip 0.7cm
$\vcenter{\xymatrix{{f(x_0) \ } \ar[r]_-{\sim }^-{[H(x_0,-)]} \ar[d]_-{\Pi (f)([\gamma ])} & { \ f'(x_0)} \ar[d]^-{\Pi (f')([\gamma ])} \\
{f(x_0) \ } \ar[r]^-{\sim }_{[H(x_0,-)]} & { \ f'(x_0)}}}$\vskip -0.55cm
\hfill $\square $
\enddemo

\subhead 2.2. Duality and Invariant Theory\endsubhead

\vskip 0.2cm
{\bf Proposition 2.2.1.} Let $\bold{K}$ be $\bold{Set}$, $\bold{Top}$ or $\bold{Diff}^+$ (spectra of smooth completion of commutative algebras with Zariski topology), $G$ be a group. Then there exists
a concrete natural dual adjunction $\xymatrix{\bold{ComAlg}^{op} \ar@/^/[r]^-F & G\text{-}\bold{K} \ar@/^/[l]^-H_-{\top }}$ with $k$ 
($\Bbb R$ or $\Bbb C$), its schizophrenic object, such that $k\in Ob\, G\text{-}\bold{K}$ has trivial action of $G$, and 
$F\circ H:G\text{-}\bold{K}\to G\text{-}\bold{K}$ is a functor 'taking the factor-space generated by equivalence relation $x\sim y$ iff 
$x,y \in \text{Closure}(\text{the same orbit})$' (it is essentially the orbit space).  \hfill $\square $ 

\vskip 0.2cm
{\bf Definition 2.2.1.} \item{$\bullet $} Adjoint object $\Cal A_X=HX$ for an object $X$ in $G\text{-}\bold{K}$ is called its 
{\bf algebra of invariants}. 
\item{$\bullet $} If $U:G\text{-}\bold{K}\to G\text{-}\bold{K}$ is an endofunctor then $\Cal A_{U(X)}$ is called an {\bf algebra of 
$U$-invariants} of object $X$. \hfill $\square $ 

\vskip 0.2cm
{\bf Remarks.}
\item{$\bullet $} For $U=(-)^n$, n-fold Cartesian product, $\Cal A_{U(X)}$ is an {\bf n points' invariants' algebra}.
\item{$\bullet $} For $\bold{K}=\bold{Diff}$, $U=\bold{Jet}^n$, $\bold{Jet}^n(X):=\{j^n_0f \, |\, f\in \bold{Diff}(k,X)\}$, set of all 
$n$-jets of all maps from $k$ to $X$ at point $0$ (with a certain manifold structure obtained from local trivializations), we get {\bf differential invariants}. 
\item{$\bullet $} $U=\bold{Jet}^{\infty }:\bold{Diff}\to \bold{Diff}^+$ does not fit to the above scheme, but everything is still correct if 
$U:G\text{-}\bold{K}\to G\text{-}\bold{K}_1$ is an extension to $G\text{-}\bold{K}_1$, a category concretely adjoint to $\bold{ComAlg}$.
\item{$\bullet $} $G$ can be, of course, $\bold{Aut}(X)$. 

\vskip 0.2cm
Accordingly to Klein's Erlangen Program every group acting on a space determines a geometry and, conversely, every geometry hides 
a group of transformations. Properties of geometric objects which do not change under all transformations are called geometric 
(or invariant, or absolute for the given $G$-space and a class of geometric objects).

\vskip 0.2cm
{\bf Equivalence problem} \cite{Car1, Car2, Vas, Olv, Gar} consists of $G$-space $X$ and two 'geometric objects' $S_1, S_2$ of the same type on space $X$. It is required
to determine if these two objects can be maped to one another by an element of $G$. An approach is to find a (complete) system of 
invariants of each object. 
 
\vskip 0.2cm
\subhead {} 2.2.1. Classification of covariant geometric objects\endsubhead

Under covariant geometric objects we mean objects like submanifold, foliation or system of differential equations, i.e., objects which
behave contravariantly from Categorical viewpoint and which can be described by a {\bf diffential ideal} $I$ ($dI\subset I$) in 
$\Lambda (X)$, exterior differential algebra of $X$.

\proclaim{\bf Proposition 2.2.1.1} Let $G$ be a Lie-like group (i.e., there exists an algebra of invariant forms on $G$). Then any $G$-equivariant map $\sigma :G\to X$ ($G$ is given with left shift action and $X$ is a left $G$-space) produces a system of invariants of differential 
ideal $I\subset \Lambda (X)$ (with generators of degree $0$ and $1$) in the folloing way:
\item{$\bullet $} Take the image $\bar {\Lambda }_{inv}:=Im\, (\Lambda _{inv}(G)\hookrightarrow \Lambda (G)\twoheadrightarrow \Lambda (G)/{\sigma ^*}(I))$,
where $\Lambda _{inv}(G)$ is a subalgebra of left-invariant forms on $G$, $\sigma ^*:\Lambda (X)\to \Lambda (G)$ is the induced map of 
exterior differential algebras, $\sigma ^*(I)$ is the smallest differential ideal in $\Lambda (G)$ containing image of $I$ under $\sigma ^*$.
\item{$\bullet $} Take module $\Lambda ^0(G)\cdot \bar {\Lambda }^1_{inv}$ generated by $1$-forms in $\bar {\Lambda }_{inv}$ over $\Lambda ^0(G)$.
There is an open set $\Cal O\subset G$ and a basis $\{\omega ^{\alpha }_{inv}\}_{\alpha \in A}\subset \bar {\Lambda }^1_{inv}$ for module 
$\Lambda ^0(G)\cdot \bar {\Lambda }^1_{inv}$ restricted on $\Cal O$, i.e., $\forall \, \omega ^i_{inv}\in \bar {\Lambda }^1_{inv}$ 
$\exists \, ! \, \text{ functions } f^i_{\alpha }\in C^{\infty }(\Cal O)$ such that $\omega ^i_{inv}=\sum\limits _{\alpha }f^i_{\alpha }\omega ^{\alpha }_{inv}$.
Form set $J_0:=\{f^i_{\alpha }\}$.
\item{$\bullet $} Take expansion of differentials $df^i_{\alpha }=\sum\limits _{\beta }f^i_{\alpha \beta }\omega ^{\beta }_{inv}$ (over $\Cal O$). Form set 
$J_1:=\{f^i_{\alpha \beta }\}$. 
\item{$\bullet $} Continue this process to get $J_2:=\{f^i_{\alpha \beta \gamma }\}, \dots , J_n:=\{f^i_{\alpha _1\dots \alpha _{n+1}}\}\dots $
Form set $J:=\bigcup\limits _{n}J_n$. Its elements are relative invariants of differential ideal $I\subset \Lambda (X)$.
\item{$\bullet $} Take algebra $\Cal A_J\subset C^{\infty }(\Cal O)$, generated by $J$, and take its smooth completion $\overline {\Cal A_J}$ (see {\bf 3.4}). 
Then ideal $\xymatrix{\bold{Rel}(\Cal A_J) \ar@{>->}[r] & \overline {\bold{Alg}(J)} \ar@{->>}[r] & \overline {\Cal A_J}}$, of all relations of $\Cal A_J$, gives absolute invariants of differential ideal $I\subset \Lambda (X)$, 
where $\overline {\bold{Alg}(J)}$ is the smooth completion of free algebra generated by $J$.
\endproclaim   
\demo\nofrills{Proof} \ \ follows from the diagrams

\hbox{\hskip 2cm
\hbox{$\xymatrix{G \ar[r]^-{l_g} \ar[d]_-{\sigma } & G \ar[d]^-{\sigma } \\
X \ar[r]_-{l_g} & X}$}
\hskip 2cm
\hbox{$\xymatrix{\Lambda _{inv}(G) & \Lambda _{inv}(G) \ar[l]_-{id} \\
\Lambda (X) \ar[u]^-{\sigma ^*} & \Lambda (X) \ar[l]^-{l_g^*} \ar[u]_-{\sigma ^*}}$}
}

and equations $\omega ^i_{inv}=\sum\limits _{\alpha }f^i_{\alpha }\omega ^{\alpha }_{inv}\, mod(\sigma ^*(I))$.
\hfill $\square $
\enddemo

{\bf Remark.} $G\text{-}\bold{Diff}(G,X)$ is in $1$-$1$-correspondence with all sections of orbit space $X_G$. So, if $X$ is homogenious then 
it is exactly the set of all points of $X$ and $\sigma :G\to X=G@>\sim >>G\times \{x_0\}@>1\times i_{x_0}>>G\times X@>\rho >>X$ is a $G$-equivariant 
map corresponding to point $x_0\in X$, where $\rho $ is the given $G$-action on $X$. 

\vskip 0.1cm
\proclaim{\bf Proposition 2.2.1.2 (Exterior differential algebra associated to a group of analityc automorphisms)} Let $X$ be an analytic $n$-dimensional manifold, $\bold{An}(X)$, its group of automorphisms, $H^{\infty }(X):=\{j^{\infty }_0f\, |\, f\in \bold{Diff}(k^n,X), \ X \text{ is analytic}, \text{ Jacobian}(f)\ne 0 \}$, 
 $\infty $-frame bundle over $X$ (with a usual topology and manifold structure). Then there is an exterior differential $k$-algebra $\Lambda _{inv}(H^{\infty }(X))$
of invariant forms on $H^{\infty }(X)$ freely generated by elements of degree 1 obtained by the following process:
\item{$\bullet $} $\omega ^i:=x^i_jdx^j$ are any 'shift' forms on $X$
\item{$\bullet $} $\omega ^i_{j}$ are most general solutions of Maurer-Cartan equations $d\omega ^i=\omega ^i_j\wedge \omega ^j$
\item{$\bullet $} $\omega ^i_{jk}$ are most general solutions of Maurer-Cartan equations $d\omega ^i_j=\omega ^i_k\wedge \omega ^k_j+\omega ^i_{jk}\wedge \omega ^k$
\item{$\bullet $} $\omega ^i_{jkl}$, $\cdots $, $\omega ^i_{i_1\dots i_n}$, $\cdots $

All forms are symmetric in lower indices. They characterize underlying space of $\bold{An}(X)$ uniquely up to analytic iso. \hfill $\square $
\endproclaim   

\vskip 0.0cm
{\bf Remark.} At each point $x_0\in X$, \, $\omega ^i=0$, \, and forms \, $\bar \omega ^i_{i_1\dots i_n}:=\bigl .\omega ^i_{i_1\dots i_n}\bigr |_{\omega ^i=0}\, , \ \, n\ge 1$, \, are free generators of exterior differential 
algebra of {\bf differential group} acting simply transitively on each fiber of $H^{\infty }(X)$.

\vskip 0.2cm
\subhead 2.2.2. Classification of smooth embeddings into Lie group\endsubhead

It is often the last step of smooth classification of geometric objects \cite{Car2, Fin, Kob}. Process of finding of differential invariants is similar to that in Proposition
2.2.1.1.

\proclaim{\bf Proposition 2.2.2.1} For a smooth embedding $f:X\to G$ of smooth manifold $X$ into Lie group $G$ a complete system of 
differential invariants of $f$ can be obtained in the following way:
\vskip 0.1cm
\item{$\bullet $} $Im(f^*:\Lambda ^1_{inv}(G)\to \Lambda (X))$ is locally free, so, take its basis $\omega _{inv}^i, i=1,\dots ,n$, $n=dim(X)$, near each point.
\item{$\bullet $} Coefficients of linear combinations $\omega ^I_{inv}=\sum \limits _{i=1}^na^I_i\omega ^i_{inv},\, I=n+1,\dots ,dim(G)$, are differential invariants 
of first order (of map $f$).
\item{$\bullet $} Coefficients of differentials of invariants of first order $da^I_i=\sum \limits _{j=1}^na^I_{ij}\omega _{inv}^j$ are 
differential invariants of second order (of map $f$).
\item{$\bullet $} $\dots $ Coefficients of differentials of invariants of $(k-1)$ order $da^I_{i_1\dots i_{k-1}}=\sum \limits _{i_k=1}^na^I_{i_1\dots i_k}\omega _{inv}^{i_k}$ 
are differential invariants of order $k$ $\dots $
\vskip 0.2cm
Such calculated invariants characterize orbit $G\cdot f$ uniquely up to 'changing parameter space' $X@>\sim >>X'$.
\endproclaim 
\demo\nofrills{Proof\ \ } is straightforward. \hfill $\square $
\enddemo

\head {\bf 3. Representable $\infty $-functors} \endhead

{\bf Definition 3.1.} $\infty $-categories $L$ and $L'$ are {\bf equivalent} if $L\sim L'$ in {\bf $\infty $-CAT}. \hfill $\square $

If equivalence $L\sim L'$ is given by functors $\xymatrix{L \ar@/^2ex/[r]^-{F}_-{\sim } & L' \ar@/^2ex/[l]^-{G}_-{\sim }}$ then  
$\forall a\in L^0 \  a\sim G\circ F(a)$,  $\forall b\in L^{'0}$\linebreak
$b\sim F\circ G(b)$  naturally in $a$ and $b$.

\vskip 0.2cm
{\bf Definition 3.2.} $\infty $-functor $F:L\to L'$ is (weakly)
\item{$\bullet $} {\bf faithful} if $\forall a,a'\in L^0$ $\forall f^n,g^n\in L^n(a,a')$ \, $F(f^n)\sim F(g^n)\Rightarrow f^n\sim g^n$, 
\item{$\bullet $} {\bf full} if $\forall a,a'\in L^0$ $\forall h^n\in L^{'n}(F(a),F(a'))$ \, $\exists f^n\in L^n(a,a')$ such that $F(f^n)\sim h^n$,
\item{$\bullet $} {\bf surjective on objects} if $\forall b\in L^{'0}$ \ $\exists a\in L^0$ such that $F(a)\sim b$.
\hfill $\square $
\vskip 0.2cm
Unlike first order equivalence there is no simple criterion of higher order equivalence. 

\proclaim{\bf {Proposition 3.1}} If functor $\xymatrix{L \ar[r]^-{F}_-{\sim } & L'}$ is an equivalence then $F$ is (weakly) faithful full and 
surjective on objects.
\endproclaim
\demo{Proof} "$\Rightarrow $" \ Regard the diagram
$$\xymatrix{a \ar@/^2ex/[rr]^-{e^n\rho _a}_-{\sim } \ar[d]_-{f^n} && G\circ F(a) \ar@/^2ex/[ll]^-{e^n\theta _a}_-{\sim } \ar[d]^-{G(F(f^n))} \\
a' \ar@/^2ex/[rr]^-{e^n\rho _{a'}}_-{\sim } && G\circ F(a') \ar@/^2ex/[ll]^-{e^n\theta _{a'}}_-{\sim }
}$$

where: $f^n\in L^n(a,a')$, $e^n\rho _{a}\in L^n(a,G(F(a)))$, $e^n\theta _a\in L^n(G(F(a)),a)$, $n\ge 0$.
\vskip 0.2cm
Take $f^n,g^n:a\to a'\in L^n(a,a')$ such that $F(f^n)\sim F(g^n)$. Then $f^n\sim e^n\theta _{a'}\circ _{n+1}G(F(f^n))\circ _{n+1}e^n\rho _a\sim e^n\theta _{a'}\circ _{n+1}G(F(g^n))\circ _{n+1}e^n\rho _a\sim g^n$, \, i.e., $F$ is faithful ($G$ is faithful by symmetry).
\vskip 0.1cm
Take $\alpha ^n:F(a)\to F(a')\in L^{'n}(F(a),F(a'))$. Then $\beta ^n:=e^n\theta _{a'}\circ _{n+1}G(\alpha ^n)\circ _{n+1}e^n\rho _a:a\to a'\in L^n(a,a')$ is such that $G(F(\beta ^n))\sim G(\alpha ^n)$. So, $F(\beta ^n)\sim \alpha ^n$ because $G$ is faithful. Therefore, $F$ is full ($G$ is full by symmetry).
\vskip 0.1cm
$F$ and $G$ are obviously surjective on objects. \hfill $\square $
\vskip 0.2cm
{\bf Remark.} The inverse direction "$\Leftarrow $" for the above proposition works only partially. Namely, 
for each $b\in L^{'0}$ choose $G(b)\in L^0$ and equivalence 
$\xymatrix{b \ar@/^2ex/[rr]^-{\rho _b}_-{\sim } && F(G(b)) \ar@/^2ex/[ll]_-{\sim }^-{\theta _b}}$ 
(which is possible since $F$ is surjective on objects), moreover, if 
$b=F(a)$ choose $G(b)=a, \ \rho _b=eb, \ \theta _b=e(F(G(b)))=eb$. For each $f^n:b\to b'\in L^{'n}(b,b')$ choose an element 
$G(f^n)\in  L^n(G(b),G(b'))$ such that $e^n\rho _{b'}\circ _{n+1}f^n\circ _{n+1}e^n\theta _b\sim F(G(f^n))$ 
(which is possible since $F$ is fully faithful). Then $G:L'\to L$ is obviously a (weak) functor. $a=G(F(a))$ is natural in $a$ 
by construction, but $b\sim F(G(b))$ is natural in $b$ for only first order arrows $\rho _b, \, \theta _b$ presenting $\sim $\, . So, $F$ 
should be somehow 'naturally surjective on objects' which does not make sense yet when functor $G$ is not defined.
\enddemo

\vskip 0.2cm
{\bf Definition 3.3.} $\infty $-functor $F:L\to L'$ is called
\item{$\bullet $} {\bf isomorphism} if it is a bijection (on sets $L,L'$) and the inverse map is a functor,
\item{$\bullet $} {\bf quasiisomorphism} if there exists a functor $G:L'\to L$ such that 
$\forall a^n\in L^n \ \, G(F(a^n))\sim a^n$ \, and \, $\forall b^n\in L^{'n} \ \, F(G(b^n))\sim b^n$, \, $n\ge 0$.  \hfill $\square $

\proclaim{\bf Proposition 3.2} Notions of (functor) isomorphism and quasiisomorphism coincide.
\endproclaim
\demo{Proof} Each isomorphism is a quasiisomorphism. Conversely, if $\xymatrix{L \ar@/^/[r]^-{F} & L' \ar@/^/[l]^-{G} }$ is 
a quasiisomorphism then $\forall a^n\in L^n$, \, $n\ge 0$, \, $G(F(ea^n))\sim ea^n$. So, $d(G(F(ea^n)))=dea^n$, i.e. $G(F(dea^n))=dea^n$ and $G(F(a^n))=a^n$ (instead of $d$, $c$ could be used). 
The same, $\forall b^n\in L^{'n}$, \, $n\ge 0$, \, $F(G(b^n))=b^n$.          \hfill $\square $
\enddemo

\vskip 0.1cm
Denote (quasi)isomorphism (equivalence) relation by \, $\simeq $\, . 

\vskip 0.35cm
\centerline{{\bf Examples} (isomorphic $\infty $-categories)}
\vskip 0.1cm
\item{1.} Assume, $\xymatrix{f^n \ar[r]^-{\simeq }_-{\alpha } & g^n }$ are 
isomorphic elements of degree $n$ (in a strict category $L$) then $L(f^n,f^n)\simeq L(g^n,g^n)$ are isomorphic $\infty $-categories.
Indeed, there is an isomorphism $F:L(f^n,f^n)\to L(g^n,g^n):x\mapsto \alpha *(x*\alpha ^{-1})$,  
where $*$ means a horizontal composite. $F$ is a functor. Its inverse is 
$G:L(g^n,g^n)\to L(f^n,f^n):y\mapsto \alpha  ^{-1}*(y*\alpha )$. 
[For $\alpha $, just equivalence, it is not true] 
\item{2.} $\infty\text{\bf -CAT}(L(-,a),F)\simeq F(a)$ (see below Yoneda Lemma).    \hfill $\square $

\vskip 0.2cm
{\bf Definition 3.4.} Two $n$-modifications $\alpha ^n, \beta ^n:L\to \infty \text{-}\bold{CAT}$, $n\ge 0$, are called 
{\bf quasiequivalent} of deepness $k$, \, $0\le k\le n+1$, (denote it by $\alpha ^n\approx _k\beta ^n$) if their corresponding components 
are quasiequivalent of deepness $k-1$, \, i.e. $\forall a\in L^0$ \, $\alpha ^n(a)\approx _{k-1}\beta ^n(a)$. 
$\approx _0$ means $\sim $ by definition. [In other words, $\alpha ^n\approx _k\beta ^n$ if all their components 
of components on deepness $k$ are equivalent, i.e. $\alpha ^n\approx _0\beta ^n$ if they are equivalent $\alpha ^n\sim \beta ^n$;
$\alpha ^n\approx _1\beta ^n$ if their components are equivalent $\forall \, a\in L^0 \ \alpha ^n(a)\sim \beta ^n(a)$; 
$\alpha ^n\approx _2\beta ^n$ if components of all components are equivalent; etc.]. 
If $\alpha ^n, \beta ^n:L\to L'$ are proper $n$-modifications (living in $\infty \text{-}\bold{CAT}$) for them only
$\approx _0$ and $\approx _1$ make sense.        \hfill   $\square $

\vskip 0.2cm
\proclaim{\bf Lemma 3.1} 
\item{$\bullet $} $\approx _k$ is an equivalence relation.
\item{$\bullet $} $\approx _{k_1}\ \Rightarrow \ \approx_{k_2}$ if $k_1\le k_2$.
\item{$\bullet $} If $\alpha ^n\approx _k\beta ^n$ then $d\alpha ^n=d\beta ^n$, \, $c\alpha ^n=c\beta ^n$.
\item{$\bullet $} If $(L_1,\approx _{k_1})$, $(L_2,\approx_{k_2})$ are two $\infty $-categories 
(not necessarily proper, i.e. living in $\infty \text{-}\bold{CAT}$) for which given equivalence relations make sense for 
all elements, and $F:L_1\to L_2$, $G:L_2\to L_1$ are maps (not necessarily functors) such that $\forall \, l_1\in L_1$ 
$G(F(l_1))\approx _{k_1}l_1$ and $\forall \, l_2\in L_2$ $F(G(l_2))\approx _{k_2}l_2$, \, and $F,G$ both preserve $d$ (or $c$) then 
$F,G$ are bijections inverse to each other.
\item{$\bullet $} For $L,L'\in Ob\, (\infty \text{-}\bold{CAT})$ and $a\in L^0$ the map 
$ev_{a}:\infty \text{-}\bold{CAT}(L,L')\to L':f^n\mapsto f^n(a)$ is a strict functor. [Similar statement holds when $L,L'$ are 
not proper, e.g. $\infty \text{-}\bold{CAT}$, but we need to formulate it a biger universe containing $\infty \text{-}\bold{CAT}$] 
\endproclaim
\demo{Proof} First two statements are obvious. Third one follows from the fact $x\sim y \, \Rightarrow \, dx=dy, \ cx=cy$ \, and that 
$d,c$ are taken componentwise. Forth statement follows by the same argument as in the proof of proposition 1.3.2. The last 
statement holds because, again, all operations in $\infty \text{-}\bold{CAT}(L,L')$ are taken componentwise.  \hfill  $\square $
\enddemo 

\vskip 0.0cm
{\bf Remark.} For the proof of Yoneda lemma a double evaluation functor is needed. For two functors $F,G:L\to \infty \text{-}\bold{CAT}$ 
take the restriction of evaluation functor $ev_a$ on the hom-set between $F$ and $G$, i.e. 
$ev_{a\, F,G}:\infty \text{-}\Bbb {CAT}(L,\infty \text{-}\bold{CAT})(F,G)\to \infty \text{-}\bold{CAT}(F(a),G(a)):f^n\mapsto f^n(a)$, 
where $\infty \text{-}\Bbb {CAT}$ is a bigger (and weaker) universe containing $\infty \text{-}\bold{CAT}$ as an object. Now, take a second 
evaluation functor $ev_x:\infty \text{-}\bold{CAT}(F(a),G(a))\to G(a):g^n\mapsto g^n(x)$, \, $x\in (F(a))^0$. Then the double 
evaluation functor is the composite $ev_x\circ _1ev_{a\, F,G}:\infty \text{-}\Bbb {CAT}(L,\infty \text{-}\bold{CAT})(F,G)\to G(a):f^n\mapsto f^n(a)(x)$.
It is a strict functor.      \hfill    $\square $

\vskip 0.2cm
$\infty \text{-}\bold{CAT}$-valued functors, natural transformations and modifications live now in a bigger universe 
$\infty \text{-}\Bbb {CAT}$, and we do not have for them appropriate definitions, yet. 

\vskip 0.2cm
{\bf Definition 3.5.} $\infty \text{-}\bold{CAT}$-valued functors, natural transformations and modifications are introduced in a 
similar way as usual ones with changing all occurrences of $\sim $ with (one degree weaker relation) $\approx _1$, i.e.
\item{$\bullet $} a map $F:L\to \infty \text{-}\bold{CAT}$ of degree $0$ is a {\bf functor} if $F$ strictly preseves $d$ and $c$,
$Fdx=dFx, \ Fcx=cFx$, and weakly up to $\approx _1$ preserves $e$ and composites, $Fex\approx _1eFx$, $F(x\circ _ky)\approx _1F(x)\circ _kF(y)$, 
\item{$\bullet $} For a given sequence of two functors $F,G:L\to \infty \text{-}\bold{CAT}$, $\dots $, two $(n-1)$-modifications 
$\alpha ^{n-1}_1, \alpha ^{n-1}_2:\alpha ^{n-2}_1\to \alpha ^{n-2}_2$ strict (or weak) {\bf $n$-modification} 
$\alpha ^{n}:\alpha ^{n-1}_1\to \alpha ^{n-1}_2$ is a map $\alpha ^{n}:L^0\to \infty \text{-}\bold{CAT}^{n+1}$ such that
$\forall a,b\in L^0$ $\alpha ^n(b)*F(-)\approx _1G(-)*\alpha ^n(a):L^{\ge n}(a,b)\to L^{'\ge n}(F(a),G(b))$ (components of values 
of functors are equivalent).      \hfill   $\square $

\vskip 0.2cm
{\bf Definition 3.6.} Covariant (contravariant) functor $F:L\to \infty \text{\bf -CAT}$ is 
\item{$\bullet $} {\bf weakly representable} if $\exists a\in L^0$
such that $L(a,-)\sim F$ \, ($L(-,a)\sim F$). It means there is an equivalence of two $\infty $-categories $L(a,b)\sim F(b)$ \, ($L(b,a)\sim F(b)$) natural in $b$,
\item{$\bullet $} {\bf strictly representable} if $\exists a\in L^0$ such that $L(a,-)\simeq F$ \, ($L(-,a)\simeq F$), 
i.e. $\forall b\in L^0$ \, $\exists $ isomorphism \, $L(a,b)\simeq F(b) \ \, (L(b,a)\simeq F(b))$ natural in $b$.   \hfill  $\square $

\vskip 0.25cm
\proclaim{\bf Lemma 3.2} For given representable $L(-,a):L^{op}\to \infty \text{-}\bold{CAT}$ and functor 
$F:L^{op}\to \infty \text{-}\bold{CAT}$
\item{$\bullet $} all natural transformations  $\tau ^0:L(-,a)\to F$  are of the form  $\forall \, b\in Ob\, L$ $b$-component is a functor $\tau ^0_b:L(b,a)\to F(b)$,  $\tau ^0_b(f^m)\sim F(f^m)(\tau ^0_a(ea))$,  $f^m\in L^m(b,a)$,
\item{$\bullet $} all $n$-modifications  $\tau ^n:L(-,a)\to F$, $n\ge 1$, are of the form  $\forall \, b\in Ob\, L$ $b$-component is a
$(n-1)$-modification $\tau ^n_b:L(b,a)\to F(b)$, $\tau ^n_b(f^0)\sim F(f^0)(\tau ^n_a(ea))$, $f^0\in L^0(b,a)$.
\endproclaim
\demo\nofrills{Proof \ } follows from the naturality square \hskip 0.35cm$\vcenter{\xymatrix{a & L(a,a) \ar[r]^-{\tau ^n_a} \ar[d]_-{L(f^m,a)} & F(a) \ar[d]^-{F(f^m)} \\
b \ar[u]^-{f^m} & L(b,a) \ar[r]_-{\tau ^n_b} & F(b)}}$\hskip 0.35cm$n\ge 0$   \hfill  $\square $
\enddemo

\proclaim{\bf Lemma 3.3} For a given $n$-cell $\beta ^n\in (F(a))^n$, $n\ge 0$, $n$-modification $\tau ^n:L(-,a)\to F$ such that
$\tau ^n_a(ea)=\beta ^n$ exists and unique up to \, $\approx _2$\, .
\endproclaim
\demo{Proof} Uniqueness follows from lemma 3.2, existence from the definition of $n$-modification \linebreak 
$\tau ^n_b(f^m):=F(f^m)(\beta ^n)$ (for $n>0$, $m=0$ only) and naturality square showing correctness of the definition \hskip 1cm
$\vcenter{\xymatrix{b & L(b,a) \ar[d]_-{L(g^k,a)} \ar[r]^-{\tau ^n_b} & F(b) \ar[d]^-{F(g^k)}\\
c \ar[u]^-{g^k} & L(c,a) \ar[r]_-{\tau ^n_c} & F(c)}}$  \vskip 0.0cm
\hskip 0.0cm$\bigl (\mu _{c,b,a}(f^m,g^k):=\mu _{c,b,a}(e^{max(m,k)-m}f^m,e^{max(m,k)-k}g^k)\bigr )$                   \hfill $\square $
\enddemo

\vskip 0.2cm
{\bf Corollary 1.} All $n$-modifications $\tau ^n:L(-,a)\to F$, $n\ge 0$, have strict form 
$\tau ^n_b(f^0)=F(f^0)(\tau ^n_a(ea))$, \, $f^0\in L^0(b,a)$.    \hfill   $\square $

\vskip 0.2cm
{\bf Corollary 2 (criterion of representability).} $\infty \text{-}\bold{CAT}$-valued presheaf $F\hskip -0.035cm:\hskip -0.035cmL^{op}\hskip -0.035cm\to \hskip -0.035cm\infty \text{-}\bold{CAT}$ is
\vskip 0.05cm 
\item{$\bullet $} {\bf strictly representable} (with representing object $a\in L^0$) \, iff \, 
there exists an object \, $\beta ^0\in (F(a))^0$ \, such that \, $\forall \, \gamma ^n\in (F(b))^n, \ n\ge 0$, \, $\exists !$ 
%a unique up to $\sim $ 
$n$-arrow \, $(f^n:b\to a)\in L^n(b,a)$ \, with \, $\gamma ^n=F(f^n)(\beta ^0)$,
\vskip 0.05cm 
\item{$\bullet $} {\bf weakly representable} (with representing object $a\in L^0$) \, iff \,
there exists an object \, $\beta ^0\in (F(a))^0$ \, such that \, $\forall b\in Ob\, L$ \, the functor \, 
$L(b,a)\to F(b):f^n\mapsto F(f^n)(\beta ^0)$ is an equivalence of categories.
%natural in $b\in Ob\, L$.
\item{}(Similar statements hold for covariant presheaf \, $F:L\to \infty \text{-}\bold{CAT}$)
\hfill $\square $ 

\vskip 0.2cm

\proclaim{\bf {Proposition 3.3 (Yoneda Lemma)}} 
For functor $F:L^{op}\to \infty\text{\bf -CAT}$ and object $a\in L^0$ there is a strict isomorphism 
$\infty \text{-}\hskip 0.0125cm\Bbb {CAT}(L(-,a),F)\simeq F(a)$ natural in $a$ and $F$.
\endproclaim
\demo{Proof} Strict functoriality of the correspondence $\tau ^n\mapsto \tau ^n_a(ea)$ is straightforward (because it is a 
double evaluation functor). The map $\beta ^n\mapsto F(-)(\beta ^n)$ is quasiinverse to the first map (with respect to $\approx _2$ and $=$ 
equivalence relations in $\infty\text{-}\Bbb {CAT}(L(-,a),F)$ and $F(a)$ respectively), and it strictly preserves $d$ and $c$. So, these both 
maps are strict isomorphisms. 

Naturality is given by \hskip 0.5cm$\vcenter{\xymatrix{a & F \ar[d]_-{\alpha ^k} & & \infty \text{-}\Bbb {CAT}(L(-,a),F) \ar[r]^-{\simeq } \ar[d]_-{\infty \text{-}\Bbb {CAT}(L(-,f^m),\alpha ^k)} & F(a) \ar[d]^-{\alpha ^k(f^m)} \\
b \ar[u]^-{f^m} & G & & \infty \text{-}\Bbb {CAT}(L(-,b),G) \ar[r]_-{\simeq } & G(b)}}$  \vskip 0.0cm
$\bigl (\text{where \, } \alpha ^k(f^m):=\mu _{F(a),F(b),G(b)}(e^{max(k,m)-k}\alpha ^k_b,e^{max(k,m)-m+1}F(f^m)), \ k,m\ge 0
\bigr )$
\hfill $\square $
\enddemo

\vskip 0.2cm
{\bf Remark.} Yoneda lemma for $\infty $-categories is similar to one for first order categories with the difference that elements 
$\beta ^n\in (F(a))^n$ of degree $n$ determine now higher degree arrows ($n$-modifications) 
$\beta ^n:L(-,a)\to F$ in $\infty \text{-}\bold{CAT}$-valued presheaves category.         \hfill    $\square $

\vskip 0.2cm
\proclaim{\bf Proposition 3.4 (Yoneda embedding)} There is Yoneda embedding $\bold Y\hskip -0.05cm:\hskip -0.05cmL\to \infty\text{\bf -CAT}^{L^{op}}\hskip -0.15cm:\alpha \mapsto L(-,\alpha )$, \ $\alpha \in L$,
which is an extension of the isomorphisms from Yoneda lemma determined on hom-sets $L(a,b)$, $a,b\in L^0$. Yoneda embedding preserves and reflects
equivalences $\sim $\, .      
\endproclaim
\demo{Proof} By Yoneda isomorphism for a given $f^n\in L^n(a,b)$ the corresponding $n$-modification is 
$L(-,b)(f^n):L(-,a)\to L(-,b)$ which is the same as $L(f^n,-):L(-,a)\to L(-,b)$, i.e. functor 
$\bold Y\hskip -0.05cm:\hskip -0.05cmL\to \infty\text{\bf -CAT}^{L^{op}}\hskip -0.15cm:\alpha \mapsto L(-,\alpha )$, \ $\alpha \in L$, 
locally coincides with isomorphisms from Yoneda lemma. By lemma 1.3 this functor preserves and reflects equivalences \, $\sim $\, .  \hfill    $\square $
\enddemo 

\vskip 0.2cm
{\bf Remark.} Under assumption that a category $\infty \text{-}\bold{CAT}$ of {\bf weak} categories, functors and $n$-modifications
exist all the above reasons remain essentially the same, i.e. Yoneda lemma and embedding seem to hold in a weak situation.   \hfill  $\square $

\head {\bf 4. (Co)limits} \endhead

{\bf Definition 4.1.} {\bf $\infty $-graph} is a graded set $G=\coprod\limits _{n\ge 0}G^n$ with two unary operations 
$d,c:\coprod\limits _{n\ge 1}G^n\to \coprod\limits _{n\ge 0}G^n$ of degree $-1$ such that $d^2=dc, \ c^2=cd$. \hfill $\square $  
\vskip 0.1cm
{\bf Definition 4.2.} {\bf $\infty $-diagram} $D:G\to L$ from $\infty $-graph $G$ to $\infty $-category $L$ is a function of degree $0$
which preserves operations $d,c$. \hfill $\square $

\proclaim{\bf {Proposition 4.1}} All diagrams from $G$ to $L$, natural transformations, modifications form $\infty $-category 
\, $\bold{Dgrm}_{G,L}$ \, in the same way as functor category $\infty \text{-}\bold {CAT}(L',L)$. \hfill $\square $
\endproclaim
For a given object $a\in L^0$ \,{\bf constant diagram} to \, $a$ \, is \, $\Delta (a):G\to L: g\mapsto e^na$ \, if \, $g\in G^n$. 
$\Delta :L\to \bold{Dgrm}_{G,L}$ is an $\infty $-functor.
\vskip 0.1cm

Denote $\{e\}\alpha :=\{\alpha , e\alpha , e^2\alpha ,...,e^n\alpha ,...\}$, \, $\alpha \in L$.

\vskip 0.1cm
{\bf Definition 4.3.} Diagram $D:G\to L$ has
\item{$\bullet $} {\bf limit} if functor $\bold{Dgrm}_{G,L}(\Delta (-),D):L^{op}\to \infty \text{\bf -CAT}$ is representable. \newline
If $\xymatrix{\nu :L(-,a) \ar[r]^-{\sim } & \bold{Dgrm}_{G,L}(\Delta (-),D)}$ is the equivalence then \newline
$\nu _a(\{e\}ea)\subset \bold{Dgrm}_{G,L}(\Delta (a),D)$ is called {\bf limit cone} over $D$, $a$ is its {\bf vertice} (or diagram {\bf limit $lim\, D$}), $\nu _a(ea)$ are its {\bf edges}, $\nu _a(e^ka),\, k>1$, are identities up to $\sim $
\item{$\bullet $} {\bf colimit} if functor $\bold{Dgrm}_{G,L}(D,\Delta (-)):L\to \infty \text{\bf -CAT}$ is representable. \newline
If $\xymatrix{\nu :L(a,-) \ar[r]^-{\sim } & \bold{Dgrm}_{G,L}(D,\Delta (-))}$ is the equivalence then \newline
$\nu _a(\{e\}ea)\subset \bold{Dgrm}_{G,L}(D,\Delta (a))$ is called {\bf colimit cocone} over $D$, $a$ is its {\bf vertice} (or diagram {\bf colimit $colim\, D$}), $\nu _a(ea)$ are its {\bf edges}, $\nu _a(e^ka),\, k>1$, are identities up to $\sim $ \hfill $\square $

\vskip 0.2cm
{\bf Remark.} Conditions on equivalence $\nu $ in the above definition can be strengthened. If it is a (natural) isomorphism then 
(co)limits are called {\bf strict} and as a rule they are different from {\bf weak} ones \cite{Bor1}.

\vskip 0.2cm
\proclaim{\bf Proposition 4.2} For strict (co)limits the following is true  \vskip 0.25cm
\item{$\bullet $} $
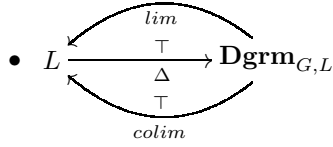
\xymatrix{L \ar[rr]_(0.49){\Delta }^(0.49){\top } & & \bold{Dgrm}_{G,L} \ar@/^5ex/[ll]^-{colim}_-{\top } \ar@/_5ex/[ll]^-{lim}}$ 
\vskip 0.1cm
\item{$\bullet $} Strict right adjoints preserve limits (strict left adjoints preserve colimits).
\endproclaim
\demo{Proof}
\item{$\bullet $} It is immediate from definition 4.3 and proposition 5.1. 
\item{$\bullet $} The argument is the same as for first oder categories (see example 13.a, point 1) [the essential thing is that
a strict adjunction is determined by (triangle) identities which are preserved under $\infty $-functors].    \hfill $\square $
\enddemo

\vskip 0.2cm
\centerline{\bf Examples}
\vskip 0.2cm
\item{$1.$} (strict binary products in $2\text{-}\bold{Top}$ and $2\text{-}\bold{CAT}$) They coincide with '$1$-dimensional' products.
Mediating $2$-cell arrow is given componentwise 
\hskip 1.5cm$\vcenter
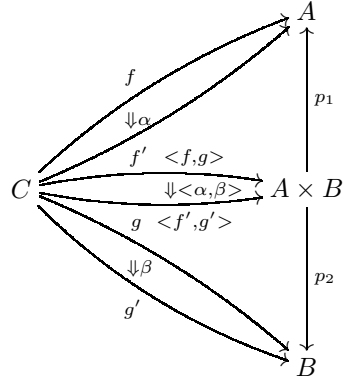
{\xymatrix{ & & & A \\
 & & & \\
C \ar@/^0.6pc/[rrruu]^(0.45){f}_(0.39){\Downarrow \alpha } \ar@/_/[rrruu]_(0.35){f'} \ar@/^/[rrrdd]^(0.35){g}_(0.43){\Downarrow \beta } \ar@/_0.7pc/[rrrdd]_(0.47){g'} \ar@/^/[rrr]^(0.6){<f,g>}_(0.65){\Downarrow <\alpha ,\beta >} \ar@/_/[rrr]_(0.6){<f',g'>} & & & A\times B \ar[uu]_-{p_1} \ar[dd]^-{p_2}\\
 & & & \\
 & & & B}}$
\vskip 0.1cm
\item{$2.$} ('equalizer' of a $2$-cell in $2\text{-}\bold{CAT}$) \cite{Bor1} For a given $2$-cell $\xymatrix{\bold{A} \ar@/^/[r]^-{F}_-{\ \Downarrow \alpha } \ar@/_/[r]_{G} & \bold{B}}$ in $2\text{-}\bold{CAT}$ 
its strict limit is a subcategory $\bold{E}\hookrightarrow \bold{A}$ such that $F(A)=G(A)$ and $\alpha _A=1_{F(A)}:F(A)\to G(A)$ (on objects), and $F(f)=G(f)$ (on arrows).
\vskip 0.1cm
\item{$3.$} (strict and weak pullbacks in $2\text{-}\bold{CAT}$) \cite{Bor1} Let $\Cal P$ be a '$2$-dimensional' graph 
$1@>x>>0@<y<<2$ with trivial $2$-cells, $F:\Cal P\to 2\text{-}\bold{CAT}$ be a $2$-functor. Then its limit is a pullback
diagram in $2\text{-}\bold{CAT}$ \hskip 0.3cm$\vcenter
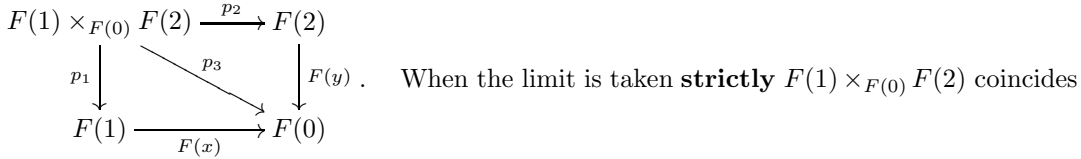
{\xymatrix{F(1)\times _{F(0)}F(2) \ar[r]^-{p_2} \ar[d]_-{p_1} \ar[dr]^-{p_3} & F(2) \ar[d]^-{F(y)} \\
F(1) \ar[r]_-{F(x)} & F(0)}}$. \hskip 0.3cm When the limit is taken {\bf strictly} $F(1)\times _{F(0)}F(2)$ coincides with '$1$-dimensional' pullback, 
i.e. $F(1)\times _{F(0)}F(2)\hookrightarrow F(1)\times F(2)$ is a subcategory consisting of objects $(A,B), \ A\in Ob\, F(1), \ B\in Ob\, F(2), \ F(x)(A)=F(y)(B)$ 
and arrows $(f,g), \ f\in Ar\, F(1), \ g\in Ar\, F(2), \ F(x)(f)=F(y)(g)$. When the limit is taken {\bf weakly} 
$F(1)\times _{F(0)}F(2)$ is not a subcategory of product $F(1)\times F(2)$. It consists of $5$-tuples $(A,B,C,f,g)$, $A\in Ob\, F(1)$, 
$B\in Ob\, F(2)$, $C\in Ob\, F(0)$, \ $f:F(x)(A)@>\sim >>C$, \ $g:F(y)(B)@>\sim >>C$ are isomorphisms, with arrows $(a,b,c)$, 
$a:A\to A'$, $b:B\to B'$, $c:C\to C'$ such that $c\circ f=f'\circ F(x)(a)$, $c\circ g=g'\circ F(y)(b)$. 
Projections $p_1, p_2, p_3$ are obvious. The pullback square commutes up to isomorphisms \ $f:F(x)\circ p_1\Rightarrow p_3$, \ $g:F(y)\circ p_2\Rightarrow p_3$.   \hfill  $\square $

\head {\bf 5. Adjunction}\endhead

{\bf Definition 5.1.} The situation $\xymatrix{L\ar@/^/[r]^-{F} & L' \ar@/^/[l]^-{G}_-{\perp }}$ 
(where $L,L'$ are $\infty $-categories, $F,G$ are $\infty $-functors) is called
\item{$\bullet $} {\bf weak $\infty $-adjunction} if there is an equivalence
$L(-,G(+))\sim L'(F(-),+):L^{op}\times L'\to \infty \text{\bf -CAT}$
(i.e. $L(a,G(b))\sim L'(F(a),b)$ natural in $a\in L^0, \, b\in L^{'0}$),
\item{$\bullet $} {\bf strict $\infty $-adjunction} if there is an isomorphism                                            
$L(-,G(+))\simeq L'(F(-),+):L^{op}\times L'\to \infty \text{\bf -CAT}$ 
(i.e. $L(a,G(b))\simeq L'(F(a),b)$ natural in $a\in L^0, \, b\in L^{'0}$).           \hfill $\square $

\proclaim{\bf {Proposition 5.1}} The following are equivalent
\item{1.} $\xymatrix{L\ar@/^/[r]^-{F} & L' \ar@/^/[l]^-{G}_-{\perp }}$ is a strict $\infty $-adjunction
\item{2.} $\forall b\in L^{'0}$ \, $L'(F(-),b)$ is strictly representable \vskip 0.1cm
\item{3.} $\forall a\in L^0$ \, $L(a,G(-))$ is strictly representable
\endproclaim
\demo{Proof} 
\item{$\bullet $}$1. \Longrightarrow  2., 3.$ \, is immediate
\item{$\bullet $}$2. \Longrightarrow  1.$ \, From criterion of strict representability (see point 1.3) it follows that $\forall \, b\in L^{'0}$ 
there exists a 'universal element' $(\beta ^0_b:F(G(b))\to b)\in L^{'0}(F(G(b)),b)$ such that $\forall \, (f^n:F(c)\to b)\in L^{'n}(F(c),b)$
$\exists \, !$ $n$-arrow $(g^n:c\to G(b))\in L^n(c,G(b))$ with $f^n=\mu _{F(c),F(G(b)),b}(e^n\beta ^0_b,F(g^n))$
\hskip 0.8cm$\vcenter{\xymatrix{G(b) & F(G(b)) \ar[r]^-{e^n\beta ^0_b} & b \\
c \ar@{-->}[u]^-{\exists ! g^n} & F(c) \ar[u]^-{F(g^n)} \ar[ur]_-{\forall f^n} & }}$ \vskip 0.0cm
Consequently, $\forall \, (f^n:b'\to b)\in L^{'n}(b',b)$ the diagram holds 
\hskip 0.0cm$\vcenter{\xymatrix{G(b) & F(G(b)) \ar[r]^-{e^n\beta ^0_b} & b \\
G(b') \ar@{-->}[u]^-{G(f^n)} & F(G(b')) \ar[r]_-{e^n\beta ^0_{b'}} \ar[u]^-{F(G(f^n))} & b' \ar[u]_-{f^n} }\hskip -0.1cm}$
\item{} It shows that assignment $Ob\, L'\ni b\mapsto G(b)\in Ob\, L$ is extendable to a functor $G:L'\to L$ (in an essentially unique way) 
and that $\beta ^0:FG\to 1_{L'}$ is a natural transformation ({\bf counit} \, $\varepsilon $ of the adjunction $F\dashv G$). 
\item{} Isomorphism $\varphi _{c,b}:L'(F(c),b)\to L(c,G(b))$ such that $\vcenter{\xymatrix{F(G(b)) \ar[r]^-{e^n\beta ^0_b} & b \\
F(c) \ar[u]^-{F(\varphi _{c,b}(f^n))} \ar[ur]_-{f^n} & }}$
is natural in $c\in Ob\, L, \, b\in Ob\, L'$ because of the naturality square 
\item{} $\vcenter{\xymatrix{c & b \ar[d]_-{f^n} & & L'(F(c),b) \ar[r]^-{\varphi _{c,b}} \ar[d]_-{L'(F(g^n),f^n)} & L(c,G(b)) \ar[d]^-{L(g^n,G(f^n))} \\
c' \ar[u]^-{g^n} & b' & & L'(F(c'),b') \ar[r]_-{\varphi _{c',b'}} & L(c',G(b'))}}$
\item{} (indeed, $\forall h^n\in L^{'n}(F(c),b)$ \, $G(f^n)*\varphi _{c,b}(h^n)*g^n\sim \varphi _{c',b'}(f^n*h^n*F(g^n))$, where $*$ is the horizontal composite,
since $e^n\beta ^0_{b'}*F(G(f^n)*\varphi _{c,b}(h^n)*g^n)\sim f^n*e^n\beta ^0_b*F(\varphi _{c,b}(h^n))*F(g^n)\sim f^n*h^n*F(g^n)$)
\item{$\bullet $}$3. \Longrightarrow  1.$ \, is similar to \, 2. $\Longrightarrow $ 1. \hfill $\square $
\enddemo

\vskip 0.1cm
{\bf Remark.} Analogous statement for a weak $\infty $-adjunction is not true. In the above proof 'universal elements' were 
essentially used.       \hfill $\square $

\vskip 0.1cm
{\bf Definition 5.2.} For a given {\bf strict} $\infty $-adjunction $\xymatrix{L\ar@/^/[r]^-{F} & L' \ar@/^/[l]^-{G}_-{\perp }}$
\item{$\bullet $} universal elements $\varepsilon _b:F(G(b))\to b$ representing functors $L'(F(-),b)$ ($b\in Ob\, L'$ is a parameter)
form natural transformation $\varepsilon :FG\to 1_{L'}$ which is called {\bf counit} of the adjunction,
\item{$\bullet $} universal elements $\eta _a:a\to G(F(a))$ representing functors $L(a,G(-))$ ($a\in Ob\, L$ is a parameter) 
form natural transformation $\eta :1_{L}\to GF$ which is called {\bf unit} of the adjunction.   \vskip 0.0cm \hfill $\square $

\vskip 0.1cm
{\bf Remark.} For a {\bf weak} $\infty $-adjunction no usefull unit and counit exist.   \hfill  $\square $

\proclaim{\bf Proposition 5.2}
\item{$\bullet $} For both weak and strict adjunctions:
composite of left adjoints is a left adjoint (composite of right adjoints is a right adjoint).  
\item{$\bullet $} For a weak (strict) adjunction a right or left adjoint is determined uniquely up to equivalence $\sim $ 
(up to isomorphism $\simeq $).                    
\endproclaim
\demo{Proof} 
\item{$\bullet $} If $\xymatrix{L \ar@/^/[r]^-{F}_-{\perp} & L' \ar@/^/[r]^-{F'}_-{\perp } \ar@/^/[l]^-{G} & L'' \ar@/^/[l]^-{G'}}$ 
then $L''(F'Fl,l'')\sim L'(Fl,G'l'')\sim L(l,GG'l'')$ (composite of natural equivalences). \ [For a strict adjunction the same reason works]
\item{$\bullet $} Assume, $L'(a,G'b)\hskip -0.05cm\sim \hskip -0.05cmL(Fa,b)\hskip -0.05cm\sim \hskip -0.05cmL'(a,Gb)$ are natural equivalences then $L'(-,G'b)\hskip -0.05cm\sim \hskip -0.05cmL'(-,Gb)$ is a natural
transformation (equivalence) natural in $b$. Then, by Yoneda embedding, $G'b\sim Gb$ naturally in $b$, i.e. $G'\sim G$.
[Again, changing $\sim $ with $\simeq $ still works].              \hfill $\square $                               
\enddemo 

\proclaim{\bf Proposition 5.3} For a strict adjunction $\xymatrix{L \ar@/^/[r]^-{F} & L' \ar@/^/[l]^-{G}_-{\perp } }$ Kan definition and 
definition via 'unit-counit' coincide, i.e. the following are equivalent
\item{$\bullet $} $\varphi _{a,b}:L(a,G(b))\simeq L'(F(a),b):\varphi ^*_{a,b}$ natural in $a\in L^0, \, b\in L^{'0}$,
\item{$\bullet $} $\exists $ natural transformations $\eta :1_L\to GF$ and $\varepsilon :FG\to 1_{L'}$ satisfaying triangle 
identities $\varepsilon F\circ _1F\eta =1_F$ and $G\varepsilon \circ _1\eta G=1_G$. 
\endproclaim
\demo{Proof} For a strict adjunction the same proof as for first order categories works.\vskip 0.1cm
\item{$\bullet $} Universal elements $\eta _a, \varepsilon _b$ for functors $L(a,G(-)),L'(F(-),b)$ mean that they are images of
$1_{F(a)}, 1_{G(b)}$ under functors $\varphi ^*_{a,F(a)},\varphi _{G(b),b}$\, , \, i.e. 
\hskip 0.5cm$\xymatrix{FGFa \ar[r]^-{\varepsilon _{Fa}} & Fa \\
Fa \ar@{-->}[u]^-{F\eta _a} \ar[ur]_-{1_{Fa}} & }$\hskip 0.5cm
$\xymatrix{Gb \ar[r]^-{\eta _{Gb}} \ar[dr]_-{1_{Gb}} & GFGb \ar@{-->}[d]^-{G\varepsilon _b} \\
 & Gb }$\vskip -0.1cm
\item{}\hskip 9.7cm (strict equalities)
\vskip 0.15cm\item{$\bullet $} Define maps $\cases \varphi _{a,b}(f^n):=e^n(\varepsilon _b)\circ _{n+1}F(f^n),  &  f^n\in L^n(a,G(b)) \\
\varphi ^*_{a,b}(g^n):=G(g^n)\circ _{n+1}e^n(\eta _a),   &  g^n\in L^{'n}(F(a),b)   \endcases $
\vskip 0.25cm\item{}They are functors $\cases \varphi _{a,b}:=\varepsilon _b*F(-):L(a,G(b))\to L'(F(a),b) &  \\
\varphi ^*_{a,b}:=G(-)*\eta _a:L'(F(a),b)\to L(a,G(b)) &   \endcases $
\hskip -0.3cmand inverses to each other:\vskip 0.15cm 
\item{}$\varphi ^*_{a,b}(\varphi _{a,b}(f^n))=\varphi ^*_{a,b}(e^n\varepsilon _b\circ _{n+1}F(f^n))=
G(e^n\varepsilon _b\circ _{n+1}F(f^n))\circ _{n+1}e^n\eta _a=e^nG(\varepsilon _b)\circ _{n+1}(GF(f^n)\circ _{n+1}e^n\eta _a)=
e^nG(\varepsilon _b)\circ _{n+1}(e^n\eta _{G(b)}\circ _{n+1}f^n)=e^{n+1}G(b)\circ _{n+1}f^n=f^n$, and similar, 
$\varphi _{a,b}(\varphi ^*_{a,b}(g^n))=g^n$.        
\item{}Naturality (e.g., of $\varphi _{a,b}$) follows from the square 
$\vcenter{\xymatrix{L(a,G(b))  \ar[r]^-{\varphi _{a,b}} \ar[d]_-{L(x^m,G(y^m))} &     L'(F(a),b)  \ar[d]^-{L'(F(x^m),y^m)}  \\
L(a',G(b'))   \ar[r]_-{\varphi _{a',b'}}       &         L'(F(a'),b')    }}$       
\vskip 0.1cm\item{}$\bigl (\varphi _{a',b'}(L(x^m,G(y^m))(f^n))=\varphi _{a',b'}(G(y^m)*f^n*x^m)=
\varepsilon _{b'}*FG(y^m)*F(f^n)*F(x^m)=y^m*\varepsilon _{b}*F(f^n)*F(x^m)=L'(F(x^m),y^m)(\varphi _{a,b}(f^n))$, where $n=0$ or 
$m=0\bigr )$.                    \hfill     $\square $
\enddemo

\vskip 0.35cm
\centerline{\bf Examples of higher order adjunction}

\vskip 0.2cm
\item{1.} Every usual 1-adjunction $\xymatrix{A \ar@/^/[r] & B \ar@/^/[l]_-{\bot }}$ is $\infty $-1-adjunction for trivial $\infty $-extentions of $A$ and $B$.
\item{2.} Gelfand-Naimark dual 1-adjunction $\xymatrix{\bold{C^*Alg}^{op} \ar@/^/[r] & \bold{CHTop} \ar@/^/[l]_-{\bot }}$ is extendable to $\infty $-2-adjunction (see point 9).
\proclaim{\hskip -0.5cm{\rm 3.}\ {\sl Quillen theorem} \cite{Mac}} Let $\bold{\Delta }$ be a category of finite linearly ordered sets, $\bold{Set}^{\bold{\Delta }^{op}}$ category of simplicial sets, $Ho(\bold{Top}):=(2\text{\bf -Top})^{(1)}$, $Ho(\bold{Set}^{\Delta ^{op}}):=(2\text{\bf -Set}^{\Delta ^{op}})^{(1)}$. Then
$$\xymatrix{\Delta \ar@{{>}->}[dr] \ar@{{>}->}[d]|-{\text{\rm Yoneda}} &  \\
\bold{Set}^{\Delta ^{op}} \ar@{-->}[r] \ar[d] & \bold{Top} \ar@/^2.5ex/[l]_-{\bot} \ar[d]\\
Ho(\bold{Set}^{\Delta ^{op}}) \ar[r] & Ho(\bold{Top}) \ar@/^2.5ex/[l]_-{\bot}
}$$
\endproclaim \hfill $\square $
\vskip 0.0cm 
So, the top adjunction is actually 2-adjunction (or $\infty $-2-adjunction).
\vskip 0.0cm 
All the above adjunctions are strict.

\head {\bf 6. Concrete duality for $\infty$-categories}\endhead

\vskip 0.0cm
Duality preseves all categorical properties. It is the one extremity of functors (in the sense of invariants) especially useful 
when the duality is concrete. It is significant that concrete duality for $\infty $-categories behaves the same as for $1$-categories.

\vskip 0.25cm
{\bf Definition 6.1.} 
\item{$\bullet $} {\bf Duality} is just an equivalence $L^{op}\sim L'$.
\item{$\bullet $} {\bf Concrete duality} over $\Bbb B\hookrightarrow \infty \text{-}\bold{CAT}$ is a duality 
$\xymatrix{L^{op} \ar@/^0.6pc/[r]^-{G}_-{\sim } & L' \ar@/^0.6pc/[l]^-{F}_-{\sim } }$
such that $\exists $ (faithful) forgetful functors $U:L\to \Bbb B$, \, $V:L'\to \Bbb B$ and objects $\tilde A\in L^0$, \, 
$\tilde B\in L^{'0}$ such that 
\vskip 0.15cm\itemitem{$\bullet $} $U(\tilde A)\sim V(\tilde B)$, 
\itemitem{$\bullet $} \hbox{$V\circ _1G\sim L(-,\tilde A)$, \ $U\circ _1F^{op}\sim L'(-,\tilde B)$} \hskip 1cm
$\vcenter{\hbox{\xymatrix{L^{op} \ar[r]^-{G} \ar[rd]_-{L(-,\tilde A)} & L' \ar[d]^-{V}\\
 & {\Bbb B}}}}$ \hskip 1cm
$\vcenter{\hbox{\xymatrix{L^{'op} \ar[r]^-{F^{op}} \ar[rd]_-{L'(-,\tilde B)} & L \ar[d]^-{U}\\
 & {\Bbb B}}}}$ 
\vskip 0.15cm
\item{} Representing objects $\tilde A\in L^0,\ \tilde B\in L^{'0}$ are called together a {\bf dualizing} or {\bf schizophrenic object} for the given concrete duality\cite{P-Th}.
\vskip 0.15cm\item{} \hskip 0cm[for a {\bf concrete dual adjunction} the definition is similar]           \hfill    $\square $ 

\vskip 0.25cm
\proclaim{\bf {Proposition 6.1} (representable forgetfuls $\Rightarrow $ concrete dual adjunction)} Let $(L,U)$, $(L',V)$ be (weakly) dually adjoint 
$\infty $-categories $\xymatrix{L^{op} \ar@/^/[r]^-{G}_-{\top } & L' \ar@/^/[l]^-{F} }$ 
with representable forgetful functors \ 
$U\sim L(A_0,-):L\to \Bbb B$, \, $V\sim L'(B_0,-):L'\to \Bbb B$ (where $\Bbb B\hookrightarrow \infty \text{-}\bold{CAT}$ is a subcategory). 
\, Then this adjunction is concrete over
$\Bbb B$ with dualizing object $(\tilde A, \tilde B)$, where $\tilde A:=F(B_0)$, \, $\tilde B:=G(A_0)$, i.e. 
\vskip 0.3cm
\item{$\bullet $} $U(\tilde A)\sim V(\tilde B)$
\item{$\bullet $} \hbox{$V\circ _1G\sim L(-,\tilde A)$, \ $U\circ _1F^{op}\sim L'(-,\tilde B)$} \hskip 1cm
$\vcenter{\hbox{\xymatrix{L^{op} \ar[r]^-{G} \ar[rd]_-{L(-,\tilde A)} & L' \ar[d]^-{V}\\
 & {\Bbb B}}}}$ \hskip 1cm
$\vcenter{\hbox{\xymatrix{L^{'op} \ar[r]^-{F^{op}} \ar[rd]_-{L'(-,\tilde B)} & L \ar[d]^-{U}\\
 & {\Bbb B}}}}$
\endproclaim

\demo{Proof}
\vskip 0.1cm
\item{$\bullet $} $U(\tilde A)=UF(B_0)\sim L(A_0,FB_0)\sim L'(B_0,GA_0)\sim VGA_0=V\tilde B$
\item{$\bullet $} $VG(-)\sim L'(B_0,G(-))\sim L(-,FB_0)=L(-,\tilde A)$ \, (and similar, $UF(-)\sim L'(-,\tilde B)$)        \hfill $\square $
\enddemo

\vskip 0.2cm
{\bf Remarks.}
\item{$\bullet $} Concrete duality as above should be called {\bf weak}. {\bf Strict} variants of definition 6.1 and 
proposition 6.1 also exist (by changing $\sim $ with isomorphism $\simeq $ and weak dual adjunction with the strict one).
\item{$\bullet $} (Weak or strict) concrete duality (dual adjunction) is given essentially by hom-functors 
which admit lifting along forgetful functors (to obtain proper values). Representing objects of these functors have equivalent 
(or isomorphic) underlying objects.      
\item{$\bullet $} For usual 1-dimensional categories $\Bbb B=\bold{Set}\hookrightarrow \infty \text{-}\bold{CAT}$ 
($\infty $-1-subcategory). For dimension $n$, as a rule, $\Bbb B=n\text{-}\bold{Cat}\hookrightarrow \infty \text{-}\bold{CAT}$ 
($\infty \text{-}n$-subcategory of small $(n-1)$-categories).    \hfill    $\square $

\subhead 6.1. Natural and non natural duality \endsubhead

\vskip 0.2cm
{\bf Definition 6.1.1.} \item{$\bullet $} For hom-set $L(A,\tilde A)$ and element $(x:A_0\to A)\in L^0(A_0,A)$ {\bf evaluation functor} 
at point $x$ is \, $ev_{A,x}:=L(x,\tilde A):L(A,\tilde A)\to L(A_0, \tilde A)$ \, ($ev_{A,x}\in \Bbb B^1\hookrightarrow \infty \text{-}\bold{CAT}^1$).
\item{}Similarly, {\bf evaluation $(n-1)$-modification} $ev_{A,x^n}$, $n=1,2,\dots $, for $x^n\in L^n(A_0,A)$ is $L(x^n,\tilde A)\in \Bbb B^n(L(A,\tilde A),L(A_0,\tilde A))$.
\item{$\bullet $} For a forgetfull functor $V:L'\to \Bbb B$ an arrow $f^n:V(Y)\to V(Y')\in \Bbb B^n(V(Y),V(Y'))$ is called 
$L'${\bf -arrow} if $\exists \Phi ^n:Y\to Y'\in L^{'n}(Y,Y')$ such that $V(\Phi ^n)=f^n$.
\item{$\bullet $} Lifting of hom-functor $V\circ G\sim L(-,\tilde A)$ is called 
{\bf initial} \cite{A-H-S} if $\forall A\in L^0\ \forall Y\in L^{'0}\ \forall f^n:V(Y)\to L(A,\tilde A)\in \Bbb B^{n}(V(Y),L(A,\tilde A))$ $f^n$ 
is an $L'$-arrow iff $\forall (x^n:A_0\to A)\in L^n(A_0,A)$ 
$ev_{A,x^n}\circ _{n+1}f^n:V(Y)\to L(A_0,\tilde A)\in \Bbb B^n(V(Y),L(A_0,\tilde A))$ is an $L'$-arrow.
\item{$\bullet $} If liftings of hom-functors $V\circ G\sim L(-,\tilde A),\ U\circ F\sim L'(-,\tilde B)$ are both initial then the concrete dual adjunction $\xymatrix{L^{op} \ar@/^/[r]^-{G}_-{\top } & L' \ar@/^/[l]^-{F} }$, if exists, is called {\bf natural} \cite{Hof, P-Th}, and otherwise, non natural. \hfill $\square $

\vskip 0.2cm
Even if $U\tilde A\sim V\tilde B$ and $\forall A\in L^0, B\in L^{'0}$ $\Bbb B$-objects $L(A,\tilde A), \, L'(B,\tilde B)$ can be lifted 
to $L', L$ hom-functors $L(-,\tilde A),\, L'(-,\tilde B)$ need not to be lifted (which happens only if lifting of the
assignments $A\mapsto L(A,\tilde A),\, B\mapsto L'(B,\tilde B)$ can be extended functorially over all cells). 
%Such possibility exists if the following condition holds.

\vskip 0.15cm
{\bf Initial lifting condition for evaluation cones} \newline 
$\{ev_{A,x^n}\in \Bbb B^n(L(A,\tilde A),L(A_0,\tilde A))\}_{x^n\in L^n(A_0,A)}^{\, n\, \in \, \Bbb N}, \, 
\{ev_{B,y^n}\in \Bbb B^n(L'(B,\tilde B),L'(B_0,\tilde B))\}_{y^n\in L^{'n}(B_0,B)}^{\, n\, \in \, \Bbb N}$ consists of the following requirements:
\item{$\bullet $} hom-categories of the form $L(A,\tilde A),\, L'(B,\tilde B)\in Ob\, (\Bbb B)$ can be lifted to $L',L$
\item{$\bullet $} evaluation cones \newline 
$\{ev_{A,x^n}\hskip -0.05cm\in \hskip -0.05cm\Bbb B^n(L(A,\tilde A),L(A_0,\tilde A))\}_{x^n\in L^n(A_0,A)}^{\, n\, \in \, \Bbb N}, \, 
\{ev_{B,y^n}\hskip -0.05cm\in \hskip -0.05cm\Bbb B^n(L'(B,\tilde B),L'(B_0,\tilde B))\}_{y^n\in L^{'n}(B_0,B)}^{\, n\, \in \, \Bbb N}$
can be lifted to $\{ev_{A,x^n}\in L^{'n}(G(A),\tilde B)\}_{x^n\in L^n(A_0,A)}^{\, n\, \in \, \Bbb N}\, , \ \{ev_{B,y^n}\in L^n(F(B),\tilde A)\}_{y^n\in L^{'n}(B_0,B)}^{\, n\, \in \, \Bbb N}$ \, in $L',L$   
\item{$\bullet $} $\forall f^n\in \Bbb B^n(VX,L(A,\tilde A))$ $f^n$ is $L'$-arrow iff $\forall x^n\hskip -0.05cm\in \hskip -0.05cmL^n(A_0,A)$ \, $\mu (ev_{A,x^n},f^n)\hskip -0.05cm\in \hskip -0.05cm\Bbb B^n(VX,L(A_0,\tilde A))$ is $L'$-arrow
(and, symmetrically, $\forall g^n\in \Bbb B^n(UY,L'(B,\tilde B))$ $g^n$ is $L$-arrow iff $\forall y^n\hskip -0.05cm\in \hskip -0.05cmL^{'n}(B_0,B)$ \, $\mu (ev_{B,y^n},g^n)\hskip -0.05cm\in \hskip -0.05cm\Bbb B^n(UY,L'(B_0,\tilde B))$ is $L$-arrow) \hfill $\square $

\vskip 0.1cm
Denote further (in the the following proof) lifted evaluation maps by $ev_{A,x}$ 
(or like that) and underlying evaluation maps in $\Bbb B$ by $|ev_{A,x}|$.

\vskip 0.2cm
\proclaim{\bf Proposition 6.1.1} If two strict $\infty $-categories $L,L'$ concrete over 
$\Bbb B\hookrightarrow \infty \text{-}\bold{CAT}$ with representable (strictly faithful) forgetful functors
$U=L(A_0,-), \ V=L'(B_0,-)$ have objects $\tilde A\in L^0, \, \tilde B\in L^{'0}$ such that
\item{$\bullet $} $U\tilde A\sim V\tilde B$
\item{$\bullet $} hom-functors $L(-,\tilde A):L^{op}\to \Bbb B, \  L'(-,\tilde B):L^{'op}\to \Bbb B$ satisfy 
{\bf initial lifting condition for evaluation cones} \newline 
%\item{$\bullet $} $\Bbb B$ has property $\Bbb B(X,\Bbb B(Y,Z))\sim \Bbb B(Y,\Bbb B(X,Z))$ (natural in all arguments) 
then $\exists $ natural {\bf strict} concrete dual adjunction \hskip 0.3cm 
$\xymatrix{L^{op} \ar@/^/[r]^-{G}_-{\top } & L' \ar@/^/[l]^-{F} }$ \hskip 0.3cm 
$L(A,FB)\underset {\text{nat. iso}}\to {\simeq }L'(B,GA)$\hskip 0.5cm
$\vcenter{\hbox{\xymatrix{L^{op} \ar[r]^-{G} \ar[rd]_-{L(-,\tilde A)} & L' \ar[d]^-{V}\\
 & {\Bbb B}}}}$ \hskip 0.5cm
$\vcenter{\hbox{\xymatrix{L^{'op} \ar[r]^-{F^{op}} \ar[rd]_-{L'(-,\tilde B)} & L \ar[d]^-{U}\\
 & {\Bbb B}}}}$ \hskip 0.5cm
with $(\tilde A,\tilde B)$, its schizophrenic object.
\endproclaim
\demo{Proof} 
\item{$\bullet $} $L(A,\tilde A), \, L'(B,\tilde B)$ are lifted to $L',L$ by condition.
\item{$\bullet $} Let $f^n\in L^n(A,A')$, then $L(f^n,\tilde A):L(A',\tilde A)\to L(A,\tilde A)$ is $L'$-arrow since
$ev_{A,a^n}\circ _{n+1}L(f^n,\tilde A):=L(a^n,\tilde A)\circ _{n+1}L(f^n,\tilde A)=L(f^n\circ _{n+1}a^n,\tilde A)=:ev_{A',f^n\circ _{n+1}a^n}$, 
which is liftable $\forall a^n\in L^n(A_0,\tilde A)$. Therefore, $L(f^n,\tilde A)$ is $L'$-arrow, and similarly, $L'(g^n,\tilde B)$ is $L$-arrow, i.e., $\exists $ 
maps $\xymatrix{L^{op} \ar@/^/[r]^-{G} & L' \ar@/^/[l]^-{F}}$, which are obviously functorial. Why do they give an adjunction?
\item{$\bullet $} (unit and counit) $1$-arrow (unit) $\eta _B:B\to GFB$ is given by 
$|\eta _B|=:V\eta _B:|B|\to |GFB|:b\mapsto [ev_{B,b}:FB\to \tilde A]$, 
$b\in |B|=L'(B_0,B)$, $|GFB|=L(FB,\tilde A)$, $|ev_{B,b}|:|FB|\to |\tilde A|$, $|FB|=L'(B,\tilde B)$, 
$|\tilde A|=L(A_0,\tilde A)\sim L'(B_0,\tilde B)$.
Why can $|\eta _B|$ be lifted to $L'$? Take composite with evaluation maps 
$|ev_{FB,c}|\circ _1|\eta _B|(b)=|ev_{FB,c}|(ev_{B,b})=|ev_{B,b}|(c)=|c|(b)$, 
where $c\in |FB|^0=L^{'0}(B,\tilde B)=L^0(A_0,FB)$, $b\in |B|^n$. So, $|ev_{FB,c}|\circ _1|\eta _B|=|c|$ is $L'$-arrow.
Therefore, $|\eta _B|$ is $L'$-arrow. Counit is given symmetrically $\varepsilon _A\to FGA$, 
$|\varepsilon _A|:|A|\to |FGA|:a\mapsto [ev_{A,a}:GA\to \tilde B]$, $|A|=L(A_0,A)$, $|FGA|=L'(GA,\tilde B)$, 
$|ev_{A,a}|:|GA|\to |\tilde B|$, $|GA|=L(A,\tilde A)$, $|\tilde B|=L'(B_0,\tilde B)\sim L(A_0,\tilde A)$.
By the same argument $|\varepsilon _A|$ is an $L$-arrow.
\item{$\bullet $} (triangle identities) $G\varepsilon _A\circ _1\eta _{GA}=1_{GA}$, $F\eta _B\circ _1\varepsilon _{FB}=1_{FB}$. 
It is sufficient to prove them for underlying maps. Since forgetful functors are faithful they will follow. \newline 
$|G\varepsilon _A|\circ _1|\eta _{GA}|\overset ?\to =|1_{GA}|$, where $|\eta _{GA}|:|GA|\to |GFGA|$, $|GA|=L(A,\tilde A)$, 
$|GFGA|=L(FGA,\tilde A)$, $\varepsilon _A:A\to FGA$, $|G\varepsilon _A|:|GFGA|\to |GA|$.\newline 
Take $(f^n:A\to \tilde A)\in |GA|=L^n(A,\tilde A)$, $a^m\in |A|=L^m(A_0,A)$. Two cases are possible 
$\cases \hskip -0.1cm(a)\ (f^n,n>0)\ \& \ (a^0):& \hskip -0.25cm||G\varepsilon _A|\circ _1|\eta _{GA}|(f^n)|(a^0)=|L(\varepsilon _A,\tilde A)(ev_{GA,f^n})|(a^0)=|ev_{GA,f^n}\circ _{n+1}\\
\hskip -0.1cm(b)\ (f^0)\ \& \ (a^n,n\ge 0):& \hskip -0.25cm||G\varepsilon _A|\circ _1|\eta _{GA}|(f^0)|(a^n)=|L(\varepsilon _A,\tilde A)(ev_{GA,f^0})|(a^n)\hskip 0.1cm=\hskip 0.1cm|ev_{GA,f^0}\circ _{1}\endcases $ \newline
$\cases \hskip -0.1cm(a) & \hskip -0.25cme^n\varepsilon _A|(a^0)=|ev_{GA,f^n}|\circ _{n+1}e^n|\varepsilon _A|(a^0)=|ev_{GA,f^n}|(ev_{A,e^na^0})=|ev_{A,e^na^0}|(f^n)=|f^n|(a^0)\hskip -0.1cm\\
\hskip -0.1cm(b) & \hskip -0.0cm\varepsilon _A|(a^n)\hskip 0.05cm=\hskip 0.05cm|ev_{GA,f^0}|\circ _{1}|\varepsilon _A|(a^n)\hskip 0.15cm=\hskip 0.15cm|ev_{GA,f^0}|(ev_{A,a^n})\hskip 0.15cm=\hskip 0.15cm|ev_{A,a^n}|(f^0)\hskip 0.15cm=\hskip 0.15cm|f^0|(a^n)\hskip -0.1cm\endcases $ \newline 
$\cases \hskip -0.1cm(a) & \hskip -0.25cm=:\mu ^L_{A_0,A,\tilde A}(f^n,e^na^0)=||1_{GA}|(f^n)|(a^0) \\
\hskip -0.1cm(b) & \hskip -0.25cm=:\mu ^L_{A_0,A,\tilde A}(e^nf^0,a^n)=||1_{GA}|(f^0)|(a^n) \endcases $

\vskip 0.1cm
Second triangle identity holds similarly.
\vskip 0.05cm
\item{$\bullet $} (naturality of $\eta _B, \varepsilon _A$) \ Again, it is sufficient to prove naturality for underlying maps \newline 
$\vcenter{\xymatrix{|B| \ar[r]^-{|\eta _B|} \ar[d]_-{|f|} & |GFB| \ar[d]|-{|GFf|=L(Ff,\tilde A)}\\
|B'| \ar[r]_-{|\eta _{B'}|} & |GFB'|}}$
\hskip 1cmTwo cases are $\cases (a)\ (b^n\in |B|^n, n\ge 0)\ \& \ (f^0\in L^{'0}(B,B')) & \\
(b)\ (b^0\in |B|^0)\ \& \ (f^n\in L^{'n}(B,B')) & \endcases $ \newline 
$(a)$ \hskip 1.5cm $\vcenter{\xymatrix{|B| \ar[r]^-{|\eta _B|} \ar[d]_-{|f^0|} & |GFB| \ar[d]|-{|GFf^0|=L(Ff^0,\tilde A)}\\
|B'| \ar[r]_-{|\eta _{B'}|} & |GFB'|}}$
\hskip 1.5cm $\vcenter{\xymatrix{b^n \ar@{|->}[r] \ar@{|->}[dd] & ev_{B,b^n} \ar@{|->}[d]\\
  & ev_{B,b^n}\circ _{n+1}e^n(Ff^0) \ar@{==>}[d]^-{=}\\
|f^0|(b^n) \ar@{|->}[r] & ev_{B',|f^0|(b^n)}}}$ \newline 
$(b)$ \hskip 1.5cm $\vcenter{\xymatrix{|B| \ar[r]^-{e^n|\eta _B|} \ar[d]_-{|f^n|} & |GFB| \ar[d]|-{|GFf^n|=L(Ff^n,\tilde A)}\\
|B'| \ar[r]_-{e^n|\eta _{B'}|} & |GFB'|}}$
\hskip 1.5cm $\vcenter{\xymatrix{b^0 \ar@{|->}[r] \ar@{|->}[dd] & ev_{B,e^nb^0} \ar@{|->}[d]\\
  & ev_{B,e^nb^0}\circ _{n+1}(Ff^n) \ar@{==>}[d]^-{=}\\
|f^n|(b^0) \ar@{|->}[r] & ev_{B',|f^n|(b^0)}}}$ \newline  
\vskip 0.1cm
\centerline{(recall \ $|f^n|(b^0)\equiv \mu (f^n,e^nb^0)$, \ $|f^0|(b^n)\equiv \mu (e^nf^0,b^n)$)} 
\vskip 0.1cm
Why \ \ $\cases \hskip -0.05cm(a)\ \ ev_{B,b^n}\circ _{n+1}e^n(Ff^0)=ev_{B',|f^0|(b^n)} & \\
\hskip -0.05cm(b)\ \ ev_{B,e^nb^0}\circ _{n+1}(Ff^n)=ev_{B',|f^n|(b^0)} & \endcases $? \newline 
\vskip 0.1cm
Take underlying maps:
\vskip 0.1cm
$\cases \hskip -0.1cm(a)& \hskip -0.2cm|ev_{B,b^n}|\circ _{n+1}e^n|Ff^0|(h^n)=|ev_{B,b^n}|(h^n\circ _{n+1}e^nf^0)=|h^n\circ _{n+1}e^nf^0|(b^n)=\\
\hskip -0.1cm(b)& \hskip -0.2cm|ev_{B,e^nb^0}|\circ _{n+1}|Ff^n|(h^0)=|ev_{B,e^nb^0}|(e^nh^0\circ _{n+1}f^n)=|e^nh^0\circ _{n+1}f^n|(e^nb^0)= \endcases $ 
\vskip 0.1cm
$\cases \hskip -0.1cm(a)& \hskip -0.2cm=|h^n|\circ _{n+1}|e^nf^0|(b^n)=|ev_{B',|f^0|(b^n)}|(h^n), \ \ h^n\in L^{'n}(B',\tilde B)\\
\hskip -0.1cm(b)& \hskip -0.2cm=e^n|h^0|\circ _{n+1}|f^n|(e^nb^0)=|ev_{B',|f^n|(b^0)}|(h^0), \ \ h^0\in L^{'0}(B',\tilde B)\endcases $ 
\vskip 0.1cm
(types of the above arrows are $Ff:FB'\to FB$, $ev_{B,b}:FB\to \tilde A\ (\text{$L$-map})$, $ev_{B',|f|(b)}FB'\to \tilde A\ (\text{$L$-map})$, 
$|ev_{B,b}|:L'(B,\tilde B)\to |\tilde B|=L'(B_0,\tilde B)$, $|ev_{B',|f|(b)}|:L'(B',\tilde B)\to |\tilde B|=L'(B_0,\tilde B)$, 
$|Ff|:L'(B',\tilde B)\to L'(B,\tilde B)$, $|Ff|=L'(f,\tilde B)$). 
\vskip 0.1cm
Therefore, $\eta _B$ is natural. Similarly, $\varepsilon _A$ is natural.
\item{$\bullet $} (iso-functors $\xymatrix{{L(A,FB)\ \ } \ar@/^0.9pc/[r]^-{\theta _{A,B}} & {\ \ L'(B,GA)} \ar@/^0.9pc/[l]^-{\theta ^*_{A,B}}}$)\newline 
Define \ $\cases \theta _{A,B}(f^n):=G(f^n)\circ _{n+1}e^n(\eta _B), & f^n\in L^n(A,FB)\\
\theta ^*_{A,B}(g^n):=F(g^n)\circ _{n+1}e^n(\varepsilon _A), & g^n\in L^{'n}(B,GA) \endcases $\newline 
Let $g^n\in L^{'n}(B,GA)$. Then $\theta _{A,B}(\theta ^*_{A,B}(g^n))\hskip -0.05cm:=G(Fg^n\circ _{n+1}e^n(\varepsilon _A))\circ _{n+1}e^n(\eta _B)=e^n(G\varepsilon _A)\circ _{n+1}GFg^n\circ _{n+1}e^n(\eta _B)\underset {\text{nat. of $\eta _B$}} \to {=} e^n(G\varepsilon _A)\circ _{n+1}e^n(\eta _{GA})\circ _{n+1}g^n\underset {\text{triangle id.}} \to {=}e^n(1_{GA})\circ _{n+1}g^n=e^{n+1}(GA)\circ _{n+1}g^n=g^n$.
Similarly, \ $\theta ^*_{A,B}(\theta _{A,B}(f^n))=f^n, \ f^n\in L^n(A,FB)$. 
$\theta _{A,B},\, \theta ^*_{A,B}$ are obviously functors. So, they are isomorphisms.
\item{$\bullet $} (naturality of $\theta _{A,B},\, \theta ^*_{A,B}$) \ We need to prove diagram \newline 
$\xymatrix{ A & B & L(A,FB) \ar[r]^-{e^n\theta _{A,B}} \ar[d]_-{L(x^n,Fy^n)} & L'(B,GA) \ar[d]^-{L'(y^n,Gx^n)}\\
A' \ar[u]^-{x^n} & B' \ar[u]^-{y^n} & L(A',FB') \ar[r]^-{e^n\theta _{A',B'}} & L'(B',GA') }$ 
\vskip 0.15cm
$L'(y^n,Gx^n)\circ _{n+1}e^n\theta _{A,B}\overset {?}\to {=}e^n\theta _{A',B'}\circ _{n+1}L(x^n,Fy^n)$
\vskip 0.25cm
Two cases are: $\cases \hskip -0.1cm(a) & \hskip -0.2cm(f^0\in L(A,FB))\ \& \ (x^n,y^n,n>0)\\ 
\hskip -0.1cm(b) & \hskip -0.2cm(f^n\in L(A,FB), n\ge 0)\ \& \ (x^0,y^0)\endcases $
\vskip 0.1cm
$(a)$\hskip 0.7cm $\vcenter{\xymatrix{f^0 \ar@{|->}[r] \ar@{|->}[dd] & e^nG(f^0)\circ _{n+1}e^n(\eta _B) \ar@{|->}[d]\\
 & Gx^n\circ _{n+1}(e^nG(f^0)\circ _{n+1}e^n(\eta _B))\circ _{n+1}y^n \ar@{==>}[d]^-{\ ?}_-{=\ } \ar@/^1.5pc/@{==>}@<3cm>[dd]_-{=}^-{({\text{$\eta _B$ is nat.}})} \\
Fy^n\circ _{n+1}e^nf^0\circ _{n+1}x^n \ar@{|->}[r] & G(Fy^n\circ _{n+1}e^nf^0\circ _{n+1}x^n)\circ _{n+1}e^n(\eta _{B'}) \ar@{==>}[d]_-{=\ }\\
 & Gx^n\circ _{n+1}e^nGf^0\circ _{n+1}GFy^n\circ _{n+1}e^n(\eta _{B'})}}$
\vskip 0.1cm
$(b)$\hskip 0.25cm $\vcenter{\xymatrix{f^n \ar@{|->}[r] \ar@{|->}[dd] & G(f^n)\circ _{n+1}e^n(\eta _B) \ar@{|->}[d]\\
 & e^nGx^0\circ _{n+1}(G(f^n)\circ _{n+1}e^n(\eta _B))\circ _{n+1}e^ny^0 \ar@{==>}[d]^-{\ ?}_-{=\ } \ar@/^1.5pc/@{==>}@<3cm>[dd]_-{=}^-{({\text{$\eta _B$ is nat.}})} \\
e^nFy^0\circ _{n+1}f^n\circ _{n+1}e^nx^0 \ar@{|->}[r] & G(e^nFy^0\circ _{n+1}f^n\circ _{n+1}e^nx^0)\circ _{n+1}e^n(\eta _{B'}) \ar@{==>}[d]_-{=\ }\\
 & e^nGx^0\circ _{n+1}Gf^n\circ _{n+1}e^nGFy^0\circ _{n+1}e^n(\eta _{B'})}}$
\vskip 0.15cm
Therefore, $L$ and $L'$ are concretely dually adjoint. This correspondence is natural (by condition) and strict ($\theta _{A,B}$ 
and $\theta ^*_{A,B}$ are isomorphisms).
\hfill $\square $
\enddemo

\vskip 0.2cm
{\bf Corollary.} Concrete natural duality is a {\bf strict} adjunction.                \hfill $\square $

\vskip 0.3cm
\centerline{{\bf Well-known dualities} \cite{P-Th, Bel, A-H-S}}
\vskip 0.25cm

All dualities below are of first order, natural \cite{P-Th}, and obtained by restriction of appropriate dual adjunctions.

\item{1.} $\bold{Vec}_k$ is dually equivalent to itself $\xymatrix{\bold{Vec}_k^{op} \ar@/^/[rr]^-{\bold{Vec}_k(-,k)} && \bold{Vec}_k
\ar@/^/[ll]^-{\bold{Vec}_k(-,k)}_-{\perp }}$, where $\bold{Vec}_k$ is a category of vector spaces over field $k$
\item{2.} $\bold{Set}^{op}\sim \text{\bf Complete Atomic Boolean Algebras}$
\item{3.} $\bold{Bool}^{op}\sim \text{\bf Boolean Spaces}$ (Stone duality), where $\bold{Bool}$ is a category of Boolean rings (every element is idempotent). It is obtained from the dual adjunction $\xymatrix{\bold{CRing} \ar@/^/[rr]^-{\bold{CRing}(-,\bold{2})} && \bold{Top} \ar@/^/[ll]^-{\bold{Top}(-,\bold{2})}_-{\perp }}$, where $\bold 2$ is two-element ring and discrete topological space. $\xymatrix{\bold{CRing}(A,\bold 2) \ar@{^{(}->}[r] & \bold 2^A}$ (subspace in Tychonoff topology)
\item{4.} $\text{\rm hom}(-,\Bbb R/\Bbb Z):\bold{CompAb}^{op}\sim \bold{Ab}$ (Pontryagin duality), where $\bold{CompAb}$, $\bold{Ab}$ are categories of compact abelian groups and abelian groups respectively
\item{5.} $\text{hom}(-,\Bbb C):\text{\bf C*Alg}^{op}\sim \bold{CHTop}$ (Gelfand-Naimark duality), where $\text{\bf C*Alg}$, $\bold{CHTop}$ are categories of commutative $\Bbb C^*$-algebras and compact Hausdorff spaces. $\xymatrix{\text{\bf C*Alg}(A,\Bbb C) \ar@{^{(}->}[r] & {\Bbb C}^A}$ (subspace in Tychonoff topology)

\head {\bf 7. Vinogradov duality} \endhead

Let $K$ be a commutative ring, $A$ a commutative algebra over $K$, $A\text{-}\bold{Mod}\hookrightarrow K\text{-}\bold{Mod}$ be categories
of modules over $A$ and $K$ respectively. 

\vskip 0.1cm
{\bf Definition 7.1.} \cite{V-K-L} For $P,Q \in Ob\, (A\text{-}\bold{Mod})$
\item{$\bullet $} $K$-linear maps \newline
$l(a):=a\, \cdot \, -, r(a):=-\, \cdot \, a, \delta (a):=l(a)-r(a):K\text{-}\bold{Mod}(P,Q)\to K\text{-}\bold{Mod}(P,Q)$ 
are called {\bf left, right multiplications} and {\bf difference operator} (by element $a\in A$),
\item{$\bullet $} $K$-linear map $\Delta :P\to Q$ is a {\bf differential operator of order $\le r$} \ if \ $\forall a_0, a_1,\dots , a_r \in A$ 
$\delta _{a_0,a_1,\dots ,a_r}(\Delta )=0$, where $\delta _{a_0,a_1,\dots ,a_r}:=\delta _{a_0}\circ \delta _{a_1}\circ \cdots \circ \delta _{a_r}$. \hfill $\square $ 

\vskip 0.1cm
\proclaim{\bf Lemma 7.1} 
\item{$\bullet $} If $\Delta _1\in K\text{-}\bold{Mod}(P,Q), \Delta _2\in K\text{-}\bold{Mod}(Q,R)$ are differential operators of order $\le r$ and $\le s$ 
respectively, then $\Delta _2\circ \Delta _1:K\text{-}\bold{Mod}(P,R)$ is a differential operator of order $\le r+s$,
\item{$\bullet $} $\forall a\in A, \ P\in Ob\, (A\text{-}\bold{Mod})$ module multiplication (by $a$) \, $l_a:P\to P:p\mapsto ap$ is a differential operator of order $0$.  \hfill $\square $ 
\endproclaim 

All differential operators between $A$-modules form a category $A\text{-}\bold{Diff}$, such that 
$A\text{-}\bold{Mod}\hookrightarrow A\text{-}\bold{Diff}\hookrightarrow K\text{-}\bold{Mod}$,  
 and first two categories have the same objects. $A\text{-}\bold{Diff}$ is enriched in $(K\text{-}\bold{Mod}, \otimes _K)$ in a proper sense and 
enriched in two different ways in $(A\text{-}\bold{Mod}, \otimes _K)$ loosing composition property to be $A$-module map. Module 
multiplication for the first enrichment $A\text{-}\bold{Diff}$ in $(A\text{-}\bold{Mod}, \otimes _K)$ is given by 
$A\times A\text{-}\bold{Diff}(P,Q)\to A\text{-}\bold{Diff}(P,Q):(a,\Delta )\mapsto l_a\circ \Delta $, for the second enrichment by
$A\times A\text{-}\bold{Diff}(P,Q)\to A\text{-}\bold{Diff}(P,Q):(a,\Delta )\mapsto \Delta \circ l_a$. Denote $A\text{-}\bold{Diff}$ with 
left module multiplication in hom-sets $l_a\circ -$ by the same name $A\text{-}\bold{Diff}$ and with right multiplication in hom-sets 
$-\circ l_a$ by $A\text{-}\bold{Diff^+}$.

\vskip 0.1cm
\proclaim{\bf Proposition 7.1} 
\item{$\bullet $} $\forall P,Q\in Ob\, (A\text{-}\bold{Mod})$ \, $A\text{-}\bold{Diff}(P,Q)=\bigcup \limits _{s=0}^{\infty }\bold{Diff}_s(P,Q)$,
\, $A\text{-}\bold{Diff^+}(P,Q)=\bigcup \limits _{s=0}^{\infty }\bold{Diff}_s^+(P,Q)$ are filtered $A$-modules by submodules of differential 
operators of order $\le s, \, s=0,1,...$,
\item{$\bullet $} $\forall P\in Ob\, (A\text{-}\bold{Mod})$ \, $A\text{-}\bold{Diff}(P,P)$ is an associative $K$-algebra.  \hfill $\square $
\endproclaim

\proclaim{\bf Proposition 7.2} 
\item{$\bullet $} $\bold{Diff}_s(P,-),\, \bold{Diff}_s^+(-,P):A\text{-}\bold{Mod}\to A\text{-}\bold{Mod}$ are $A$-linear functors,
\item{$\bullet $} $\forall P\in Ob\, (A\text{-}\bold{Mod})$ functor $\bold{Diff}_s^+(-,P)$ is representable by object 
$\bold{Diff}_s^+(P):=\bold{Diff}_s^+(A,P)$, i.e. $\forall Q\in Ob\, (A\text{-}\bold{Mod})$ $A\text{-}\bold{Mod}(Q,\bold{Diff}_s^+(P))@>\sim >>\bold{Diff}_s^+(Q,P)$, 
\item{$\bullet $} $\forall P\in Ob\, (A\text{-}\bold{Mod})$ functor $\bold{Diff}_s(P,-)$ is representable by object 
$\bold{Jet}^s(P):=A\otimes _KP\, \text{\rm mod} \, \mu ^{s+1}$, where $\mu ^{s+1}$ is a submodule of $A\otimes _KP$ generated by elements 
$\delta ^{a_0}\circ \cdots \circ \delta ^{a_{s+1}}(a\otimes p)$ {\rm [$\delta ^{b}(a\otimes p):=ab\otimes p-a\otimes bp$]}, i.e. $\forall Q\in Ob\, (A\text{-}\bold{Mod})$ $A\text{-}\bold{Mod}(\bold{Jet}^s(P),Q)@>\sim >>\bold{Diff}^s(P,Q)$,
\item{$\bullet $} inclusion $A\text{-}\bold{Mod}\hookrightarrow A\text{-}\bold{Diff}^+$ is an (enriched) left adjoint with counit
$ev:\bold{Diff}^+(P)\to P:\Delta \mapsto \Delta (1)$, i.e. $\forall \Delta \in \bold{Diff}^+(Q,P)$ 
$\exists ! f_{\Delta }\in A\text{-}\bold{Mod}(Q,\bold{Diff}^+(P))$ such that 
$$\xymatrix{\bold{Diff}^+(P) \ar[r]^-{ev} & P\\
Q \ar@{-->}[u]^-{f_{\Delta }} \ar[ru]_-{\Delta } & }$$
and this correspondence is $A$-linear, \ $f_{\Delta }:q\mapsto (a\mapsto \Delta (aq))$,
\item{$\bullet $} inclusion $A\text{-}\bold{Mod}\hookrightarrow A\text{-}\bold{Diff}$ is an (enriched) right adjoint with unit 
$j^{\infty }:P\to \bold{Jet}^{\infty }(P):p\mapsto 1\otimes p\, \text{\rm mod}\, \mu ^{\infty }$ {\rm [$\mu ^{\infty }:=\bigcap \limits _{s=0}^{\infty }\mu ^s$]}, i.e. 
$\forall \Delta \in \bold{Diff}(P,Q)$ 
$\exists ! f^{\Delta }\in A\text{-}\bold{Mod}(\bold{Jet}^{\infty }(P),Q)$ such that 
$$\xymatrix{P \ar[r]^-{j^{\infty }} \ar[dr]_-{\Delta } & \bold{Jet}^{\infty }(P) \ar@{-->}[d]^-{f^{\Delta }}\\
 & Q  }$$
and this correspondence is $A$-linear, \ $f^{\Delta }:(a\otimes p)\, \text{\rm mod}\, \mu^{\infty } \mapsto a\Delta (p)$,
\item{$\bullet $} subcategory $A\text{-}\bold{Mod}$ is reflective and coreflective in $A\text{-}\bold{Diff}$ (enriched in $K\text{-}\bold{Mod}$). \hfill $\square $
\endproclaim

\vskip 0.15cm
$\forall s\in \Bbb N$ introduce two full subcategories of $A\text{-}\bold{Mod}$: 
\item{$\bullet $} $A\text{-}\bold{Mod}\text{-}\bold{Diff}_s$, consisting of all $A$-modules of type $\bold{Diff}_s(P,A), \ P\in Ob\, (A\text{-}\bold{Mod})$, 
and \linebreak 
$A$-module $A$,
\item{$\bullet $} $A\text{-}\bold{Mod}\text{-}\bold{Jet}^s$, consisting of all $A$-modules of type $\bold{Jet}^s(P), \ P\in Ob\, (A\text{-}\bold{Mod})$, 
and \linebreak
$A$-module $A$.

\vskip 0.3cm
\proclaim{\bf Proposition 7.3 (Vinogradov Duality)} \newline
There is a duality  
$\xymatrix{A\text{-}\bold{Mod}\text{-}\bold{Diff}_s^{op} \ar@/^/[r]^-{\sim } & A\text{-}\bold{Mod}\text{-}\bold{Jet}^s \ar@/^/[l]^-{\sim }}$, $s\in \Bbb N$, 
concrete over $A\text{-}\bold{Mod}$, namely, 
$\bold{Diff}_s(P,A)\simeq A\text{-}\bold{Mod}(\bold{Jet}^s(P),A)$, \ \ 
$\bold{Jet}^s(P)\simeq A\text{-}\bold{Mod}(\bold{Diff}^s(P,A),A)$. $A$ is a schizophrenic object. \hfill $\square $ 
\endproclaim 

\vskip 0.1cm
{\bf Remarks.}
\vskip 0.15cm
\item{$\bullet $} The above proposition states a formal analogue of duality between differential operators and jets over a fixed manifold $X$. 
Geometric modules of sections of vector bundles over $X$ correspond to modules $P$ over $C^{\infty }(X)$ with property 
$\bigcap \limits _{x\in X}\mu _xP=0$, where $\mu _x$ is a maximal ideal at point $x\in X$. Functors $\bold{Diff}_s(-,A)$ and $\bold{Jet}^s(-)$ preserve module property to be geometric \cite{V-K-L}.
\item{$\bullet $} This duality is an alternative (algebraic) way to introduce jet-bundles in Geometry (instead of classical approach
due to Grothendieck and Ehresmann as equivalence classes of maps tangent of order $s$ at a point). When $A=C^{\infty }(X)$ and $P$ is a geometric 
module realizable as a vector bundle $V(P)$ over $X$ then $\bold{Jet}^s(P)$ is realizable as $\bold{Jet}^s(V(P))$ over $X$ in 
classical sense \cite{V-K-L, Vin1, Vin2}.   \hfill $\square $

\head {\bf 8. Duality for differential equations}\endhead
\vskip 0.0cm
\proclaim{\bf {Proposition 8.1}} Let \, {\bf UAlg} be a category of universal algebras with representable forgetful functor.
Then every topological algebra ${\goth A}$ is a schizophrenic object (see \cite{P-Th}), and so, yields
a natural dual adjunction between {\bf UAlg} and {\bf Top}.
\endproclaim
\demo{Proof}
\item{$\bullet $} Initial topology on {\bf UAlg}(B,$\goth A$) gives initial lifting with respect to evaluation maps
$ev_{B,\, b}:\text{\bf UAlg}(B,\goth A)\to |\goth A|$, \ $b\in |B|$.
\item{$\bullet $} Algebra of continuous functions $\text{\bf Top}(X,\goth A)$ is initial with respect to evaluation maps
$ev_{X,\, x}:\text{\bf Top}(X,\goth A)\to |\goth A|$, $x\in |X|$ (which are obviously homomorphisms) since operations in \linebreak
$\text{\bf Top}(X,\goth A)$ are pointwise and each arrow $f\in \text{\bf Top}(X,\goth A)$ is completely determined by all its values
$ev_{X,\, x}(f)=|f|(x)$, \ $x\in |X|$. So that, if $g:|B|\to \text{\bf Top}(X,\goth A)$ is a {\bf Set}-map such that $\forall x\in |X|$
$ev_{X,\, x}\circ g$ is a homomorphism ($\omega _n(ev_{X,\, x}\circ g)b_1,...,(ev_{X,\, x}\circ g)b_n=
ev_{X,\, x}\circ g\omega _nb_1,...,b_n=ev_{X,\, x}\omega _ngb_1,...,gb_n$, where $\omega _n$ is an $n$-ary operation.
First equality holds because $ev_{X,\, x}\circ g$ is a homomorphism, second equality because $ev_{X,\, x}$ is a homomorphism), then $g$
is a homomorphism since two maps whose values coincide at each point coincide themselves.
\hfill $\square $
\enddemo

\proclaim{\bf Corollary} Take {\bf UAlg}$=k\text{-}\Lambda \text{-}\bold{Alg}$ (category of exterior differential algebras over a field $k$ ($\Bbb R$ or $\Bbb C$)
which presents differential equations). Take $\goth A=\Lambda (C^{\infty }({\Bbb R}^n))$ or $\Lambda (C^{\omega }({\Bbb C}^n))$
(which presents a parameter space) with a topology not weaker than $jet^{\infty }$. Then there exists a natural dual adjunction
$\xymatrix{k\text{-}\Lambda \text{-}\bold{Alg}^{op} \ar@/_/[r]^-{\perp } & {\text{\bf Top}} \ar@/_/[l]}$ (between differential equations and their solution spaces). \vskip 0cm \hfill $\square $
\endproclaim
\vskip -0.1cm
{\bf Remark.} If we regard full category $k\text{-}\Lambda \text{-}\bold{Alg}$ whose forgetful functor is representable we will get a lot of extra 'points'
which do not have geometric sense. Only graded maps of degree $0$ to $\goth A$
have geometric sense (they present integral manifolds of dimension not bigger than $n$). In this case representation of exterior
differential algebras when it exists will be not via their solution spaces but via much bigger spaces. If we restrict $k\text{-}\Lambda \text{-}\bold{Alg}$ to only
graded morphisms of degree $0$ then forgetful functor is not representable. But the notion of 'schizophrenic object' still makes sense
and theorem for natural dual adjunction \cite{P-Th} still holds. So, there is a representation of exterior differential algebras via
their usual solution spaces. \hfill $\square $
\vskip 0.2cm
Denote concrete subcategories of $\bold{Top}$ dual to categories $k\text{-}\bold{Alg}$ (algebras over $k$) and 
$k\text{-}\Lambda \text{-}\bold{Alg}$ (exterior differential algebras over $k$ with graded degree $0$ morphisms) 
by ${\text{\tt alg-}}\bold{Sol}$ and ${\text{\tt diff-}}\bold{Sol}$ respectively, i.e., 
$k\text{-}\bold{Alg}^{op}\sim {\text{\tt alg-}}\bold{Sol}$, $k\text{-}\Lambda \text{-}\bold{Alg}^{op}\sim {\text{\tt diff-}}\bold{Sol}$. 
In particular, ${\text{\tt alg-}}\bold{Sol}$ contains all algebraic and all smooth $k$-manifolds 
($k=\Bbb R$ or $\Bbb C$), ${\text{\tt diff-}}\bold{Sol}$ contains all spaces like ${\text{\tt alg-}}\bold{Sol}(k^n,X)$ 
(with representing object $\goth A=\Lambda (C^{\infty }(k^n))$).

\vskip 0.3cm
\proclaim{\bf {Lemma 8.1} ({\rm rough structure of} $\text{\tt diff-}\bold{Sol}$)}
\item{$\bullet $} $Ob\, (\text{\tt diff-}\bold{Sol})$ are pairs $(X,\coprod\limits _{i=1}^n\Cal F_i)$ where $X\hskip -0.1cm:=k\text{-}\Lambda \text{-}\bold{Alg}(D,k)=k\text{-}\bold{Alg}(D,k)\in Ob\, (\text{\tt alg-}\bold{Sol})$, $\Cal F_i\subset \text{\tt alg-}\bold{Sol}(k^i,X),\ 1\le i\le n$ [$\Cal F_i$ are not arbitrary subspaces of $\text{\tt alg-}\bold{Sol}(k^i,X)$].
\item{$\bullet $} $Ar\, (\text{\tt diff-}\bold{Sol})$ are pairs $(f,\coprod\limits _{i=1}^n\text{\tt alg-}\bold{Sol}(k^i,f)):(X,\coprod\limits _{i=1}^n\Cal F_i)\to (X',\coprod\limits _{i=1}^n\Cal F'_i)$ where $f:X\to X'\in Ar\, (\text{\tt alg-}\bold{Sol})$, $\text{\tt alg-}\bold{Sol}(k^i,f):\Cal F_i\to \Cal F'_i, \ 1\le i\le n$. \hfill $\square $
\endproclaim

\vskip 0.1cm
\proclaim{\bf {Proposition 8.2}} There are following adjunctions
\item{$\bullet $} $\xymatrix{k\text{-}\bold{Alg} \ar@/^/[r]^-{\ \Lambda _k}_-{\bot } & k\text{-}\Lambda \text{-}\bold{Alg} \ar@/^/[l]^-{p_0}}$ where
$\Lambda _k$ is the {\bf free exterior differential algebra functor} (see 4.2.1), \, $p_0$ is the projection onto subalgebra of
$0$-elements,
\vskip 0.1cm
\item{$\bullet $} $\xymatrix{\text{\tt alg-}\bold{Sol} \ar@/^/[r]^-{\text{ \rm hom}(k^n,-)}_-{\top } & \text{\tt diff-}\bold{Sol} \ar@/^/[l]^-{b}}$
where \, $b$ \, is taking the base space,

such that
$$\xymatrix{ k\text{-}\Lambda \text{-}\bold{Alg}^{op} \ar@/^2ex/[rr]^-{F}_-{\sim } \ar@/^2ex/[dd]^-{p_0^{op}} && \text{\tt diff-}\bold{Sol} \ar@/^2ex/[ll]^-{F'}_-{\sim } \ar@/^2ex/[dd]^-{b}\\
 && \\
k\text{-}\bold{Alg}^{op} \ar@/^2ex/[rr]^-{G}_-{\sim } \ar@/^2ex/[uu]^-{\Lambda _k^{op}}_-{\ \vdash } && \text{\tt alg-}\bold{Sol} \ar@/^2ex/[ll]^-{G'}_-{\sim } \ar@/^2ex/[uu]^-{\text{hom}(k^n,-)}_-{\ \vdash }
}$$
\endproclaim
\demo{Proof}
\item{$\bullet $} $\xymatrix{k\text{-}\Lambda \text{-}\bold{Alg}(\Lambda _k(A),D) \ar[r]^-{\sim }  & k\text{-}\bold{Alg}(A,p_0(D))\\
\rho  \ar@{{|}->}[r]^-{\sim } \ar[u]^-{\in } & \rho _0 \ar[u]_-{\in }
}$ \hskip 1cm(natural in $A$ and $D$)\newline
where  $\rho _0$ is the $0$-component of graded degree $0$ homomorphism $\rho =\bigoplus\limits_{i\ge 0}\rho _i$.
\item{$\bullet $} $\xymatrix{\text{\tt diff-}\bold{Sol}(S,\text{\rm hom}(k^n,X)) \ar[r]^-{\sim } & \text{\tt alg-}\bold{Sol}(b(S),X) \\
f \ar@{{|}->}[r]^-{\sim } \ar[u]^-{\in } & f \ar[u]_-{\in }
}$ \hskip 1cm(natural in $S$ and $X$)\newline
where:  $S$ is a pair $(b(S),\coprod\limits _{i=1}^n\Cal F_i), \ \Cal F_i\subset \text{hom}(k^i,b(S)), \ 1\le i\le n$, \
right $f:b(S)\to X$ is a usual map, and left $f:=(f,\coprod\limits _{i=1}^n\text{hom}(k^i,f)):(b(S),\coprod\limits _{i=1}^n\Cal F_i)\to (X,\coprod\limits _{i=1}^n\text{hom}(k^i,X))$.
\vskip 0.2cm
The above square of adjunctions is immediate.
\hfill $\square $
\enddemo

\vskip 0.2cm
\subhead 8.1. Cartan involution \endsubhead

\vskip 0.1cm
For systems in Cartan involution a (single) solution can be calculated recursively beginning from smallest $0$ dimension. By Cartan's theorem \cite{BC3G, Car1, Fin, Vas}
every system can be made into such a form by sufficient number of differential prolongations \cite{BC3G, Car1, Fin, Vas}. There is a 
cohomological criterion for systems to be in the involution. 

\vskip 0.1cm
{\bf Definition 8.1.1.} Let $\Cal A\in Ob\, (k\text{-}\Lambda \text{-}\bold{Alg})$, $\goth A_n$ be $\Lambda _{\Bbb R}(C^{\infty }(\Bbb R^n))$  or  $\Lambda _{\Bbb C}(C^{\omega }(\Bbb C^n))$, \, $n\ge 0$. 
\item{$\bullet $} Any (differential homomorphism of degree $0$) $\rho :\Cal A\to \goth A_n$ is called an {\bf integral manifold} of $\Cal A$ (of dimension not bigger than $n$).
\item{$\bullet $} $deg\, (\rho :\Cal A\to \goth A_n)=m, \ 0\le m\le n$, \ iff $\rho $ can be factored through a $\gamma :\Cal A\to \goth A_m$, i.e.,
$\vcenter{\xymatrix{\Cal A \ar@{..>}[r]^-{\gamma } \ar[dr]_-{\rho } & \goth A_m \ar@{..>}[d] \\
 & \goth A_n}}$ \, and $m$ is the smallest such number.
\item{$\bullet $} $deg\, (\Cal A)=n$ \, iff maximal degree of integral manifolds of \, $\Cal A$ \, is $n$.
\item{$\bullet $} $\Cal A$, $deg(\Cal A)=n$, is in {\bf Cartan involution} iff for each $m$-dimensional integral manifold $\rho :\Cal A\to \goth A_m$, $m<n$, 
there exists an $(m+1)$-dimensional integral manifold $\beta :\Cal A\to \goth A_{m+1}$ which {\bf contains} $\rho $, i.e., \hskip 0.1cm
$\vcenter{\xymatrix{\Cal A \ar@{..>}[r]^-{\exists \, \beta } \ar[dr]_-{\forall \, \rho } & \goth A_{m+1} \ar@{..>}[d] \\
 & \goth A_{m}}}$  \hfill $\square $ 

\vskip 0.25cm
{\bf Remarks.}
\item{$\bullet $} $\goth A_0$ is just $k$ ($\Bbb R$ or $\Bbb C$) with trivial differential. $\rho :\Cal A\to \goth A_0$ corresponds to a point $b\, (\rho ):b\, (\Cal A)\to k$.
Each point of $\Cal A$ is a $0$-dimensional integral manifold.
\item{$\bullet $} Original Cartan's definition was for classical algebras (factor-algebras of $\Lambda _{\Bbb R}(C^{\omega }(\Bbb R^N))$) 
and in terms of 'infinitesimal integral elements' (nondifferential homomorphisms of degree $0$ \, $f:\Cal A\to \Lambda _{k}(d\tau ^1,...,d\tau ^N)$) \cite{BC3G, Car1, Fin}. For that case both definitions coincide.
\item{$\bullet $} By a number of differential prolongations (adding new jet-variables with obvious relations) every classical 
system can be made into Cartan involution form (E. Cartan's theorem). 
\item{$\bullet $} Integration step (constructing a $1$ dimension bigger integral manifold) is done by appropriate 'Cauchy characteristics'.  \hfill $\square $

\vskip 0.3cm
\proclaim{\bf Proposition 8.1.1} Let $\Cal A$ be a factor-algebra of $\Lambda _{\Bbb R}(C^{\omega }(\Bbb R^N))$, $deg(\Cal A)=n$, 
corresponding to a system of differential equations $\vcenter{\xymatrix{\Cal E^q \ar@{>->}[r] \ar@{->>}[dr] & X\equiv b(F(\Cal A)) \ar@{->>}[d]^-{\pi } \ar@{==}[r] & \bold{Jet}^q(\Bbb R^{n+k}) \ar@{->>}[dl]^-{\pi } \\
 & \Bbb R^n & }}$, $dim(X)=N$.  
Then $\Cal A$ is in Cartan involution iff 
the following {\bf Spencer $\delta $-complex} is acyclic
\vskip 0.3cm
$\xymatrix{0 \ar[r] & g^{(r)} \ar[r]^-{\delta } & g^{(r-1)}\otimes \Lambda ^1(\Bbb R^n) \ar[r]^-{\delta } & g^{(r-2)}\otimes \Lambda ^2(\Bbb R^n) \ar[r]^-{\delta } & \cdots  \\
\cdots \ar[r]^-{\delta } & g^{(r-n)}\otimes \Lambda ^n(\Bbb R^n) \ar[r] & 0 & & }$
\vskip 0.2cm
where $g^{(r)}:=T(\bold{Jet}^r(\Cal E^q))\bigcap V\pi ^{q+r}_{q+r-1}\hookrightarrow S_{q+r}(T_*\Bbb R^n)\otimes V\pi $ is $r$-th prolongation of symbol $g$,  
$\pi ^{q+r}_{q+r-1}:\bold{Jet}^{q+r}(\Bbb R^{n+k})\to \bold{Jet}^{q+r-1}(\Bbb R^{n+k})$ is a natural projection of jet-bundles,
$V$ is taking 'vertical' subbundle, $S_p$ is $p$-th symmetric power,\newline
$\delta (\alpha _1\cdots \alpha _{q+r-l}\otimes v\otimes \beta _1\wedge \cdots \wedge \beta _l):=\sum \limits _{i=1}^{q+r-l}\alpha _1\cdots \hat {\alpha _i}\cdots \alpha _{q+r-l}\otimes v\otimes \alpha _i\wedge \beta _1\wedge \cdots \wedge \beta _l$,
\ $v$ is a section of $V\pi $.
\endproclaim 
\demo{Proof} See \cite{A-V-L, Sei, Vin2, V-K-L, Ver}   \hfill $\square $
\enddemo 

\vskip 0.1cm
Original Cartan's involutivity test was in terms of certain dimensions of 'infinitesimal integral elements'. The above theorem is 
due to J.P. Serre \cite{A-V-L, La-Se}.

\head {\bf 9. Gelfand-Naimark 2-duality}\endhead

Gelfand-Naimark duality is extendable to 2-duality over homotopies, which implies that cohomology theory for either C*-algebras or
compact Hausdorff spaces is automatically cohomology theory for the dual.

Let $\xymatrix{\bold {C^{*}Alg^{op}} \ar@/^/[r]^{F}& \bold {CHTop}\ar@/^/[l]^{G}_{\perp}}$ be the usual Gelfand-Naimark duality between commutative $C^*$-algebras and compact Hausdorff spaces. Both categories are strict 2-categories with homotopy classes of homotopies as 2-cells (homotopy of $C^*$-algebras is a homotopy in $\bold {Top}$ each instance of which is a $C^*$-algebra homomorphism). The reasonable quastion is: can it be extended to a 2-duality? The answer is yes.

By definition 
$$\xymatrix{
{\bold {C^*Alg}}(A,B)\times |A| \ar[r]^-{ev} & |B| & {\bold {C^*Alg}}(B,\Bbb C)\times {\bold {C^*Alg}}(A,B)\ar[r]^-{c_{A,B,\Bbb C}} & {\bold {C^*Alg}}(A,\Bbb C)\\
|I|\times |A|\ar[u]^{f\times 1} \ar[ur]_{\bar f} && {\bold {C^*Alg}}(B,\Bbb C)\times |I| \ar[u]^{1\times f} \ar[ur]_{F(\bar f)}\\
}$$
So that, if $f:|I|\times |A|\to |B|$ is a homotopy in $\bold {C^*Alg}$, then its image in $\bold {CHTop}$ is $F(\bar f):|F(B)|\times |I|\to |F(A)|$ (where $|\ |$ denotes underlying set or map).

We need to prove that such extended $F$ preserves 2-categorical structure (for $G$ proof is symmetric).

\centerline{\bf Preserving homotopies}

\proclaim{\bf Lemma 9.1}If $B$ is locally compact then $\xymatrix{{\bold {Top}}(B,C)\times {\bold {Top}}(A,B)\ar[r]^-{c_{A,B,C}} & {\bold {Top}}(A,C)}$ is continuous {\rm (with} compact-open {\rm topology in all hom-sets)}.
\endproclaim
\demo\nofrills{Proof \ }is standard. Let $f=g\circ h=c_{A,B,C}(g,h)$. Take $U^K$ be a (subbase) nbhd of $f$. Sufficient to show that $\exists$ (subbase) nbhds $U_1^{K^1}\ni g,\  U_2^{K^2}\ni h$, s.t. $U_1^{K^1}\circ U_2^{K^2}=c_{A,B,C}(U_1^{K^1},U_2^{K^2})\subset  U^{K}$. Take $U_1=U,\ K_2=K,\ K_1$ be a compact nbhd of $h(K)$, s.t. $K_1\subset g^{-1}(U)$ ($K_1$ exists by local compactness of $B$), $U_2=\text{int}(K_1)$. \hfill $\square $
\enddemo
\proclaim{\bf Corollary} If $A$ is locally compact then $ev_{A,B}:\bold {Top}(A,B)\times |A|\to |B|$ is continuous.
\endproclaim 
\demo{Proof} Each space $A$ is homeomorphic to $\bold {Top}(1,A)$ (with compact-open topology), and $ev_{A,B}$ corresponds to $c_{1,A,B}$. \hfill $\square $
\enddemo
\proclaim{\bf Lemma 9.2}  $\bullet $ Initial topology on $|F(A)|=\bold {C^*Alg}(A,\Bbb C)$ w.r.t. evaluation maps $\forall a\in A$  $\xymatrix{{\bold {C^*Alg}}(A,\Bbb C)\times 1
\ar[r]_{1\times a} & {\bold {C^*Alg}}(A,\Bbb C)\times |A|\ar[r]_-{ev} & |\Bbb C|}$ is point-open.

\hskip 1.9cm$\bullet$ Initial topology on $|G(X)|=\bold {CHTop}(X,\Bbb C)$ w.r.t. evaluation maps $\forall x\in X$  $\xymatrix{{\bold {CHTop}}(X,\Bbb C)\times 1
\ar[r]_{1\times x} & {\bold {CHTop}}(X,\Bbb C)\times |X|\ar[r]_-{ev} & |\Bbb C|}$ is compact-open.
\endproclaim
\demo{Proof} See \cite{P-Th}, \cite{Joh}, \cite{Eng}. \hfill $\square$
\enddemo

\proclaim{\bf Lemma 9.3} If ${\Cal A, \Cal B\subset \bold {LCTop}}$ are naturally dual subcategories of locally compact spaces (let $D$ be a dualizing object) then if $\Cal A(X,D)$ has compact-open topology (as initilal topology w.r.t. evaluation maps) then initial topology of $|X|\cong \Cal B(\Cal A(X,D),D)$ is compact-open as well.
\endproclaim
\demo{Proof} Evaluation map $ev:\Cal A(X,D)\times |X|\to |D|$ is continuous (since $X$ is locally compact and $\Cal A(X,D)$ has compact-open topology). It implies that initial (point-open) topology on $|X|\cong \Cal B(\Cal A(X,D),D)$ is actually compact-open [by assumption, topology of $|X|$ is initial w.r.t. all maps ${'f'}:|X|@>\sim >>1\times |X|@>f\times 1>>\Cal A(X,D)\times |X|@>ev>>|D|$. It means that topology on $|X|\cong \Cal B(\Cal A(X,D),D)$ is point-open since subbase open sets in point-open and initial topologies are the same $U^{'f'}:=\{ x\in |X|\bigm | {'f'}(x)\in \underset {open}\to U\subset D\} ={'f'}^{-1}(U)$].

We need to show that $\{x\in |X|\bigm | \forall f\in \underset {compact}\to K\subset \Cal A(X,D).\hphantom{ol} 'f'(x)\in \underset {open}\to U\subset D\}=\underset {f\in K}\to \bigcap {{'f'}^{-1}(U)}$ is open in point-open topology on $|X|$.

Take $x\in \underset {f\in K}\to \bigcap {'f'}^{-1}(U)$, then $ev\, (K,x)\subset U$. By continuity of $ev$, $\forall y\in K.\hphantom{o} \exists \underset {open}\to {V_y}\ni ~ y.\mathbreak
\hphantom{o} \exists \underset {open}\to {W_y}\ni x$, s.t. $ev\, (V_y,W_y)\subset U$. $\underset {open \atop covering}\to {\underset {y \in K}\to \bigcup V_y}\supset K$, so, by compactness, $\underset {j=1,\dots,n}\to \bigcup V_{y_j}\supset K$. Therefore, \,$ev\, \bigl (V_{y_j},\underset {j=1,\dots,n}\to \bigcap W_{y_j}\bigr )\subset U$, \,$ev\, \bigl (\underset {j=1,\dots,n}\to \bigcup V_{y_j},\underset {j=1,\dots,n}\to \bigcap W_{y_j}\bigr )\subset U$, \, $ev\, \bigl (K,\underset {j=1,\dots,n}\to \bigcap W_{y_j}\bigr )\subset U$, i.e., $x$ is internal. \hfill $\square$
\enddemo

\proclaim{\bf Corollary} Gelfand-Naimark duality preserves homotopies.
\endproclaim
\demo{Proof} $|A|=\bold {CHTop}(X,\Bbb C)$ has compact-open topology. $|X|=\bold {C^*Alg}(A,\Bbb C)$ has point-open topology, so, by {\bf {Lemma 3}} compact-open topology.

Multiplication $c_{A,B,\Bbb C}$ is continuous (since all hom-sets have compact-open topology). Therefore, $F(\bar f)$ is continuous.

[In inverse direction $G:{\bold {CHTop}}\to {\bold {C^*Alg}}$ there are no problem because $\bold {CHTop}(X,\Bbb C)$ has compact-open
topology. See also \cite{Loo}]. \hfill $\square$
\enddemo

\vskip 0.2cm
\centerline{\bf Preserving homotopy relation between homotopies}
\vskip 0.3cm

Let $\Bar {\Bar f}:|I|\times |I|\times |A|\to |B|$ be continuous, s.t. $\Bar {\Bar f}(0,t,a)={\Bar f}_0(t,a)$, \, $\Bar {\Bar f}(1,t,a)={\Bar f}_1(t,a)$.

$$\xymatrix{{\bold {C^*Alg}}(A,B)\times |A|\ar[r]^-{ev} & |B| & {\bold {C^*Alg}}(B,\Bbb C)\times {\bold {C^*Alg}}(A,B)\ar[r]^-{c_{A,B,\Bbb C}} & {\bold {C^*Alg}}(A,\Bbb C)\\
|I|\times |I|\times |A| \ar[u]^{\Bar {\Bar f}^T\times 1_{|A|}} \ar[ur]^{\Bar {\Bar f}} & & {\bold {C^*Alg}}(B,\Bbb C)\times |I|\times |I| \ar[u]^{1\times {\Bar {\Bar f}}^T} \ar[ur]^{F({\Bar {\Bar f}})}\\
1\times |I|\times |A| \ar@/^/[u]^{0\times 1_{|I|\times |A|}} \ar@/_/[u]_{1\times 1_{|I|\times |A|}} & \\
|I|\times |A| \ar[u]_{\sim}^{<!,1_{|I|}>\times 1_{|A|}} \ar@/_1pc/[uuur]^(0.7){{\bar f}_0} \ar@/_1.3pc/@<-1.2ex>[uuur]_(0.7){{\bar f}_1} & & {\bold {C^*Alg}}(B,\Bbb C)\times |I| \ar@/^/[uu]^{1\times ((0\times 1_{|I|})\circ <!,1_{|I|}>)} \ar@/_/[uu]_{1\times ((1\times 1_{|I|})\circ <!,1_{|I|}>)} \ar@/_3.5pc/[uuur]^(0.7){F({\bar f}_0)} \ar@/_3.7pc/@<-1.2ex>[uuur]_(0.7){F({\bar f}_1)}\\
}$$
So, $F({\Bar {\Bar f}})$ is a homotopy from $F({{\Bar f}_0})$ to $F({{\Bar f}_1})$. $F({\Bar {\Bar f}})$ is continuous since $c_{A,B,\Bbb C}$ is continuous in compact-open topology. ${\bold {C^*Alg}}(B,\Bbb C)$ has compact-open topology by {\bf Lemma 4.1.3}.
\vskip 5mm

\centerline{\bf Preserving unit 2-cells $i_f$}
\vskip -5mm
$$\xymatrix{{\bold {C^*Alg}}(A,B)\times |A|\ar[r]^-{ev} & |B| & {\bold {C^*Alg}}(B,\Bbb C)\times {\bold {C^*Alg}}(A,B)\ar[r]^-{c_{A,B,\Bbb C}} & {\bold {C^*Alg}}(A,\Bbb C)\\
1\times |A| \ar[u]^{f'\times 1} \ar[r]^-{\sim}_-{p_2} & |A| \ar[u]_{f} & {\bold {C^*Alg}}(B,\Bbb C)\times 1 \ar[u]^{1\times f'} \ar[r]^(0.5){\sim }_(0.5){p_1} \ar[ur]^{F(f\circ p_2)} & {\bold {C^*Alg}}(B,\Bbb C) \ar[u]_{-\circ f}\\
|I|\times |A| \ar[u]^{!\times 1} \ar[ur]_{p_2} && {\bold {C^*Alg}}(B,\Bbb C)\times |I| \ar[u]^{1\times !} \ar[ur]_{p_1} \ar@/_0.8pc/[uur]^(0.28){F(i_f)}\\
}$$

So, if $i_f=f\circ p_2\circ (!\times 1_{|A|})=f\circ p_2$, then $F(i_f)=F(f)\circ p_1=i_{F(f)}$.

\vskip 0.2cm
\centerline{\bf Preserving composites $i_g*{\bar f}:|I|\times |A|@>{\bar f}>>|B|@>g>>|C|$}
\centerline{\bf and ${\bar f*i_h:|I|\times |A'|@>{1\times h}>>|I|\times |A|}@>{\bar f}>>|B|$}

$$\xymatrix{{\bold {C^*Alg}}(A,C)\times |A|\ar[r]^-{ev} & |C| &{\bold {C^*Alg}}(C,\Bbb C)\times |I|\ar[dr]^{F(g\circ {\bar f})} \ar[d]|{1\times ({\bold {C^*Alg}}(A,g)\circ f)}   \ar@{->}@/_3.2pc/@<-8ex>[ddd]|{{\bold {C^*Alg}}(g,\Bbb C)\times 1} &\\
{\bold {C^*Alg}}(A,B)\times |A|\ar[u]^{(g\circ -)\times 1} \ar[r]^-{ev} & |B|\ar[u]_{g}& {\bold {C^*Alg}}(C,\Bbb C)\times {\bold {C^*Alg}}(A,C) \ar[r]_-{c_{A,C,\Bbb C}} & {\bold {C^*Alg}}(A,\Bbb C)\\
|I|\times |A| \ar[u]^{f\times 1} \ar[ur]_{\bar f}& &{\bold {C^*Alg}}(B,\Bbb C)\times {\bold {C^*Alg}}(A,B) \ar[r]^-{c{A,B,\Bbb C}} &{\bold {C^*Alg}}(A,\Bbb C) \ar[u]^{1}_{\sim }\\
&&{\bold {C^*Alg}}(B,\Bbb C)\times |I| \ar[u]|{1\times f} \ar[ur]_{F({\bar f})}&\\
}$$
$g\circ {\bar f}$ is a homotopy corresponding to ${\bold {C^*Alg}}(A,g)\circ f$. Outer perimeter of the right diagram commutes because of definition of $F({\bar f})$, $F(g\circ {\bar f})$ and associativity low [if $(s,t)\in {\bold {C^*Alg}}(C,\Bbb C)\times |I|$ then $s\circ (g\circ f(t))=(s\circ g)\circ f(t)$]. So, $F(g\circ {\bar f})=F({\bar f})\circ (F(g)\times 1_{|I|})$, i.e., $F(i_g*{\bar f})=F({\bar f})*i_{F(g)}$.

$$\xymatrix{ {\bold {C^*Alg}}(A,B)\times |A| \ar[r]^-{ev} & |B| & {\bold {C^*Alg}}(B,\Bbb C)\times {\bold {C^*Alg}}(A,B) \ar[r]^-{c_{A,B,\Bbb C}} & {\bold {C^*Alg}}(A,\Bbb C) \ar[dd]|{{\bold {C^*Alg}}(h,\Bbb C)}\\
|I|\times |A| \ar[u]^{f\times 1} \ar[ur]^{\bar f} & & {\bold {C^*Alg}}(B,\Bbb C)\times |I| \ar[u]^{1\times f} \ar[ur]^{F({\bar f})} \ar[d]_{1\times ({\bar f}\circ (1\times h))^T} \ar[dr]^{F({\bar f}\circ (1\times h))}&\\
|I|\times |A'| \ar[u]^{1\times h} \ar[uur]_{{\bar f}\circ (1\times h)} \ar[d]_{({\bar f}\circ (1\times h))^T\times 1} \ar[dr]|{{\bar f}\circ (1\times h)= \atop {ev\circ (f\times 1)\circ (1\times h)= \atop ev\circ (f\times h)}}& &{\bold {C^*Alg}}(B,\Bbb C)\times {\bold {C^*Alg}}(A',B) \ar[r]_-{c_{A',B,\Bbb C}} & {\bold {C^*Alg}}(A',\Bbb C)\\
{\bold {C^*Alg}}(A',B)\times |A'| \ar[r]_-{ev} & |B| \ar[uuu]_(0.4){\sim }^(0.4){1}\\
}$$

Right internal triangle of the right diagram commutes since if $(g,t)\in {\bold {C^*Alg}}(B,\Bbb C)\times |I|$ then ${\bold {C^*Alg}}(h,\Bbb C)\circ c_{A,B,\Bbb C}\circ (1\times f)(g,t)=(g\circ f(t))\circ h=g\circ (f(t)\circ h)=c_{A',B,\Bbb C}(g,f(t)\circ h)=c_{A',B,\Bbb C}(g,({\bar f}\circ (1\times h))^T(t))=c_{A',B,\Bbb C}\circ (1\times ({\bar f}\circ (1\times h))^T)(g,t)$. So, $F({\bar f}*i_h)=F({\bar f}\circ (1\times h))=F(h)\circ F({\bar f})=i_{F(h)}*F({\bar f})$.

\vskip 0.3cm
\centerline{\bf Preserving vertical composites}
\vskip 0.2cm

We need to show if ${\bar f}:{\bar f}\circ i_0\simeq {\bar f} \circ i_1$ and ${\bar g}:{\bar g}\circ i_0\simeq {\bar g} \circ i_1$ are homotopies in ${\bold {C^*Alg}}$ s.t. ${\bar f}\circ i_1={\bar g} \circ i_0$ then $F({\bar g}\odot {\bar f})=F({\bar g})\odot F(\bar f)$.

By definition, vertical composite ${\bar g}\odot {\bar f}$ is

$$\xymatrix@C=3pc{|A|\times |[0,{1\over 2}]| \ar@{^{(}->}@<-2ex>[dr]_{1\times i} & |A|\times |I| \ar[l]_-{1\times \alpha }^-{\sim } \ar[dr]^{\bar f}&\\
& |A|\times |I| \ar@{-->}[r]^{\exists !\, {\bar g}\odot {\bar f}} & |B|\\
|A|\times |[{1\over 2},1]| \ar@{^{(}->}@<1.5ex>[ur]^{1\times j} & |A|\times |I| \ar[l]_-{1\times \beta }^-{\sim } \ar[ur]_{\bar g}&\\
}$$

\hbox{
$$\xymatrix@C=0.5pc{
{\bold {C^*Alg}}(A,B)\times |A|\ar[rr]^-{ev} & & |B|\\
|I|\times |A| \ar[u]|{({\bar g}\odot {\bar f})^T\times 1} \ar[urr]_{{\bar g}\odot {\bar f}} & &\\
|[{1\over 2},1]|\times |A| \ar@{^{(}->}[u]^{j\times 1} & |[0,{1\over 2}]|\times |A| \ar@{^{(}->}@<-1ex>[ul]_{i\times 1}&\\
|I|\times |A| \ar@/^2pc/@<3ex>[uuu]|{f\times 1} \ar[ur]_{\alpha \times 1} \ar[uuurr]_(0.7){\bar f}&&\\
|I|\times |A| \ar@/^2pc/[uu]|(0.25){\beta \times 1} \ar@/^3pc/@<4ex>[uuuu]|{g\times 1} \ar@/_4pc/[uuuurr]_{\bar g}&&\\
}$$

$$\xymatrix@C=1.5pc{{\bold {C^*Alg}}(B,\Bbb C)\times {\bold {C^*Alg}}(A,B) \ar[r]^-{c_{A,B,\Bbb C}}& {\bold {C^*Alg}}(A,\Bbb C)\\
&\\
{\bold {C^*Alg}}(B,\Bbb C)\times |I| \ar[uu]^{1\times ({\bar g}\odot {\bar f})^T} \ar[uur]^{F({\bar g}\odot {\bar f})}& \\
&\\
{\bold {C^*Alg}}(B,\Bbb C)\times |I| \ar[uu]^{1\times (j\circ \beta)} \ar[uuuur]_{F({\bar g})} & {\bold {C^*Alg}}(B,\Bbb C)\times |I| \ar[uul]_(0.4){1\times (i\circ \alpha )} \ar[uuuu]_{F(\bar f)}\\
}$$
}

By uniqueness $f\equiv {\bar f}^T=({\bar g}\odot {\bar f})^T\circ i\circ \alpha$, \  $g\equiv {\bar g}^T=({\bar g}\odot {\bar f})^T\circ j\circ \beta$.\pagebreak

So, $\cases F({\bar g}\odot {\bar f})\circ (1\times (i\circ \alpha ))=F({\bar f})\\ F({\bar g}\odot {\bar f})\circ (1\times (j\circ \beta ))=F({\bar g})\endcases$. \ It means $F({\bar g}\odot {\bar f})=F({\bar g})\odot F({\bar f})$.

\centerline{\bf Preserving horisontal composites $\xymatrix{A \ar@/^/[r]^{f_0}_(0.4){{\bar f}\, \Downarrow } \ar@/_/[r]_{f_1}& B\ar@/^/[r]^{g_0}_(0.4){{\bar g}\, \Downarrow } \ar@/_/[r]_{g_1}&C}$}

${\bar g}*{\bar f}:=({\bar g}*i_{f_1})\odot (i_{g_0}*{\bar f})\simeq (i_{g_1}*{\bar f})\odot ({\bar g}*i_{f_0})$ (homotopic homotopies). \newline
$F({\bar g}*{\bar f})=F({\bar g}*i_{f_1})\odot F(i_{g_0}*{\bar f})=(i_{F(f_1)}*F({\bar g}))\odot (F({\bar f})*i_{F(g_0)})\simeq F({\bar f})*F({\bar g})$. \

\vskip 0.2cm
Proposition 9.1 completes the proof of Gelfand-Naimark 2-duality 
$\xymatrix{\bold {C^{*}Alg^{op}} \ar@/^/[r]^{F}& \bold {CHTop}\ar@/^/[l]^{G}_{\perp}}$.

\vskip 0.0cm
\proclaim{\bf Proposition 9.1} If $\xymatrix{\bold{C} \ar@/^/[r]^-{F} & \bold{D} \ar@/^/[l]^-{G} }$ are two strict $n$-categories
and two strict $n$-functors in the opposite directions such that the restriction
$\xymatrix{\bold{C}^{\le 1} \ar@/^/[r]^-{F^{\le 1}}_-{\perp } & \bold{D}^{\le 1} \ar@/^/[l]^-{G^{\le 1}} }$
is an adjunction with unit $\eta :1_{\bold{C}^{\le 1}}\to G^{\le 1}F^{\le 1}$ and counit 
$\varepsilon :F^{\le 1}G^{\le 1}\to 1_{\bold{D}^{\le 1}}$ which are still natural transformations for the extension 
(i.e. $\eta :1_{\bold{C}}\to GF$ and $\varepsilon :FG\to 1_{\bold{D}}$ are natural transformations) then the extended situation 
$\xymatrix{\bold{C} \ar@/^/[r]^-{F}_-{\perp } & \bold{D} \ar@/^/[l]^-{G} }$ is a strict adjunction.
\endproclaim 
\demo{Proof} A strict adjunction is completely determined by its 'unit-counit' (proposition 1.5.3). $\eta :1_{\bold{C}}\to GF$ and
$\varepsilon :FG\to 1_{\bold{D}}$ are natural transformations and satisfy triangle identities
$\varepsilon F\circ _1F\eta =1_F$ and $G\varepsilon \circ _1\eta G=1_G$ (because, e.g. $\varepsilon F=\varepsilon F^{\le 1}$, $1_{F}=1_{F^{\le 1}}$ (set-theoretically), etc.)  \vskip 0.0cm  \hfill   $\square $
\enddemo

\vskip 0.2cm
{\bf Corollary.} Any $1$-adjunction between a category of topological algebras and 
a subcategory of topological spaces is a $2$-adjunction if it can be extended functorially over 2-cells in the way that each instance
of the image of a homotopy is the image of this instance of the preimage-homotopy. 
\demo{Proof} Under given conditions unit and counit of 1-adjunction are automatically natural transformations for the extension. 
E.g., take unit $\eta $.  Naturality square $\vcenter{\xymatrix{A  \ar[r]^-{\eta _A}  \ar[d]_-{f^1} &   GFA  \ar[d]^-{GFf^1}  \\
B   \ar[r]_-{\eta _B}  &     GFB  }}$,
where $f^1:A\times I\to B$ is a homotopy, holds because each instance of it holds (since $\eta $ is a unit of 1-adjunction), i.e.
$\forall t\in I$ $\eta _B\circ f^1(-,t)=GF(f^1(-,t))\circ \eta _A$, it means $\eta _B\circ f^1=GF(f^1)\circ (\eta _A\times I)$, i.e.
$\eta _B*f^1=GF(f^1)*\eta _A$.            \hfill            $\square $
\enddemo

Gelfand-Naimark case is one of the above corollary.   End of proof of Gelfand-Naimark 2-duality.   \vskip 0.0cm  \hfill  $\square $

\vskip 0.2cm
{\bf Remark.} There are 'forgetful' functors $\bold {C^{*}Alg}\to 2\text{-}\bold{Set}$ and 
$\bold {CHTop}\to 2\text{-}\bold{Set}$ (where $2\text{-}\bold{Set}$ is the usual $\bold{Set}$ with just one iso-2-cell for each 
pair of maps with the same domain and codomain) but they are not faithful and forget too much in order $2\text{-}\bold{Set}$ 
could be an underlying category of Gelfand-Naimark 2-duality.     \hfill        $\square $

\vskip 0.3cm
\proclaim{\bf Proposition 9.2} $\bullet $ Gelfand-Naimark 2-duality is concrete over $2\text{-}\bold{Cat}$ ($2\text{-}\bold{Cat}$ is the 
usual 2-category of (small) categories, functors and natural transformations), i.e. $\exists $ (faithful) forgetful functors 
$U:\bold {C^{*}Alg}\to 2\text{-}\bold{Cat}$ and $V:\bold {CHTop}\to 2\text{-}\bold{Cat}$ such that 
$\vcenter{\xymatrix{\bold {C^{*}Alg^{op}}  \ar[r]^-{F} \ar[dr]_-{\bold {C^{*}Alg}(-,\Bbb C)\hskip 0.2cm} &   \bold {CHTop} \ar[d]^-{V} \\
    &          2\text{-}\bold{Cat}  }}$
and \hskip 0.15cm
$\vcenter{\xymatrix{\bold {CHTop^{op}} \ar[r]^-{G^{op}} \ar[dr]_-{\bold {CHTop}(-,\Bbb C)\hskip 0.2cm}  &    \bold {C^{*}Alg}  \ar[d]^-{U} \\
    &          2\text{-}\bold{Cat}  }}$ \hskip 0.1cm
where $U$ and $V$ are composites of inclusion and fundamental groupoid functors  
($U:\bold {C^{*}Alg}\hookrightarrow 2\text{-}\bold{Top}@>2\text{-}\bold{Top}(1,-)>>2\text{-}\bold{Cat}$ and  
$V:\bold {CHTop}\hookrightarrow 2\text{-}\bold{Top}@>2\text{-}\bold{Top}(1,-)>>2\text{-}\bold{Cat}$).  
\vskip 0.1cm\item{$\bullet $} This duality is natural, i.e. lifting of hom-functors $\bold {C^{*}Alg}(-,\Bbb C)$, $\bold {CHTop}(-,\Bbb C)$ along
$V$ and $U$ is initial.                           \hfill    $\square $
\endproclaim

\vskip 0.2cm
{\bf Remark.} 2-duality allows us to transfer (co)homology theories from one side to another.
Under a reasonable assumption that K-theory was determined in a universal way we could get {\bf M. Atiyah theorem} that 
{\it K-groups of $C^*$-algebras and compact Hausdorff spaces coincide}. The problem, however, is that K-groups were determined technically
(not universally). But, there is a theorem by J. Cuntz \cite{Weg} that K-theory is universally determied on a large subcategory of 
$C^*$-algebras.            \hfill        $\square $

\head {\bf Bibliography}\endhead

\vskip 0.2cm
\refstyle{A}
\widestnumber\key{AAAAA}

\ref\key A-H-S
\by J. Adamek, H. Herrlich, G.E. Strecker
\book Abstract and Concrete Categories. The Joy of Cats
\yr online edition, 2004
\endref

\ref\key A-V-L
\by D.V. Alekseevskiy, A.M. Vinogradov, V.V. Lychagin
\book Main Ideas and Concepts of Differential Geometry
\yr 1988
\publ Moscow
\lang Russian
\endref

\ref\key Bel
\by J.L. Bell
\book Toposes and Local Set Theories: An Introduction
\yr 1988
\publ Clarendon Press, Oxford
\endref

\ref\key Bi-Cr
\by R.L. Bishop, R.J. Crittenden
\book Geometry of Manifolds
\publ Academic Press, New York and London
\yr 1964
\endref

\ref\key Bor1
\by F. Borceux
\book Handbook of Categorical Algebra 1. Basic Category Theory
\yr 1994
\publ Cambridge University Press
\endref

\ref\key Bor2
\by F. Borceux
\book Handbook of Categorical Algebra 2. Categories and Structures
\yr 1994
\publ Cambridge University Press
\endref

\ref\key Bor3
\by F. Borceux
\book Handbook of Categorical Algebra 3. Categories of Sheaves
\yr 1994
\publ Cambridge University Press
\endref

\ref\key Bru
\by U. Bruzzo
\book Introduction to Algebraic Topology and Algebraic Geometry
\yr 2002
\publ International School for Advanced Studies, Trieste
\endref

\ref\key BC3G
\by R.L. Bryant, S.S. Chern, R.B. Gardner, H.L. Goldschmidt, P.A. Griffiths
\book Exterior Differential Systems
\yr 1991
\publ Springer-Verlag
\endref

\ref\key Car1
\by E. Cartan
\book Exterior Differential Systems and their Geometric Applications
\yr 1962
\publ Moscow State University
\lang translated into Russian
\endref

\ref\key Car2
\by E. Cartan
\book Moving Frame Method and Theory of Groups of Transformations 
\yr 1963
\publ Moscow State University
\lang translated into Russian
\endref

\ref\key C-C-L
\by S.S. Chern, W.H.Chen, K.S. Lam
\book Lectures on Differential Geometry
\yr 2000
\publ World Scientific
\endref

\ref\key Cl-D
\by D.M. Clark, B.A. Davey
\book Natural Dualities for the Working Algebraist
\yr 1997
\publ Cambridge University Press
\endref

\ref\key C-L
\by E. Cheng, A. Lauda
\book Higher-Dimensional Categories: an illustrated guide book
\publ University of Cambridge
\yr 2004
\endref 

\ref\key D-N-F
\by B.A. Dubrovin, S.P. Novikov, A.T. Fomenko
\book Modern Geometry
\publ Moscow
\yr 1979
\lang Russian
\endref

\ref\key Eng
\by R. Engelking
\book General Topology
\yr 1977
\endref

\ref\key ELOS
\by L.E. Evtushik, U.G. Lumiste, N.M. Ostianu, A.P. Shirokov
\book Differential Geometric Structures on Manifolds
\yr 1979
\publ Moscow
\lang Russian
\endref

\ref\key Fin
\by S.P. Finikov
\book The Method of Cartan's Exterior Forms
\yr 1948
\publ Moscow
\lang Russian
\endref

\ref\key Gar
\by R.B. Gardner
\book The Method of Equivalence and Its Applications
\yr 1989
\publ Society for Industrial and Applied Mathematics
\endref

\ref\key Gol
\by R. Goldblatt
\book Topoi: The Categorial Analysis of Logic
\yr 1984
\publ North-Holland
\endref

\ref\key Hof
\by D. Hofmann
\paper Natural Dualities 
\endref

\ref\key Jac
\by B. Jacobs
\book Categorical Logic and Type Theory
\publ Elsevier, North-Holland
\yr 2001
\endref

\ref\key Joh
\by P.T. Johnstone
\book Stone spaces
\publ Cambridge University Press
\yr 1982
\endref

\ref\key Kob
\by S. Kobayashi
\book Groups of Transformations in Differential Geometry
\publ Moscow
\lang translated into Russian
\yr 1986
\endref

\ref\key Koc
\by J. Kock
\paper Weak Identity Arrows in Higher Categories
\inbook IMRP International Mathematics Research Papers, Volume 2006, Article ID 69163
\pages 1-54
\yr 2006
\endref 

\ref\key Kra
\by I.S. Krasil'shchik
\paper Calculus over Commutative Algebras: a Concise User Guide
\inbook The Diffiety Institute Preprint Series
\yr 1996
\endref 

\ref\key La-Se
\by L.A. Lambe and W.M. Seiler
\book Differential Equations, Spencer Cohomology, and Computing Resolutions
\publ online article
\endref

\ref\key Lap
\by G.F. Laptev
\paper The Main Infinitesimal Structures of Higher Order on a Smooth Manifold
\inbook Proceedings of Geometric Seminar, vol.\, 1
\publaddr Moscow
\yr 1966
\lang Russian
\pages 139-190
\endref

\ref\key Lei
\by T. Leinster
\book Higher Operads, Higher Categories
\publ Cambridge University Press
\yr 2003
\endref

\ref\key Lic
\by A. Lichnerowicz
\book Theorie Globale des Connexions et des Groupes d'Holonomie
\publ Roma, Edizioni Cremonese
\yr 1955
\lang translated into Russian
\endref

\ref\key Loo
\by L.H. Loomis
\book An Introduction to Abstract Harmonic Analysis
\publ D. Van Nostrand Company
\yr 1953
\endref

\ref\key Mac
\by S. MacLane
\book Categories for the Working Mathematician
\publ Springer-Verlag
\yr 1971
\endref

\ref\key M-M
\by S. MacLane, I. Moerdijk
\book Sheaves in Geometry and Logic
\publ Springer-Verlag
\yr 1992
\endref

\ref\key M-T
\by I. Madsen, J. Tornehave
\book From Calculus to Cohomology
\publ Cambridge University Press
\yr 1996
\endref

\ref\key Man
\by U.I. Manin
\book Lectures in Algebraic Geometry. Part 1. Affine Schemes
\publ Moscow State University
\yr 1970
\lang Russian
\endref

\ref\key May
\by J.P. May
\book A Concise Course in Algebraic Topology
\publ The University of Chicago Press
\yr 1999
\endref

\ref\key Nes
\by J. Nestruev
\book Smooth manifolds and observables
\publ Moscow
\yr 2003
\lang Russian
\endref

\ref\key Olv
\by P.J. Olver
\book Equivalence, Invariants and Symmetry
\publ Cambridge University Press
\yr 1995
\endref

\ref\key P-Th
\by H.-E. Porst, W. Tholen
\paper Concrete dualities
\inbook Category Theory at Work\eds H. Herrlich. H.-E. Porst
\publ Heldermann Verlag 
\publaddr Berlin \yr 1991
\pages 111-136
\endref

\ref\key Pos
\by M.M. Postnikov
\book Lectures in Algebraic Topology. Foundations of Homotopy Theory
\publ Moscow
\yr 1984
\lang Russian
\endref

\ref\key Sat
\by H. Sato
\book Algebraic Topology: An Intuitive Approach
\publ American Mathematical Society
\yr 1996
\endref

\ref\key Sei
\by W.M. Seiler
\book Differential Equations, Spencer Cohomology, and Pommaret Bases
\publ Heidelberg University, online lecture notes
\endref

\ref\key Sch
\by J.T. Schwartz
\book Differential Geometry and Topology
\publ New York
\yr 1968
\endref

\ref\key Ste
\by S. Sternberg
\book Lectures on Differential Geometry
\publ Prentice Hall
\yr 1964
\endref

\ref\key Str1
\by T. Streicher
\book Introduction to Category Theory and Categorical Logic
\publ online lecture notes
\yr 2003
\endref

\ref\key Vas0
\by A.M. Vasiliev
\paper Differential Algebra as a Technics of Differential Geometry
\inbook Proceedings of Geometric Seminar, vol.\, 1
\publaddr Moscow
\yr 1966
\lang Russian
\endref

\ref\key Vas
\by A.M. Vasiliev
\book Theory of Differential-Geometric Structures
\yr 1987
\publ Moscow State University
\lang Russian
\endref

\ref\key Ver
\by A.M. Verbovetskii
\book Geometry of Finite Jets and Differential Equations: additional chapters
\publ Moscow Independent University
\yr 1999
\lang Russian
\endref 

\ref\key Vin1
\by A.M. Vinogradov
\paper The Logic Algebra for the Theory of Linear Differential Operators
\publ Dokl. Akad. Nauk SSSR (Soviet Math. Dokl.), vol. 13, No. 4
\yr 1972
\lang Russian
\endref

\ref\key Vin2
\by A.M. Vinogradov
\book Geometry of Nonlinear Differential Equations
\publ Moscow
\yr 1979
\lang Russian
\endref

\ref\key V-K-L
\by A.M. Vinogradov, I.C. Krasil'schik, V.V. Lychagin
\book Introduction to Geometry of Nonliniear Differential Equations
\publ Moscow
\yr 1986
\lang Russian
\endref

\ref\key Weg
\by N.E. Wegge-Olsen
\book K-Theory and $C^*$-Algebras, A Friendly Approach
\publ Oxford University Press
\yr 1993
\endref

\ref\key S-W
\by R. Sulanke, P. Wintgen
\book Differential Geometry and Fiber Bundles
\yr 1972
\publ Berlin
\lang translated into Russian
\endref

\enddocument